\documentclass{article}

\usepackage{amsmath, amsthm, amssymb, amsfonts, mathtools}
\usepackage{thmtools}
\usepackage{graphicx}
\usepackage{setspace}
\usepackage{geometry}
\usepackage{float}
\usepackage{hyperref}
\usepackage[utf8]{inputenc}
\usepackage[english]{babel}
\usepackage{framed}
\usepackage[all]{xy}
\usepackage{mathrsfs}
\usepackage{calligra}

\newtheorem{thm}{Theorem}[section]
\newtheorem{cor}[thm]{Corollary}
\newtheorem{lem}[thm]{Lemma}
\newtheorem{prop}[thm]{Proposition}

\newtheorem{rem}[thm]{Remark}
\newtheorem{conj}[thm]{Conjecture}
\newtheorem{hyp}[thm]{Hypothesis}

\numberwithin{equation}{subsection}

\begin{document}


\title{p-adic Hodge parameters in the crystalline representations of GSp$_4$}
		
\date{\today}
\author{Han Xiaozheng}
\maketitle
\abstract{This article gives a generalization of the work of \cite{Ding} in the context of $\mathrm{GSp}_4(\mathbb{Q}_p)$, where $p$ is an odd prime number. Let $\rho$ be a 4-dimensional generic non-critical crystalline representations of the absolute Galois group of $\mathbb{Q}_p$ of regular Hodge-Tate weights which is valued in $\mathrm{GSp}_4(E)$, where $E$ is a finite extension of $\mathbb{Q}_p$, we associate to $\rho$ an explicit locally analytic $E$-representation $\pi_\mathrm{min}(\rho)$ of $\mathrm{GSp}_4(\mathbb{Q}_p)$, which encodes enough information to determines $\rho$. Moreover, under certain settings, this construction follows the local-global compatibility.}
\tableofcontents
\newpage


\section{Introduction}
The locally analytic $p$-adic Langlands program for $\mathrm{GSp}_4(\mathbb{Q}_p)$ aims at building a correspondence between $p$-adic continuous representations valued in $\mathrm{GSp}_4$ of the abosulte Galois group $\mathrm{Gal}_{\mathbb{Q}_p}$ of $\mathbb{Q}_p$ and certain locally analytic representations of $\mathrm{GSp}_4(\mathbb{Q}_p)$. In particular, it is expected to match the parameters on both sides via the conjectural correspondence. Let $p$ be an odd prime number and $E/\mathbb{Q}_p$ be a finite extension. Let $\rho: \mathrm{Gal}_{\mathbb{Q}_p}\rightarrow \mathrm{GSp}_4(E)$ be a $\mathrm{GSp}_4$-valued crystalline representation. In this paper, if $\rho$ is generic, non-critical and have distinct Hodge-Tate weights, I construct a locally analytic representation $\pi_\mathrm{min}(\rho)$ of $\mathrm{GSp}_4(\mathbb{Q}_p)$ encoding all the information of the given $\mathrm{GSp}_4$-valued Galois representation $\rho$ and prove such construction following local-global compatibility in certain settings. This is a generalization of the work of Y. Ding in \cite{Ding}.

By Fontaine’s theory, $\rho$ is equivalent to the associated filtered $\varphi$-module with $\mathrm{GSp}_4$-structure, $D_{\mathrm{cris}}(\rho)$. (See \S 3.1 for the explicit definition of "with $\mathrm{GSp}_4$-structure".) More precisely, $D_{\mathrm{cris}}(\rho)$ is a $4$-dimensional symplectic space over $E$, equipped with a $\varphi$-operator $\varphi\in\mathrm{GSp}_4(D_{\mathrm{cris}}(\rho))$ and a Hodge filtration 
$$D_{\mathrm{cris}}(\rho)=\mathrm{Fil}^{-\infty}D_{\mathrm{cris}}(\rho)\supset\cdots\supset\mathrm{Fil}^{0}D_{\mathrm{cris}}(\rho)\supset\mathrm{Fil}^{1}D_{\mathrm{cris}}(\rho)\supset\cdots\supset\mathrm{Fil}^{\infty}D_{\mathrm{cris}}(\rho)=0$$
which is anisotropic, i.e., for a unique integer $h_0\in\mathbb Z$ and all $i\in\mathbb Z$, $\mathrm{Fil}^{i}D_{\mathrm{cris}}(\rho)=(\mathrm{Fil}^{h_0+1-i}D_{\mathrm{cris}}(\rho))^\perp$.

Assume $\boldsymbol{\alpha}=(\alpha_1,\alpha_2,\alpha_3,\alpha_4)$ are eigenvalues of the $\varphi$-operator, $\alpha_1\alpha_4=\alpha_2\alpha_3$. We assume the $\varphi$-action is generic (and we simply call such $\rho$ generic), which means $\alpha_i\alpha_j^{-1}\neq 1,p^{\pm1}$ for $1\leq i<j\leq 4$. Consider the smooth character of $T(\mathbb{Q}_p)$ $\phi:=\frac{\mathrm{unr}(\alpha_1)}{\mathrm{unr}(\alpha_3)}(p_1)\frac{\mathrm{unr}(\alpha_1)}{\mathrm{unr}(\alpha_2)}(p_2)\mathrm{unr}(\alpha_4)(p_3)$. (See \S2.3 for $p_1,p_2,p_3$.)

Then the smooth irreducible representation $\pi_{\mathrm{sm}}(\phi):=(\mathrm{Ind}^{\mathrm{GSp}_4(\mathbb{Q}_p)}_{\overline B(\mathbb{Q}_p)}\phi\vert p_1\vert_p^{-2}\vert p_2\vert^{-1})^\infty$ encodes the information of the $\varphi$-action of $D_{\mathrm{cris}}(\rho)$. Here $\overline{B}$ is the Borel subgroup consisting of lower triangular matrices.

Assume $\boldsymbol{h}=(h_1,h_2,h_3,h_4)$ are the Hodge-Tate weights of $\rho$, in other words, the Hodge filtration jumps at $-h_i$ for $1\leq i\leq 4$. We assume $\rho$ has distinct Hodge-Tate weights, i.e., $h_1>h_2>h_3>h_4$. Consider the algebraic character of $T(\mathbb{Q}_p)$ $\lambda:=p_1^{h_1-h_3-2}p_2^{h_1-h_2-1}p_3^{h_4}$. It is a dominant weight with respect to the Borel subgroup consisting of upper triangular matrices. Then the irreducible algebraic representation $L(\lambda)$ encodes the information of the Hodge-Tate weights. Then the locally algebraic representation $\pi_{\mathrm{alg}}(\phi,\boldsymbol{h}):=\pi_{\mathrm{sm}}(\phi)\otimes_EL(\lambda)$ encodes the information of both $\varphi$-action and Hodge-Tate weights of $D_{\mathrm{cris}}(\rho)$.

Because of Hodge-Tate weights are distinct, the underlying flag of the Hodge filtration is a complete anisotropic flag. Its relative position with $\varphi$-eigenlines is parameterized by $E$-points in $T\backslash\mathrm{GSp}_4/B$, which is called the Hodge parameter of $\rho$. We assume $\rho$ is non-critical, i.e., the Hodge filtration is in general position with $\varphi$-eigenlines. It can be represented by $2$ scalars as \eqref{standard form}.

In \S3, I study the $\mathrm{GSp}_4$-valued deformation theory of $\rho$ to detect the Hodge parameter. In \S3.2, by showing that the $\mathrm{GSp}_4$-valued deformation functor of $\rho$ is a subfunctor of the deformation functor of $\rho$ with forgetting the $\mathrm{GSp}_4$-structure, I give an analogue of \cite[\S2.3.1]{Ding} to the case of $\mathrm{GSp}_4$. In \S3.3, using Fontaine's functor $D_{\mathrm{pdR}}$, I prove that the linear relations of the $8$ tranguline subspaces of the tangent space of the $\mathrm{GSp}_4$-valued deformation functor of $\rho$ implies the Hodge parameter of $\rho$ (Thm. \ref{3.17}).

In \S4.1, I calculate a series of extension groups of locally analytic representations which contribute to the construction of $\pi_\mathrm{min}(\rho)$. In \S4.2, by combining the Galois side (\S3) and the automorphic side (\S4.1), I give an explicit construction of the locally analytic representation $\pi_\mathrm{min}(\rho)$ which encodes all information of $\rho$ (Thm. \ref{4.14} and Thm. \ref{4.18}). $\pi_\mathrm{min}(\rho)$ has $3$ layers in total. The socle of $\pi_\mathrm{min}(\rho)$ is isomorphic to $\pi_{\mathrm{alg}}(\phi,\boldsymbol{h})$. The cosocle of $\pi_\mathrm{min}(\rho)$ is isomorphic to $\pi_{\mathrm{alg}}(\phi,\boldsymbol{h})^{\oplus 2} $.

In \S5, I prove that $\pi_\mathrm{min}(\rho)$ satisfies local-global compatibility (Cor. \ref{5.19}) in the setting of \S5.2 by studying the tangent space of eigenvariety. 

At the end of the introduction, let me remark that the local part (i.e. the construction of $\pi_{\mathrm{min}}(\rho)$) is supposed to be directly generalized to the context of abstract reductive groups with a "good" embedding to the general linear groups, using the language of $p$-adic Hodge theory with $G$-structure, with no essential difficulty. (Rem. \ref{3.19}) I am working on this generalization and I believe that it will come soon.

\section{Review of structure of $\mathrm{GSp}_4$ and notations}

    \subsection{Some subgroups}
We denote $$J=\begin{bmatrix}&&&1\\&&1&\\&-1&&\\-1&&&\end{bmatrix}.$$
For ring F,$$\mathrm{GSp}_4(F)=\{A\in\mathrm{GL}_4(F)\vert A^TJA=sim(A)J\ for\ some\ sim(A)\in F^\times \}.$$
It's a fact that $\mathrm{sim}(A)^2=\mathrm{det}(A)$.
Its center is $Z(F)=F^\times I$, here $I$ is the identity matrix. We fix a Borel subgroup $B$ which consists of all upper triangular matrices in $G$. Its maximal torus is $$T(F)=\{\mathrm{diag}(a,b,cb^{-1},ca^{-1})\vert a,b,c\in F^\times\}.$$ Here $c$ is just the scalar $\mathrm{sim}(\mathrm{diag}(a,b,cb^{-1},ca^{-1}))$.

For $M\in\mathrm{GL}_2(F)$, we denote $M'=\begin{bmatrix}0&1\\1&0\end{bmatrix}(M^T)^{-1}\begin{bmatrix}0&1\\1&0\end{bmatrix}$.
We denote $$L_P(F)=\{\begin{bmatrix}M&\\&a M' \end{bmatrix}\vert M\in\mathrm{GL}_2(F), a \in F^\times\}$$and
$$L_Q(F)=\{\begin{bmatrix}a&&\\&M&\\&&det(M)a^{-1} \end{bmatrix}\vert M\in\mathrm{GL}_2(F), a \in F^\times\}.$$

$P:=L_PB$,(respectively, $Q:=L_QB$) is called as the Siegel (respectively, Klingen) parabolic subgroup of $G$.
We also denote by $\overline{B}$ the subgroup of $G$ consisting of all lower triangular matrices, and denote $\overline{P}=L_P\overline{B}$, $\overline{Q}=L_Q\overline{B}$.

Let $\mathfrak{gsp}_4$(respectively, $\mathfrak{b}$, $\mathfrak{t}$, etc.) denote the Lie algebra of $\mathrm{GSp}_4$(respectively, $B$, $T$, etc.). 
$$\mathfrak{t}(F)=\{\mathrm{diag}(a,b,c-b,c-a)\vert a,b,c\in F\}.$$
$$\mathfrak{gsp}_4(F)=\{A\in \mathrm{M}_4(F) \vert A^TJ+JA=f(A)J, for\ some\ f(A)\in F\}.$$
If 2 is invertible in $F$, then by considering the trace of $J^{-1}A^TJ+A=f(A)I$ we see $f(A)=\frac{tr(A)}{2}$, for all $A\in \mathfrak{gsp}_4(F)$.

    \subsection{Anisotropic flags}
Assume $F$ is a field, char$F\neq 2$. $V$ is a finite dimensional symplectic space over $F$. We say a flag $0\subset V_0\subset V_1 \subset \cdots \subset V_m\subset V$ is anisotropic, if $V_i=V_{m-i}^\perp$, for all $0\leq i\leq m$. We say a subspace $W\subset V$ is anisotropic if $W\subset W^\perp$, that is, the symplectic form of $V$ vanishes after restricted on $W$. 

$\mathrm{GSp}_4(F)$ acts on the set of all anisotropic flags of $(F^4,J)$. $B(F)$ is the stabilizer subgroup of $0\subsetneq <e_1>\subsetneq <e_1,e_2>\subsetneq <e_1,e_2,e_3>\subsetneq F^4$. $P(F)$ is the stabilizer subgroup of $0\subsetneq<e_1,e_2>\subsetneq F^4$. $Q(F)$ is the stabilizer subgroup of $0\subsetneq<e_1>\subsetneq<e_1,e_2,e_3>\subsetneq F^4$.
$\mathrm{GSp}_4(F)/B(F)$ is naturally bijective to the set of all anisotropic complete flags of the symplectic space $(F^4, J)$. $\mathrm{GSp}_4(F)/P(F)$ is naturally bijective to the set of all anisotropic flags in the form of $0\subsetneq V_2\subsetneq V$ (by definition, dim$V_2$=2), or, equivalently, to the set of all 2-dimensional anisotropic subspaces. $\mathrm{GSp}_4(F)/Q(F)$ is naturally bijective to the set of all anisotropic flags in the form of $0\subsetneq V_1\subsetneq V_3\subsetneq V$ (by definition, dim$V_1$=1, dim$V_3$=3), or, equivalently, to the set of all 1-dimensional anisotropic subspaces.

    \subsection{Root data and Weyl group}
Let's describe the root data $(X,\Delta,X^\lor,\Delta^\lor)$ of $(G,B)$. It is of $B_2$=$D_2$ type and has exactly two simple roots, a short one and a long one. Denote $p_1$(respectively, $p_2$, $p_3$): $T\rightarrow \mathbb{G}_m$, $\textrm{diag}(a,b,cb^{-1},ca^{-1})\mapsto a$(respectively, $b$, $c$). Notice that $p_3$ is just the similitude character sim. Clearly, the weight lattice $X=p_1^\mathbb{Z}p_2^\mathbb{Z}p_3^\mathbb{Z}$.The set of simple roots is $\Delta=\{\alpha,\beta\}$, where $\alpha=p_1p_2^{-1}$ and $\beta=p_2^2p_3^{-1}$. The set of roots is $\Phi(G,B)=\{\alpha^{\pm1},\beta^{\pm1},(\alpha\beta)^{\pm1},(\alpha^2\beta)^{\pm1}\}$. $\alpha^\lor(x)=\textrm{diag}(x,x^{-1},x,x^{-1}).$$ \beta^\lor(x)=\textrm{diag}(1,x,x^{-1},1)$. $(\alpha\beta)^\lor=\alpha^\lor\beta^{\lor2}$. $(\alpha^2\beta)^\lor=\alpha^\lor\beta^\lor$. $(\alpha,\alpha^\lor)=(\beta,\beta^\lor)=2$. $(\alpha,\beta^\lor)=-1$. $(\beta,\alpha^\lor)=-2$. $(\text{sim}, \alpha^\lor)=(\text{sim}, \beta^\lor)=0$. $\alpha$ is the short root. $\beta$ is the long root. An weight $p_1^{n_1}p_2^{n_2}p_3^{n_3}$ is dominant (respectively, strictly dominant) if and only if $n_1\geq n_2\geq 0$ (respectively, $n_1>n_2>0$).

Denote $W=N_G(T)/T$ as the Weyl group of $G$. Then $W=<s_1,s_2\vert s_1^2=s_2^2=(s_1s_2)^4=1>$, where$$s_1=\begin{bmatrix}&1&&\\1&&&\\&&&1\\&&1&\end{bmatrix},s_2=\begin{bmatrix}1&&&\\&&1&\\&-1&&\\&&&1\end{bmatrix}.$$The longest element in $W$ is $s_1s_2s_1s_2=s_2s_1s_2s_1\triangleq s_0$.
We can also view $W$ as the subgroup $D_4=\{s\in S_4\vert s(1)+s(4)=s(2)+s(3)=5\}$ of $S_4$ in the following meaning: 
$$w\mathrm{diag}(x_1,x_2,x_3,x_4)w^{-1}=\mathrm{diag}(x_{w^{-1}}(1),x_{w^{-1}}(2),x_{w^{-1}}(3),x_{w^{-1}}(4)).$$

Now let's describe the action of $W$ on $X$. $s_1$ (respectively, $s_2$) is the simple reflection corresponding to $\alpha$ (respectively, $\beta$). $\textrm{sim}=p_3$ is fixed by $W$. $s_1(\alpha)=\alpha^{-1}$. $s_1(\beta)=\alpha^2\beta$. $s_2(\alpha)=\alpha\beta$. $s_2(\beta)=\beta^{-1}$. 

$P$ (respectively, $Q$) is the maximal parabolic subgroup associated to $\beta$ (respectively, $\alpha$). $W_{P}=<\alpha>$. $W_{Q}=<\beta>$.

    \subsection{Classical LLC for $\mathrm{GSp}_4$}
A convenient reference to the classical local Langlands correspondence of $\textrm{GSp}_4$ is \cite{localnewforms}. In this section, I only introduce the LLC of generic principal series. In this case, the LLC is a 1-1 correspondence, although $\textrm{GSp}_4$ has some L-packets with 2 elements.

Let $E$ be a finite extension of $\mathbb{Q}_p$. Consider the following Weil-Deligne representation valued in $\textrm{GSp}_4(E)$:
$$\textrm{diag}(\chi_1\circ\textrm{rec}^{-1},\chi_2\circ\textrm{rec}^{-1},\chi_3\circ\textrm{rec}^{-1},\chi_4\circ\textrm{rec}^{-1}):W_{\mathbb{Q}_p}\rightarrow\textrm{GSp}_4(E),$$$$N=0.$$
$\chi_i\in \textrm{Hom}_{\textrm{sm}}(\mathbb{Q}_p^\times, E^\times)$, $i\in\{1,2,3,4\}$. By definition, it implies $\chi_1\chi_4=\chi_2\chi_3$. We say it is generic if $\chi_i\chi_j^{-1}\neq 1, \vert\cdot \vert_p^{\pm 1}$, for all $1\leq i<j\leq 4$. 

We assume it is generic, then the corresponding irreducible smooth representaion of $\textrm{GSp}_4(\mathbb{Q}_p)$ over $E$ is $$(\text{Ind}_{\overline B(\mathbb{Q}_p)}^{G(\mathbb{Q}_p)}\frac{\phi_1}{\phi_3}(p_1)\frac{\phi_1}{\phi_2}(p_2)\phi_4(p_3)\vert p_1\vert_p^{-2}\vert p_2\vert_p^{-1} \vert p_3\vert_p^{\frac{3}{2}})^\infty.$$

Let's say more about this correspondence. For a commutative ring $R$, consider the following isomorphism of Abelian groups:
\begin{equation}\label{L}
\begin{aligned}
    \mathscr L_R: \mathrm{Hom}(\mathbb Q_p^\times,T(R))&\rightarrow \mathrm{Hom}(T(\mathbb{Q}_p),R^\times),\\
    \mathrm{diag}(\chi_1,\chi_2,\chi_3,\chi_4)&\mapsto \frac{\chi_1}{\chi_3}(p_1)\frac{\chi_1}{\chi_2}(p_2)\chi_4(p_3).
\end{aligned}
\end{equation}
Simply denote $\mathscr L=\mathscr L_E$. We define $w\mapsto \check w$ to be the automorphism of $W$ sending $s_1$ and $s_2$ to each other. By definition, we have $w\mathscr L_R=\mathscr L_R\check w$.
Restricting to the set of algebraic homomorphisms on both side, $\mathscr L$ comes from an isomorphism between the root data of $\mathrm{GSp}_4$ and its dual. Assume $\chi\in \mathrm{Hom}_{\textrm{sm}}(\mathbb Q_p^\times,T(E))$ such that the $\mathrm{GSp}_4$-valued Weil-Deligne representation $\chi\circ\mathrm{rec}^{-1}$ is generic, then the corresponding locally analytic representation is $\mathrm{Ind}^{\mathrm{GSp}_4(\mathbb{Q}_p)}_{\overline B(\mathbb{Q}_p)}(\vert p_1\vert_p^{-2}\vert p_2\vert_p^{-1} \vert p_3\vert_p^{\frac{3}{2}}\mathscr L(\chi))^\infty$, as showing above.

\section{Crystalline ($\varphi,\Gamma$)-modules with $\mathrm{GSp}_4$-structure}

    \subsection{Preliminaries and notations}
Let $E$ be a finite extension of $\mathbb{Q}_p$. Let $\mathcal{R}_E$ be the $E$-coefficient Robba ring. For continuous character $\chi:\mathbb{Q}_p^\times\rightarrow E^\times$, denote by $\mathcal{R}_E(\chi)$ the associated rank one $(\varphi,\Gamma)$-module. We define that a $(\varphi,\Gamma)$-module over Robba ring with $\mathrm{GSp}_4$-structure is a free of rank 4 (as module of $\mathcal{R}_E$) $(\varphi,\Gamma)$-module $D$, equiped with a free of rank 1 (as module of $\mathcal{R}_E$) $(\varphi,\Gamma)$-module $\mathrm{sim}(D)$ and an $(\varphi,\Gamma)$-equivariant alternate $\mathcal{R}_E$-bilinear form $r_D:D\wedge D\rightarrow \mathrm{sim}(D)$, such that the induced map $s_D:D\rightarrow \mathrm{End}_{\mathcal{R}_E}(D,\mathrm{sim}(D))=D^\lor\otimes\mathrm{sim}(D),\ s_D(x)(y)=r_D(x\wedge y)$, is a ($\varphi,\Gamma$)-module isomorphism.  
A filtered vector space, filtered $\varphi$-module, Galois representation or Weil-Deligne representation, etc., with $\mathrm{GSp_4}$-structure is defined similarly. Notice that a vector space with $\mathrm{GSp}_4$-structure is nothing else but a 4-dimensional symplectic space. Let $G$ be a topological group. An isomorphic class of continuous representation over $E$ of $G$ with $\mathrm{GSp}_4$-structure is equivalent to an isomorphic class of continuous representations of $G$ valued in $\mathrm{GSp}_4(E)$, by taking matrices in a fixed standard basis (By a standard basis of a vector space with $\mathrm{GSp}_4$-structure ($V,r_V$), we mean an ordered basis in which the matrix of $r_V$ is linear dependent on $J$).

By definition, it is easy to see that,
\begin{lem}\label{3.1}
Assume $V$ is a filtered vector space with $\mathrm{GSp}_4$-structure. Let h be the integer such that $\mathrm{Fil}^h\mathrm{sim}(V)=\mathrm{sim}(V)$ and $\mathrm{Fil}^{h+1}\mathrm{sim}(V)=0$. Then $\mathrm{Fil}^iV=(\mathrm{Fil}^{h+1-i}V)^\perp$. (Here $(-)^\perp$ means the annihilator subspace in $V$ by $r_V$.)
\end{lem}\qed

Roughly speaking, if there is a functor between two $E$-categories mentioned above preserving duality, wedge product, and "rank", then it also carries $\mathrm{GSp}_4$-structure. For example, $D_{\textrm{cris}}$, $D_{\textrm{dR}}$, $D_{\textrm{rig}}$ and so on.

Assume $D$ is a crystalline $(\varphi,\Gamma)$-module over $\mathcal{R}_E$ with $\mathrm{GSp}_4$-structure. Then $D_{\textrm{cris}}(M)$ is a filtered $\varphi$-module over $E$ with $\mathrm{GSp}_4$-structure. We call $D$ is generic if the $\varphi$-operator on $D_{\textrm{cris}}(D)$ has 4 distinct eigenvalues, and the quotient of each two different eigenvalues is not $p$.  In this case, we may choose eigenvectors $e_1$, $e_2$, $e_3$, $e_4$ such that they form a standard basis. $D_{\textrm{cris}}(D)$=$\oplus_{i=1}^4 Ee_i$. Assume $\varphi e_i= \alpha_ie_i$. The matrix of $\varphi$ in this basis should be in $\mathrm{GSp_4}(E)$, so that $\alpha_1\alpha_4=\alpha_2\alpha_3\triangleq \alpha_0$. By Lem. \ref{3.1}, The Hodge flag is anisotropic and the Hodge-Tate weights $h_1\geq h_2\geq h_3\geq h_4$ satisfy $h_1+h_4=h_2+h_3\triangleq h_0$. $\mathrm{sim}(D_{\textrm{cris}}(D))=(D_{\textrm{cris}}(\mathrm{sim}D))$ has $\varphi$-eigenvalue $\alpha_0$ and Hodge-Tate weight $h_0$. We call $D$ is non-critical if $h_1> h_2> h_3> h_4$ and the Hodge flag of $D_{\textrm{cris}}(D)$ is in general position with the basis $\{e_1,e_2,e_3,e_4\}$. We define a refinement of $D$ as an ordering of $\{\alpha_1,\alpha_2,\alpha_3,\alpha_4\}$ in the form of $(\alpha_{w^{-1}(1)},\alpha_{w^{-1}(2)},\alpha_{w^{-1}(3)},\alpha_{w^{-1}(4)})$, $w\in W=<s\in S_4\vert s(1)+s(4)=s(2)+s(3)=5>$.

Assume $\boldsymbol{\alpha}=(\alpha_1,\alpha_2,\alpha_3,\alpha_4)\in (E^{\times})^4$ satiesfies $\alpha_1\alpha_4=\alpha_2\alpha_3$ and $\alpha_i\alpha_j^{-1}\neq 1,p^{\pm 1}$ for $i\neq j$. Assume $\boldsymbol{h}=(h_1,h_2,h_3,h_4)\in \mathbb{Z}^4$ satiesfies $h_1> h_2> h_3> h_4$ and $h_1+h_4=h_2+h_3$. Let $\mathrm{GSp}_4$-$\Phi_{\mathrm{nc}}(\boldsymbol\alpha,\boldsymbol h)$ denote the collection of all non-critical crystalline $(\varphi,\Gamma)$-modules over $E$ with $\mathrm{GSp}_4$-structure with Hodge-Tate weights $\boldsymbol h$ and a refinement $\boldsymbol\alpha$. Notice that if $D\in \Phi_{\mathrm{nc}}(\boldsymbol\alpha,\boldsymbol h)$, then $\mathrm{End}_{(\varphi,\Gamma)}(D,D)=E$, so up to isomorphism there is at most one $\mathrm{GSp}_4$-structure on $D$. That means acturally $\mathrm{GSp}_4$-$\Phi_{\mathrm{nc}}(\boldsymbol\alpha,\boldsymbol h)$ is a subset of $\Phi_{\mathrm{nc}}(\boldsymbol\alpha,\boldsymbol h)$. Until the end of this subsection, we assume $D\in \mathrm{GSp}_4$-$\Phi_{\mathrm{nc}}(\boldsymbol\alpha,\boldsymbol h)$. The isomorphism class of $D$ is uniquely determined by the relative position of Hodge anisotropic flag with the $\varphi$-eigenspaces in $D_{\textrm{cris}}(D)$, which is parameterized by an Zarisky open set of $T(E)\backslash\mathrm{GSp}_4(E)/B(E)$, $E$ points of a 2 dimensional variety. More explicitly, one can uniquely write $D_{\textrm{cris}}(D)$ in the following form, 
\begin{equation}\label{standard form}
\mathrm{Fil}^iD_{\textrm{cris}}(D)=
\begin{cases}
D_{\textrm{cris}}(D),\ \ i\leq -h_1\\
<a_De_1-e_2+e_3-e_4,b_De_1+(b_D+1)e_2-e_3,e_1+e_2>,\\
\ \ \ \ \ \ \ \ \ \ \ \ \ \ \ \ \ \ \ \ \ \ \ \ \ \ \ \ \ \ \ \ \ \ \ \ \ \ \ \ \ \ \ \ \ \ \ \ \ \ \ \ \ \ \ \ \ \ \ -h_1<i\leq -h_2\\
<a_De_1-e_2+e_3-e_4,b_De_1+(b_D+1)e_2-e_3>,\\
\ \ \ \ \ \ \ \ \ \ \ \ \ \ \ \ \ \ \ \ \ \ \ \ \ \ \ \ \ \ \ \ \ \ \ \ \ \ \ \ \ \ \ \ \ \ \ \ \ \ \ \ \ \ \ \ \ \ \ -h_2<i\leq -h_3\\
<a_De_1-e_2+e_3-e_4>,\ \ -h_3<i\leq -h_4\\
0,\ \ i>-h_4
\end{cases},
\end{equation}
$a_D,b_D\in E, a_Db_D(b_D+1)(a_D+b_D)(a_Db_D+a_D+b_D)\neq 0$, $(e_1,e_2,e_3,e_4)$ is a standard basis and $e_i$ is an $\alpha_i$-eigenvector of $\varphi$. We call the scalars $(a_D,b_D)$ as the Hodge parameters of $D$. $(a_D,b_D)$ is uniquely determined by $D$ and together with $\boldsymbol{\alpha}$ and $\boldsymbol{h}$, $(a_D,b_D)$ uniquely determined the isomorphism class of $D$. Remark that the Hodge parameters rely on a choice of refinement.

    \subsection{Deformations with $\mathrm{GSp}_4$-structure}
    \subsubsection{Deformation functors}
$\boldsymbol\alpha,\boldsymbol h$ are defined above. Assume $(D,r_D)\in \mathrm{GSp}_4$-$\Phi_{\mathrm{nc}}(\boldsymbol\alpha,\boldsymbol h)$. We denote by $\mathcal{C}_E$ the category of all Artinian local $E$-algebras with residue field $E$. We denote by $X_D$ (respctively, $\mathrm{GSp}_4\textrm{-}X_D$)$:\mathcal{C}_E\rightarrow (Sets)$ the ($\varphi,\Gamma$)-module (respectively, ($\varphi,\Gamma$)-module with $\mathrm{GSp_4}$-structure) deformation functor of $D$. By forgetting the $\mathrm{GSp_4}$-structure, we naturally obtain a nature transformation $\mathrm{GSp}_4\textrm{-}X_D\rightarrow X_D$. 

By \cite[Prop. 2.4.1]{Selmergp}, for all $w\in S_4$, there exists a unique triangulation $\mathcal{F}_w:\ 0\subsetneq \mathcal{F}_{w,1}\subsetneq \mathcal{F}_{w,2}\subsetneq \mathcal{F}_{w,3}\subsetneq D$, such that $\mathcal{F}_{w,i}/\mathcal{F}_{w,i-1}\simeq\mathcal{R}_E(z^{h_i}\mathrm{unr}(\alpha_{w^{-1}(i)}))$. By construction in the proof of \cite[Prop. 2.4.1]{Selmergp}, one can easy to see:
\begin{prop}\label{3.2}
If $w\in W$, then $\mathcal{F}=\mathcal{F}^\perp$, that means $\mathcal{F}_{w,i}=
(\mathcal{F}_{w,4-i})^\perp$, $0\leq i\leq 4$. Here $(-)^\perp$ means the perspective submodule in $D$ by $r_D$. 
\end{prop}\qed

We say a filtration $\mathcal{F}$ of $D$ is anisitropic, if one of the following proposition holds for some $w\in W$, 
\\(i) $\mathcal{F}=\mathcal{F}_w$, in this case $\mathcal{F}$ is a triangulation;
\\(ii) $\mathcal{F}=0\subsetneq\mathcal{F}_{w,2}\subsetneq D$, in this case $\mathcal{F}$ is called a $P$-filtration (or Siegel filtration);
\\(iii)  $\mathcal{F}=0\subsetneq\mathcal{F}_{w,1}\subsetneq\mathcal{F}_{w,3}\subsetneq D$, in this case $\mathcal{F}$ is called a $Q$-filtration (or Klingen filtration);
\\(iv) $\mathcal{F}$ is trival. \\
If $\mathcal{F}$ is an anisotropic filtration, by Prop. \ref{3.2}, $\mathcal{F}=\mathcal{{F}}^\perp$. We define a functor $\mathrm{GSp}_4$-$X_{D,\mathcal{F}}: \mathcal{C}_E\rightarrow(Sets)$, $A\mapsto\{(D_A,r_{D_A},\mathcal{F}_A)\vert (D_A,r_{D_A})$ is a ($\varphi,\Gamma$)-module over $\mathcal{R}_A$ lifting $(D, r_D)$, $\mathcal{F}_A=\mathcal{F}_A^\perp$ is a filtration of $D_A$ lifting $\mathcal{F}$ consisting of ($\varphi,\Gamma$)-submodules of $D_A$ which are free direct summands as $R_A$-submodules\}/\{isomorphisms\}.
\begin{prop}\label{3.3}
    $\mathrm{GSp}_4\textrm{-}X_{D,\mathcal{F}}$ is a subfunctor of $X_{D,\mathcal{F}}$(\cite[\S3]{che2011}).
\end{prop}
\proof Assume $A\in\mathcal{C}_E$, $A\neq E$. Then there exists a non-zero element $a\in A$, such that $m_Aa=0$. $(a)=aA=Ea$. For free $\mathcal{R}_A$-module $M$, denote by $i_{aM}$, or more simply, $i_a$, the nature isomorphism $M/m_AM\rightarrow aM$. Denote by $p_a:A\rightarrow A/(a)$ the quotient map. By induction to the length, we only need to proof that for any $(D_{A/(a)},r_{A/(a)},\mathcal{F}_{A/(a)})\in\mathrm{GSp}_4$-$X_{D,\mathcal{F}}(A/(a))$, $$\mathrm{GSp}_4\text{-}X_{D,\mathcal{F}}(p_a)^{-1}(D_{A/(a)},r_{A/(a)},\mathcal{F}_{A/(a)})\rightarrow X_{D,\mathcal{F}}(p_a)^{-1}(D_{A/(a)},\mathcal{F}_{A/(a)})$$ is injective.

If $\mathrm{GSp}_4\text{-}X_{D,\mathcal{F}}(p_a)^{-1}(D_{A/(a)},r_{A/(a)},\mathcal{F}_{A/(a)})$ is empty, then there is nothing to prove. Otherwise we fix a lifting $(D_A,r_A,\mathcal{F}_A)\in\mathrm{GSp}_4\text{-}X_{D,\mathcal{F}}(p_a)^{-1}(D_{A/(a)},r_{A/(a)},\mathcal{F}_{A/(a)})$. Denote by $*_A$ (resp., $\mathrm{sim}(*_A)$),$*\in\{\varphi\}\cup\Gamma$, the $*$-operator on $D_A$ (resp., $\mathrm{sim}(D_A)$). 
Assume $*_A'$, $*\in\{\varphi\}\cup\Gamma$ is another probable group of $*$-operators on the $\mathcal{R}_A$-module $D_A$ which lifts the $(\varphi,\Gamma)$-action on $D_{A/(a)}$ and preserves $\mathcal{F}_A$. 
Then there exists a unique $c_*\in \mathrm{End}_{\mathcal{F}}(D,D)$, such that $*_A'=*_A+i_ac_*i_a^{-1}\circ a*_A$. By an argument similar to the proof of \cite[Prop. 3.6(ii)]{che2011}, $(*_A')_{*\in\{\varphi\}\cup\Gamma}\mapsto\{c_*\}_{*\in\{\varphi\}\cup\Gamma}$ induces a bijection $X_{D,\mathcal{F}}(p_a)^{-1}(D_{A/(a)},\mathcal{F}_{A/(a)})\rightarrow H^1(\mathrm{End}_{\mathcal{F}}(D))$.

Moreover,assume there exist a group of operators $\mathrm{sim}(*_A')$, $*\in\{\varphi\}\cup\Gamma$ on $\mathrm{sim}(D_A)$, such that together with $*_A'$, giving a lifting of ($\varphi,\Gamma$)-module's structure from the bottom of \eqref{diag1} to the top.
\begin{equation}
    \label{diag1}
    \xymatrix{
\wedge^2D_{A}\ar[r]^{r_A}\ar[d]&\mathrm{sim}(D_{A})\ar[d]\\
\wedge^2D_{A/(a)}\ar[r]^{r_{A/(a)}}&\mathrm{sim}(D_{A/(a)})
    }  
\end{equation}
Assume $\mathrm{sim}(*_A')=\mathrm{sim}(*_A)+i_at_ai_a^{-1}\circ a\mathrm{sim}(*_A)$, $*\in\{\varphi\}\cup\Gamma$, $t_a\in \mathrm{End}_{\mathcal{R}_E}(\mathrm{sim}(D))$. Then for all $*\in\{\varphi\}\cup\Gamma$, $x,y\in D_A$, we have
$$r_A(*_A'x,*_A'y)=\mathrm{sim}(*_A')r_A(x,y),$$
$$r_A(*_Ax+i_ac_*i_a^{-1}\circ a*_Ax,*_Ay+i_ac_*i_a^{-1}\circ a*_Ay)=(\mathrm{sim}(*_A)+i_at_ai_a^{-1}\circ a\mathrm{sim}(*_A))r_A(x,y),$$
$$r_A(i_ac_*i_a^{-1}\circ a*_Ax,*_Ay)+r_A(*_Ax+i_ac_*i_a^{-1}\circ a*_Ay)=i_at_ai_a^{-1}\circ a\mathrm{sim}(*_A))r_A(x,y),$$
$$i_ar_D(c_**_D\overline x,*_D\overline y)+i_ar_D(*_D\overline x,c_**_D\overline y)=i_at_*(*_D\overline x,*_D\overline y).$$
\begin{equation}\label{Lie(D)}
    r_D(c_*\overline x,\overline y)+r_D(\overline x,c_*\overline y)=t_*r_D(\overline x,\overline y),\ \ \ *\in\{\varphi\}\cup\Gamma.
\end{equation}
By substituting special x,y (such as, in a standard basis) into \eqref{Lie(D)}, we see that $t_*=\frac{\mathrm{tr}c_*}{2}$.
We denote by $\mathrm{Lie}(D)$ the subsets of elements satisfying \eqref{Lie(D)} in $\mathrm{End}(D)$. A direct calculating shows it is a saturated ($\varphi,\Gamma$)-submodule. $\mathrm{Lie}_\mathcal{F}(D):=\mathrm{Lie}(D)\cap\mathrm{End}_{\mathcal{F}}(D)$. By an similar argument to the proof of \cite[Prop. 3.6(ii)]{che2011}, $(*_A')_{*\in\{\varphi\}\cup\Gamma}\mapsto\{c_*\}_{*\in\{\varphi\}\cup\Gamma}$ induces a bijection $$\mathrm{GSp}_4\textrm{-}X_{D,\mathcal{F}}(p_a)^{-1}(D_{A/(a)},r_{A/(a)},\mathcal{F}_{A/(a)})\xrightarrow{\sim}H^1(\mathrm{Lie}_{\mathcal{F}}(D)).$$
    Then the proposition follows from Lem. \ref{3.4}(iv).
\qed

For $f\in\mathrm{End}(D)$, denote by $f^*:=s_D^{-1}(f^\lor\otimes \mathrm{sim}(D))s_D$ the adjoint operator of $f$. We define
\begin{equation}
        \mathbf s: \mathrm{End}(D)\rightarrow\mathrm{End}(D),\ \
    f\mapsto -f^*+\frac{\mathrm{tr}(f)}{2}\mathrm{\textrm{id}}_D.\label{t}
\end{equation}
 
\begin{lem}\label{3.4}
(i) $\mathbf s$ is a $(\varphi,\Gamma)$-equivariant involution.

(ii) The fixed point set of $\mathbf s$ is $\mathrm{Lie}(D)$.

(iii) $\mathrm{End}_{\mathcal{F}}(D)$ is $\mathbf s$-invariant.

(iv) $H^1(\mathrm{Lie}_{\mathcal{F}}(D))\rightarrow H^1(\mathrm{End}_{\mathcal{F}}(D))$ is injective.

These shows that $\frac{1+\mathbf s}{2}$ is a $(\varphi,\Gamma)$-projection from $\mathrm{End}_{\mathcal{F}}(D)$ to $\mathrm{Lie}_{\mathcal{F}}(D)$.
\end{lem}
\proof
(i)(ii)(iii) is direct by definition. Because of (i)(ii)(iii), $\mathrm{Lie}_{\mathcal{F}}(D)$ is a direct summand of $\mathrm{End}_{\mathcal{F}}(D)$, which proves (iv).

\qed
\begin{lem}\label{3.5}
Let $D_A$ be a $(\varphi,\Gamma)$-module over $\mathcal{R}_A$ lifting $D$. ($D_A$ may not lift the $\mathrm{GSp}_4$-structure.) Let $D_0$ be a saturated $(\varphi,\Gamma)$-submodule of $D$. Then $D_A$ has at most one $(\varphi,\Gamma)$-submodule which is free and direct summand as an $\mathcal{R}_A$-module projecting onto $D_0$.  
\end{lem} 
\proof 
Denote by $D_1$ the $\mathcal{R}_E$-saturated $(\varphi,\Gamma)$-submodule over $\mathcal{R}_A$ of $D_A$ generated by $\sum\limits_{f\in\mathrm{Hom}_{(\varphi,\Gamma)}(D_0,D_A)}\mathrm{Im}f$. Let $D_2$ be a $(\varphi,\Gamma)$-submodule which is free and direct summand as an $\mathcal{R}_A$-module projecting onto $D_0$. We only need to show $D_2=D_1$.
By assumption, $\mathrm{Hom}_{(\varphi,\Gamma)}(D_0,D/D_0)=0$. Notice that $D_A/D_2$ is a lifting of $D/D_0$, hence it is a successful extension of $\mathrm{dim}_E A$ copies of $D/D_0$. Then $\mathrm{Hom}_{(\varphi,\Gamma)}(D_0,D_A/D_2)=0$. Then $\mathrm{Hom}_{(\varphi,\Gamma)}(D_0,D_A)=\mathrm{Hom}_{(\varphi,\Gamma)}(D_0,D_2)$. Because $D_2$ is a lifting of $D_0$, $D_2$ is the $\mathcal{R}_E$-saturated $(\varphi,\Gamma)$-submodule over $\mathcal{R}_A$ of $D_2$ generated by $\sum\limits_{f\in\mathrm{Hom}_{(\varphi,\Gamma)}(D_0,D_2)}\mathrm{Im}f$, which equals to $D_1$
\qed

\begin{prop}\label{3.6}
(i) $X_{D,\mathcal{F}}$ is a subfunctor of $X_D$. 

(ii) $\mathrm{GSp}_4\textrm{-}X_{D,\mathcal{F}}$ is a subfunctor of $\mathrm{GSp}_4\textrm{-}X_D$.

(iii) $\mathrm{GSp}_4\textrm{-}X_{D,\mathcal{F}}=\mathrm{GSp}_4\textrm{-}X_{D}\times_{X_D}X_{D,\mathcal{F}}$.
\end{prop}
\proof (i) This is \cite[Prop. 3.6(i)]{che2011}. (Or by Lem. \ref{3.5}.)

(ii) Combine Prop. \ref{3.3} and Lem. \ref{3.5}.

(iii) Assume $(D_A, r_A)\in\mathrm{GSp}_4\textrm{-}X_{D}(A)$ and $(D_A, \mathcal{F}_A)\in X_{D,\mathcal{F}}(A)$. Then $\mathcal{F}_A^\perp$ projects to $\mathcal{F}^\perp=\mathcal{F}$, too. By Lem. \ref{3.5}, $\mathcal{F}_A=\mathcal{F}_A^\perp$. So $(D_A, r_A, \mathcal{F}_A)\in\mathrm{GSp}_4\textrm{-}X_{D,\mathcal{F}}(A)$.\qed

We simply denote $\mathrm{GSp}_4\textrm{-}X_{D,w}=\mathrm{GSp}_4\textrm{-}X_{D,\mathcal{F}_w}$, $X_{D,w}=X_{D,\mathcal{F}_w}$. We denote $X_{D,g}$ (resp., $\mathrm{GSp}_4\textrm{-}X_{D,g}$) the de Rham deformation subfunctor of $X_D$ (resp., $\mathrm{GSp}_4\textrm{-}X_D$). Then by definition, $\mathrm{GSp}_4\textrm{-}X_{D,g}=\mathrm{GSp}_4\textrm{-}X_{D}\times_{X_D}X_{D,g}$.
\begin{prop}\label{3.7}
$X_{D,*}$, $\mathrm{GSp}_4\textrm{-}X_{D}$, $*\in\{\emptyset,\mathcal{F},g\}$ are all pro-representable and formally smooth .
\end{prop}
\proof This claim for $X_{D,*}$, $*\in\{\emptyset,\mathcal{F},g\}$ is guaranteed by \cite[Prop. 3.15]{che2011}. For $\mathrm{GSp}_4\textrm{-}X_{D}$, it is guaranteed by a standard proceeding using Schlessinger's criterion\cite{schl}. Then the whole claim holds using Prop. \ref{3.6}(iii).
    \subsubsection{Tangent subspaces}
We write $\mathrm{Ext}^i$ (and $\mathrm{Hom}=\mathrm{Ext}^0$) and $H^i$ without $(\varphi,\Gamma)$ in the subscript for the $i$-th extension group and cohomology group of $(\varphi,\Gamma)$-modules (cf.\cite{Liu1}). For a $(\varphi,\Gamma)$-module $M$, we identify elements in $\mathrm{Ext}^1(M,M)$ with deformations of $M$ over $\mathcal{R}_{E[\epsilon]/(\epsilon^2)}$. Indeed, the $E[\epsilon]/(\epsilon^2)$ structure on $\widetilde{M}\in\mathrm{Ext}^1(M,M)$ is given by letting $\epsilon$ act via $\widetilde{M}\twoheadrightarrow M\hookrightarrow \widetilde{M}$. We identify $\mathrm{Ext}^1(M,M)$ with the tangent space of $X_M$ in this way.

Assume $D\in \mathrm{GSp}_4\textrm{-}\Phi_{\mathrm{nc}}(\boldsymbol\alpha,\boldsymbol h)$. We denote by $\mathrm{Ext}^G(D,D)\subset \mathrm{Ext}^1(D,D)$ the subspace of deformations with $\mathrm{GSp}_4$-structure. Assume $\mathcal{F}$ is an anisotropic filtration of $D$, we denote by $\mathrm{Ext}^1_\mathcal{F}(D,D)$ the subspace of $\mathcal{F}$-deformations. For $w\in W\subset S_4$, we donote by $\mathrm{Ext}^1_w(D,D)$ the subspace of trianguline deformations with respect to $w$. We denote by $\mathrm{Ext}^1_g(D,D)$ (resp., $\mathrm{Ext}^1_{g'}(D,D)$) the subspace of de Rham (resp., twisted de Rham) deformations. Denote  $\mathrm{Ext}^G_*(D,D):=\mathrm{Ext}^G(D,D)\cap\mathrm{Ext}^1_*(D,D)$, for $*=w,\mathcal{F},g$ or $g'$.

\begin{prop}\label{3.8}
    There are nature isomorphisms which coincide with the inclusions on both sides:
    (i)$\mathrm{Ext}^G_\mathcal{F}(D,D)\tilde\rightarrow H^1(\mathrm{Lie}_\mathcal{F}(D))$;
    (ii)$\mathrm{Ext}^1_\mathcal{F}(D,D)\tilde\rightarrow H^1(\mathrm{End}_\mathcal{F}(D))$.
\end{prop}
\proof (ii) is \cite[Prop. 3.6(ii)]{che2011}. (i) is proved by a similar argument of the proof of \cite[Prop. 3.6(ii)]{che2011} and Prop. \ref{3.3}. It's easy to see these isomorphisms coincide with inclusions by their definition.\qed
\begin{cor}\label{3.9}
    $\mathrm{dim}_E\mathrm{Ext}^G(D,D)=12$.
    
    $\mathrm{dim}_E\mathrm{Ext}^G_w(D,D)=8$, for $w\in W$. 
    
    $\mathrm{dim}_E\mathrm{Ext}^G_\mathcal{F}(D,D)=9$, for any Siegel or Klingen filtration $\mathcal{F}$.
\end{cor}
\proof Combine Prop. \ref{3.8} and \cite[Thm. 4.3]{Liu1}.\qed

Assume $\widetilde{D}\in\mathrm{Ext}^1(D,D)$, then there exists a unique $\psi_0\in \mathrm{Hom}(\mathbb{Q}_p^\times,E)$, such that $\mathrm{det}(\widetilde D)=\mathrm{det}(D)\otimes_{\mathcal{R}_E}\mathcal{R}_{E[\epsilon]/(\epsilon^2)}(1+2\psi_0\epsilon)\in \mathrm{Ext}^1(\mathrm{det}(D),\mathrm{det}(D))$. Define $\mathrm{sim}(\widetilde D):=\mathrm{sim}(D)\otimes_{\mathcal{R}_E}\mathcal{R}_{E[\epsilon]/(\epsilon^2)}(1+\epsilon\psi_0)\in\mathrm{Ext}^1(\mathrm{sim}(D),\mathrm{sim}(D))$. (Notice that if $\widetilde{D}\in\mathrm{Ext}^G(D,D)$, then this definition of $\mathrm{sim}(\widetilde D)$ coincides with the original one.) $\widetilde D^\lor\otimes_{\mathcal{R}_{E[\epsilon]/(\epsilon^2)}}\mathrm{sim}(\widetilde D)$ is a deformation of $D^\lor\otimes\mathrm{sim}(D)$. Its pullback via $s_D:D\tilde\rightarrow D^\lor\otimes\mathrm{sim}(D)$, which we denote by $\tau (\widetilde D)$, is again a deformation of $D$. 
\begin{prop}\label{3.10}
    (i) $\tau$ is an involution of $\mathrm{Ext}^1(D,D)$.
    
    (ii) The fixed point set of $\tau$ is $\mathrm{Ext}^G(D,D)$.

    (iii) $\mathrm{Ext}^1_{\mathcal{F}}(D,D)$ is $\tau$-invariant.

    These shows that $\frac{1+\tau}{2}$ is a projection from $\mathrm{Ext}^1_{\mathcal{F}}(D,D)$ to $\mathrm{Ext}^G_{\mathcal{F}}(D,D)$.
\end{prop}
\proof By definition, the following square commutes:
$$\xymatrix{
\mathrm{Ext}^1(D,D)\ar[r]^\tau\ar[d]_\sim &\mathrm{Ext}^1(D,D)\ar[d]_\sim\\
H^1(\mathrm{End}(D))\ar[r]^{H^1(\mathbf s)}&H^1(\mathrm{End}(D))}.$$
Here $\mathbf s$ is defined in \eqref{t}. Combine this square and Lem. \ref{3.4}(i)(ii)(iii).\qed

Assume $w\in W$, $\widetilde D\in \mathrm{Ext}^1_w(D,D)$. Assume $\widetilde D$ is a successive extension of $$\mathcal{R}_{E[\epsilon]/(\epsilon^2)}(\mathrm{unr}(\alpha_{w^{-1}(i)})z^{h_i}(1+\epsilon\psi_i)), \psi_i\in \mathrm{Hom}(\mathbb{Q}_p^\times,E), 1\leq i\leq 4.$$ \cite[\S 2.3.1]{Ding} defines a linear map $$\kappa_{w,\mathrm{GL}_4}: \mathrm{Ext}^1_w(D,D)\rightarrow \mathrm{Hom}(\mathbb{Q}_p^\times, \mathfrak{t}_{\mathrm{GL}_4}(E)),$$$$\widetilde D\mapsto \mathrm{diag}(\psi_1,\psi_2,\psi_3,\psi_4).$$(In fact, the original notation in \cite{Ding} of this map is $\kappa_w$, but in this article we keep this notation for our case. In fact, the original definition in \cite{Ding} takes $\mathrm{Hom}(T_{\mathrm{GL}_4}(\mathbb{Q}_p^\times),E)$ as the target of this map. Without harm, we replace it with $\mathrm{Hom}(\mathbb{Q}_p^\times, \mathfrak{t}_{\mathrm{GL}_4}(E))$ for convenience.) 

\begin{prop}\label{3.11}
    $\kappa_{w,\mathrm{GL}_4}(\mathrm{Ext}^G_w(D,D))=\mathrm{Hom}(\mathbb{Q}_p^\times, \mathfrak{t}(E))$.
\end{prop}
\proof By definition we easily see $\kappa_{w,\mathrm{GL}_4}(\mathrm{Ext}^G_w(D,D))\subset\mathrm{Hom}(\mathbb{Q}_p^\times, \mathfrak{t}(E))$. On the other side, assume $\psi\in \mathrm{Hom}(\mathbb{Q}_p^\times, \mathfrak{t}(E))$. By \cite[Prop. 2.10(2)]{Ding}, we can choose a $\widetilde D\in \mathrm{Ext}^1_w(D,D)$, such that $\kappa_{w,\mathrm{GL}_4}(\widetilde D)=\psi$. By definition of $\tau$, $\tau \widetilde D$ is a successive extension of $\mathcal{R}_{E[\epsilon]/(\epsilon^2)}(\mathrm{unr}(\alpha_{w^{-1}(i)})z^{h_i}(1-\epsilon\psi_{5-i}+\frac{1}{2}(\psi_1+\psi_2+\psi_3+\psi_4)))=\mathcal{R}_{E[\epsilon]/(\epsilon^2)}(\mathrm{unr}(\alpha_{w^{-1}(i)})z^{h_i}(1+\epsilon\psi_i))$, so it still maps to $\psi$ via $\kappa_{w,\mathrm{GL}_4}$. So by Prop. \ref{3.10}, $\frac{1+\tau}{2}\widetilde D$ is a preimage of $\psi$ in $\mathrm{Ext}^G_w(D,D)$.\qed 

We denote $\kappa_w: \mathrm{Ext}^G_w(D,D)\twoheadrightarrow \mathrm{Hom}(\mathbb{Q}_p^\times, \mathfrak{t}(E))$ by the restriction of $\kappa_{w,\mathrm{GL}_4}$ .
Define $\mathrm{Hom}_{g'}(\mathbb{Q}_p^\times, \mathfrak{t}(E)):=\mathrm{Hom}_{\textrm{sm}}(\mathbb{Q}_p^\times, \mathfrak{t}(E))+\mathrm{Hom}(\mathbb{Q}_p^\times,\mathfrak{z}(E))$.

\begin{prop}\label{3.12}
    (i) $\mathrm{Ext}^G_g(D,D)$ is the preimage of $\mathrm{Hom}_{\textrm{sm}}(\mathbb{Q}_p^\times, \mathfrak{t}(E))$ via $\kappa_w$, for all $w\in W$.

    (ii) $\mathrm{Ext}^G_{g'}(D,D)$ is the preimage of $\mathrm{Hom}_{g'}(\mathbb{Q}_p^\times, \mathfrak{t}(E))$, for all $w\in W$.

    (iii) $\mathrm{ker}\kappa_w$ doesn't depend on $w\in W$.
\end{prop}
\proof
Restrict \cite[Prop. 2.10(3)]{Ding} and its corollary about $\mathrm{Ext}^G_{g'}(D,D)$ in \cite[\S 2.3.1]{Ding} into the subspace $\mathrm{Ext}^G_w(D,D)$ for deducing (i)(ii).

\cite[\S 2.3.1]{Ding} suggests that $\mathrm{ker}\kappa_{w,\mathrm{GL}_4}$ doesn't depend on $w$. By (i), $\mathrm{ker}\kappa_w=\mathrm{ker}\kappa_{w,\mathrm{GL}_4}\cap \mathrm{Ext}^G_g(D,D)$. (iii) is proved.
\qed

\cite{Ding} defines $\mathrm{Ext}_0^1(D,D):=\mathrm{ker}\kappa_{w,\mathrm{GL}_4}$. Similarly, we define $\mathrm{Ext}_0^G(D,D):=\mathrm{ker}\kappa_w$. By Prop. \ref{3.12}(i), $\mathrm{Ext}_0^G(D,D)=\mathrm{Ext}_0^1(D,D)\cap\mathrm{Ext}^G_g(D,D)=\mathrm{Ext}_0^1(D,D)\cap\mathrm{Ext}^G(D,D)$.
\begin{cor}\label{3.13}
    $\mathrm{dim}_E\mathrm{Ext}_0^G(D,D)=2$. $\mathrm{dim}_E\mathrm{Ext}_g^G(D,D)=5$. $\mathrm{dim}_E\mathrm{Ext}_{g'}^G(D,D)=6$.
\end{cor}\qed
\begin{prop} \label{3.14}
For $w_1,w_2\in W$, the following square commutes:
    $$\xymatrix{
        \mathrm{Ext}_{g'}^G(D,D)\ar@{=}[d] \ar[r]^-{\kappa_{w_1}}_-\sim & \mathrm{Hom}_{g'}(\mathbb{Q}_p^\times, \mathfrak{t}(E))\ar[d]_{w_2w_1^{-1}}^\sim\\
        \mathrm{Ext}_{g'}^G(D,D)\ar[r]^-{\kappa_{w_2}}_-\sim & \mathrm{Hom}_{g'}(\mathbb{Q}_p^\times, \mathfrak{t}(E))    
    }.$$
    
\end{prop}
\proof Restrict the corollary of \cite[Lem. 2.11]{Ding} about $\mathrm{Ext}_{g'}^1(D,D)$ in \cite[\S 2.3.1]{Ding} into the subspace $\mathrm{Ext}_{g'}^G(D,D)$.\qed

Assume $\mathcal{F}_P: 0\subsetneq \mathcal{F}_{w,2} \subsetneq D$ is a $P$-filtration. \cite[\S 2.3.1]{Ding} defines $\mathrm{Ext}^1_{\mathcal{F}_P,g'}(D,D)$ to be the preimage of $\mathrm{Ext}^1_{g'}(\mathcal{F}_{w,2},\mathcal{F}_{w,2})\times\mathrm{Ext}^1_{g'}(D/\mathcal{F}_{w,2},D/\mathcal{F}_{w,2})$ via the nature map $\mathrm{Ext}^1_{\mathcal{F}_P}(D,D)\rightarrow\mathrm{Ext}^1(D_2,D_2)\times\mathrm{Ext}^1(D/D_2,D/D_2)$. We define $\mathrm{Ext}^G_{\mathcal{F}_P,g'}(D,D):=\mathrm{Ext}^G_{\mathcal{F}_P}(D,D)\cap\mathrm{Ext}^1_{\mathcal{F}_P,g'}(D,D)$. Assume $\mathcal{F}_Q: 0\subsetneq \mathcal{F}_{w,1}\subsetneq \mathcal{F}_{w,3} \subsetneq D$ is a $Q$-filtration. Define $\mathrm{Ext}^G_{\mathcal{F}_Q,g'}(D,D)$ in a similar way. 

Define $\mathrm{Hom}_{X,g'}(\mathbb{Q}_p^\times, \mathfrak{t}(E)):=\mathrm{Hom}_{\textrm{sm}}(\mathbb{Q}_p^\times, \mathfrak{t}(E))+\mathrm{Hom}(\mathbb{Q}_p^\times,\mathfrak{z}_{\mathfrak{l}_X}(E))$, here $X=P$ or $Q$.

\begin{prop}\label{3.15}
Let ${\delta_X}:=\begin{cases}
    1,\ X=P\\
    2,\ X=Q
\end{cases}$. For $X\in\{P,Q\}$,

(i) $\mathrm{Ext}^G_{\mathcal{F}_X,g'}(D,D)$ is the preimage of $\mathrm{Hom}_{X,g'}(\mathbb{Q}_p^\times, \mathfrak{t}(E))$ via both $\kappa_w$
 and $\kappa_{s_{{\delta_X}}w}$;
 
(ii) $$\xymatrix{
\mathrm{Ext}_{\mathcal{F}_X,g'}^G(D,D)\ar@{=}[d] \ar[r]^-{\kappa_{w}}_-\sim & \mathrm{Hom}_{X,g'}(\mathbb{Q}_p^\times, \mathfrak{t}(E))\ar[d]_{s_{{\delta_X}}}^\sim\\
\mathrm{Ext}_{\mathcal{F}_X,g'}^G(D,D)\ar[r]^-{\kappa_{s_{{\delta_X}}w}}_-\sim & \mathrm{Hom}_{X,g'}(\mathbb{Q}_p^\times, \mathfrak{t}(E))
}$$
is commutative.
\end{prop}
    \proof Notice that the triangulines of $D$ compatible with $\mathcal{F}_P$ (resp., $\mathcal{F}_Q$) are exactly $\mathcal{F}_w$ and $\mathcal{F}_{s_1w}$ (resp., $\mathcal{F}_w$ and $\mathcal{F}_{s_2w}$). Then the proposition is proved by restricting \cite[Cor. 2.15]{Ding} in $\mathrm{Ext}^G(D,D)$.

The following theorem is a special case of \cite[Cor. 7.13]{density}.
\begin{thm}\label{3.16}
The nature map $\bigoplus\limits_{w\in W}\mathrm{Ext}^G_w(D,D)\rightarrow \mathrm{Ext}^G(D,D)$ is surjective.
\end{thm}
\qed

Denote $\overline{\mathrm{Ext}^\star_*}(D,D):=\mathrm{Ext}^\star_*(D,D)/\mathrm{Ext}^\star_0(D,D)$ for $\star=1,G$ and $*=w,g,g',etc.$. By Prop. \ref{3.12}, for $w\in W$, $\kappa_w$ induces isomorphisms: 
\begin{equation}\label{kappa}
    \begin{aligned}
    \overline{\mathrm{Ext}^G_w}(D,D)&\simeq \mathrm{Hom}(\mathbb{Q}_p^\times, \mathfrak{t}(E)),\\ 
    \overline{\mathrm{Ext}^G_g}(D,D)&\simeq \mathrm{Hom}_{\textrm{sm}}(\mathbb{Q}_p^\times, \mathfrak{t}(E)),\\
    \overline{\mathrm{Ext}^G_{g'}}(D,D)&\simeq \mathrm{Hom}_{g'}(\mathbb{Q}_p^\times, \mathfrak{t}(E)),\\
    \overline{\mathrm{Ext}^G_{\mathcal{F}_X,g'}}(D,D)&\simeq \mathrm{Hom}_{X,g'}(\mathbb{Q}_p^\times, \mathfrak{t}(E)), X\in\{P,Q\}.
    \end{aligned}
\end{equation} We still denote these isomorphisms by $\kappa_w$.
Denote by $j:\mathop{\bigoplus}\limits_{w\in W}\mathrm{Hom}(\mathbb{Q}_p^\times,\mathfrak{t}(E))_w\twoheadrightarrow \overline{\mathrm{Ext}}^G(D,D)$, the composition of the natural map $\mathop{\bigoplus}\limits_{w\in W}\overline{\mathrm{Ext}}^G_w(D,D)\twoheadrightarrow \overline{\mathrm{Ext}}^G(D,D)$ with $\mathop{\boxplus}\limits_{w\in W}\kappa_w^{-1}$. Here $\mathrm{Hom}(\mathbb{Q}_p^\times,\mathfrak{t}(E))_w$ is just a copy of $\mathrm{Hom}(\mathbb{Q}_p^\times,\mathfrak{t}(E))$.
We will prove the following theorem in the next subsection.
\begin{thm}\label{3.17}
$\mathrm{ker} j$ determines the Hodge parameters of $D$.
\end{thm}

    \subsection{Remeet Hodge parameters}
Assume $A\in \mathcal{C}_E$, a $\mathbb{G}_a$-representation over $A$ is defined as a finite free $A$-module $V_A$, equipped with a nilpotent $A$-endomorphism $\nu_{V_A}$. For an almost de Rham ($\varphi,\Gamma$)-module M over $\mathcal{R}_A$, $D_{\textrm{pdR}}(M):=(B_{\textrm{pdR}}\otimes_{B_{\textrm{dR}}^+} W^{+}_{\textrm{dR}}(M))^{\mathrm{Gal}_{\mathbb{Q}_p}}$ is a filtered $\mathbb{G}_a$-representation over $A$. We denote the equipped nilpotent operater of $D_{\textrm{pdR}}(M)$ by $\nu_{\textrm{pdR}}(M)$. $\nu_{\textrm{pdR}}(M)=0$ iff $M$ is de Rham. If $M$ is over $\mathcal{R}_E$ and crystalline, then by definition $D_{\textrm{pdR}}(M)=D_{\textrm{cris}}(M)$, which is an equality of filtered $E$-vector spaces. (For details, see \cite[\S 3]{model} or \cite{font}.)

Assume $D\in \mathrm{GSp}_4\textrm{-}\Phi_{\mathrm{nc}}(\boldsymbol\alpha,\boldsymbol h)$. For $\widetilde D\in \mathrm{Ext}^1(D,D)$, $D_{\textrm{pdR}}(\widetilde D)$ is a deformation of $D_{\textrm{pdR}}(D)=D_{\textrm{cris}}(D)$ over $E[\epsilon]/(\epsilon^2)$, as a filterd $\mathbb{G}_a$-representation. Denote by $i_\epsilon$ the nature isomorphism $\epsilon D_{\textrm{pdR}}(\widetilde D)\rightarrow D_{\textrm{pdR}}(\widetilde D)/\epsilon D_{\textrm{pdR}}(\widetilde D)=D_{\textrm{cris}}(D)$. Then there exists a unique $\nu_{\widetilde D}\in \mathrm{End}_{F_H}(D_{\textrm{cris}}(D))$ such that $\nu_{\textrm{pdR}}(\widetilde D)=i_\epsilon^{-1}\nu_{\widetilde D}i_\epsilon\epsilon$. Here $F_H$ means the Hodge flag of $D_{\textrm{cris}}(D)$. $\nu(\widetilde D)=0$ iff $\widetilde D$ is de Rham. Moreover, if $\widetilde D\in\mathrm{Ext}^G(D,D)$, then $D_{\textrm{pdR}}(\widetilde D)$ also admits a lifting of the $\mathrm{GSp_4}$-structure, and it means $\nu(\widetilde D)\in \mathfrak{gsp}(D_{\textrm{cris}}(D))$.

Thus we obtain a linear map $\nu$: $\overline{\mathrm{Ext}^G(D,D)}\rightarrow \mathfrak{gsp}_{F_H}(D_{\textrm{cris}}(D))$. $\overline{\mathrm{Ext}}^G_g(D,D)$ is its kernel. By considering dimensions, we see $\nu$ is surjective.

We choose an ordered basis $(e_1,e_2,e_3,e_4)$ of $D_{\textrm{cris}}(D)$ satisfying \eqref{standard form}. We denote by $F_w$ the complete flag $0\subsetneq<e_{w^{-1}(1)}>\subsetneq<e_{w^{-1}(1)},e_{w^{-1}(2)}>\subsetneq...\subsetneq D_{\textrm{cris}}(D)$. We denote by $F_w^i$ (resp., $F_H^i$) the $i$-dimensional term in $F_w$(resp., $F_H$). Because $D$ is non-critical, $F_H$ is in general position with $F_w$ for all $w\in S_4$. $D_{\textrm{cris}}(D)=\bigoplus\limits_{i=1}^4F_w^i\cap F_H^{5-i}$, each direct summand has dimension $1$.

\begin{lem}\label{3.18}
    Assume $w\in S_4$, $\widetilde D\in \mathrm{Ext}^1_w(D,D)$. Assume $\widetilde D$ is a successful extension of $\mathcal{R}_{E[\epsilon]/(\epsilon^2)}(\mathrm{unr}(\alpha_{w^{-1}(i)})z^{h_i}(1+\epsilon(\psi_i^{\textrm{sm}}+a_i\mathrm{log}_p)))$, $1\leq i\leq 4$. Then $\nu(\widetilde D)$ is diagonalizable and acts as the scalar $a_i$ on $F_w^i\cap F_H^{5-i}$, for $1\leq i\leq 4$. This uniquely determines $\nu(\widetilde D)$.
\end{lem}
\proof Without lose of generality, we may assume $\psi_i^{\textrm{sm}}=0,\ a_i=1$, by Prop. \ref{3.12}(1). $\widetilde D$ has a filtration lifting $\mathcal{F}_w$, so $D_{\textrm{pdR}}(\widetilde D)$ has a flag lifting $D_{\textrm{pdR}}(\mathcal{F}_w)=F_w$. This means $\nu(\widetilde D)$
 preserves $F_w$. $\nu(\widetilde D)$ preserves both $F_w$ and $F_H$, so $F_H^{5-i}\cap F_w^{i}$, $1\leq i\leq 4$ are eigenspaces of $\nu(\widetilde D)$. That the eigenvalues equal to 1 follows from the fact
 $D_{\textrm{pdR}}(\mathcal{R}_{E[\epsilon]/(\epsilon^2)}(1+\mathrm{log}_p))=(E[\epsilon]/(\epsilon^2),\epsilon)$ (viewed as a deformation of $D_{\textrm{pdR}}(\mathcal{R}_E)=(E,0)$).\qed

 Now we can prove Thm. \ref{3.17}.

\proof of Thm. \ref{3.17}. 
Consider the following commutative diagram.
$$\xymatrix@C=1.1em{
\bigoplus\limits_{w\in W}\mathrm{Hom}(\mathbb{Q}_p^\times,\mathfrak{t}(E))_w
\ar@{->>}[d]_{\boxplus_{w\in W}l}
\ar@/^2em/@{->>}[rr]^j
&\bigoplus\limits_{w\in W}\overline{\mathrm{Ext}}^G_w(D,D)
\ar@{->>}[r]
\ar@{->>}[d]
\ar[l]^-\sim_-{\boxplus_{w\in W}\kappa_w}
&\overline{\mathrm{Ext}}^G(D,D)
\ar@{->>}[d]\\
\bigoplus\limits_{w\in W}\mathfrak{t}(E)_w
\ar@{-->>}[rrd]_{\exists! \overline j}
&\bigoplus\limits_{w\in W}{\overline{\mathrm{Ext}}^G_w(D,D)}/{\overline{\mathrm{Ext}}^G_g(D,D)}
\ar@{->>}[r]
\ar[l]^-\sim_-{\boxplus_{w\in W}\overline{\kappa_w}}
&{\overline{\mathrm{Ext}}^G(D,D)}/{\overline{\mathrm{Ext}}^G_g(D,D)}
\ar[d]^\nu_\sim\\
&&\mathfrak{gsp}_{F_H}(D_{\textrm{cris}}(D))
}$$
Assume the deformation $\widetilde{D}$ in Lem. \ref{3.18} is in $\mathrm{Ext}^G(D,D)$, then it maps to $\mathrm{diag}(a_1,\allowbreak a_2,\allowbreak a_3,\allowbreak a_4)$ via $\overline{\kappa_w}$. $l:\mathrm{Hom}(\mathbb{Q}_p^\times,\mathfrak{t}(E))\rightarrow \mathfrak{t}(E)$ is given by $l(\psi)=\psi(\mathrm{exp}_p(1))$. We only need to prove that $\mathrm{ker}\overline j$ determines $a_D$, $b_D$.

For $x\in\mathfrak{t}(E), w'\in W$, we denote by $(x)_{w'}\in\mathfrak{t}(E)_{w'}\subset \bigoplus\limits_{w\in W}\mathfrak{t}(E)_w$ the copy of $x$ in $\mathfrak{t}(E)_{w'}$. Define $t_1=\mathrm{diag}(-1,-1,1,1)$, $t_2=\mathrm{diag}(-1,0,0,1)$. Consider the following 8 elements in $\bigoplus\limits_{w\in W}\mathfrak{t}(E)_w$:
$$f_1=(t_1)_{\textrm{id}},\ f_2=(t_1)_{s_2},\ f_3=(t_1)_{s_0},\ f_4=(t_1)_{s_2s_1},$$
$$g_1=(t_2)_{\textrm{id}},\ g_2=(t_2)_{s_1},\ g_3=(t_2)_{s_1s_2},\ g_4=(t_2)_{s_0}.$$

By Lem. \ref{3.18}, we can explicitly calculate the image of them via $\overline j$. Choose a standard basis $\mathcal{B}:=(v_1,v_2,v_3,v_4)$ of $D_{\textrm{cris}}(D)$, where
$$v_1=a_De_1-e_2+e_3-e_4;$$
$$v_2=b_De_1+(b_D+1)e_2-e_3;$$
$$v_3=e_1+e_2;$$
$$v_4=e_1.$$
Then,
$$[f_1]_\mathcal{B}=\begin{bmatrix}
1&&&\\
&1&&\\
&&-1&\\
&&&-1
\end{bmatrix},
[f_2]_\mathcal{B}=\begin{bmatrix}
1&&&\\
&1&\frac{2}{b+1}&\\
&&-1&\\
&&&-1
\end{bmatrix},$$
$$[f_3]_\mathcal{B}=\begin{bmatrix}
1&&\frac{2}{ab+a+b}&\frac{2(b+1)}{ab+a+b}\\
&1&\frac{2(a+1)}{ab+a+b}&\frac{2}{ab+a+b}\\
&&-1&\\
&&&-1
\end{bmatrix},
[f_4]_\mathcal{B}=\begin{bmatrix}
1&&\frac{2}{a+b}&\frac{2}{a+b}\\
&1&\frac{2}{a+b}&\frac{2}{a+b}\\
&&-1&\\
&&&-1
\end{bmatrix},$$
$$[g_1]_\mathcal{B}=\begin{bmatrix}
1&&&\\
&0&&\\
&&0&\\
&&&-1
\end{bmatrix},
[g_2]_\mathcal{B}=\begin{bmatrix}
1&-1&&\\
&0&&\\
&&0&1\\
&&&-1
\end{bmatrix},$$
$$[g_3]_\mathcal{B}=\begin{bmatrix}
1&-(b+1)&-1&\\
&0&&-1\\
&&0&b+1\\
&&&-1
\end{bmatrix},
[g_4]_\mathcal{B}=\begin{bmatrix}
1&\frac{b}{a}&\frac{1}{a}&\frac{2}{a}\\
&0&&\frac{1}{a}\\
&&0&-\frac{b}{a}\\
&&&-1
\end{bmatrix}.$$
(Here we omit the subscript of $a_D$, $b_D$.)

Then the projection of $\mathrm{ker}\overline j\cap<f_1,f_2,f_3,f_4,g_1,g_2,g_3>$
 on $<g_2,g_3>$ is the line $<(b_D+1)g_2-g_3>$. Thus $\mathrm{ker}\overline j$ determines $b_D$. Similarly, the projection of $\mathrm{ker}\overline j\cap<f_1,f_2,f_3,f_4,g_1,g_2,g_4>$
 on $<g_2,g_4>$ is the line $<b_Dg_2+a_Dg_4>$. Thus $\mathrm{ker}\overline j$ determines $a_D$, too.\qed
\begin{rem}\label{3.19}
We give a reductive-group-style description of $\overline j$. This implies our method in the local theory can be generalized to much more general reductive group (which is what I am working on and will come soon). In this remark, denote $G=\mathrm{GSp}_4(D_{\textrm{cris}}(D))$. Denote by $B$ the Borel subgroup preserving $F_{\textrm{id}}$ and by $N$ its unipotent subgroup. Denote by $T$ the maximal torus preserving $Ee_i,\ 1\leq i\leq 4$. Denote by $B_H$ the Borel subgroup preserving $F_H$. Here, all the groups are viewd as Lie groups over $E$. 

By assumption, $wBw^{-1}\cap B_H$ is a maximal torus for all $w\in W$. There exists a unique $n_{w,H}\in wNw^{-1}$, such that $n_{w,H}Tn_{w,H}^{-1}=wBw^{-1}\cap B_H$. We denote $f_w:T\rightarrow G, t\mapsto n_{w,H}wtw^{-1}n_{w,H}^{-1}$. We identify $\mathfrak{t}(E)$ with $\mathrm{Lie}(T)$. Then $\overline j$=$\Sigma_{w\in W} T_ef_w:\bigoplus\limits_{w\in W}\mathrm{Lie}(T)\rightarrow \mathrm{Lie}(G)$. We have proven its kernel uniquely determines $F_H$ (or equivalently, $B_H$) up to conjugations by elements in $T$.
\end{rem}
\begin{rem}\label{3.20}
By the direct calculation above, we can also show $\overline j$ is surjective. This gives another proof of Thm. \ref{3.16}.
\end{rem}

\section{Locally analytic crystalline representations of $\mathrm{GSp}_4(\mathbb{Q}_p)$}
    \subsection{Locally analytic representations of $\mathrm{GSp}_4(\mathbb{Q}_p)$ and extensions}
    \subsubsection{Preliminaries and notations}
We retain the notations in \S2. For weight $\mu$ of $T$, denote by $\overline M(\mu):=\mathrm{U}(\mathfrak{gsp}_4)\otimes_{\mathrm U(\overline{\mathfrak b})}\mu$, and $\overline L(\mu)$ its unique irreducible quotient. If $\mu$ is anti-dominant, then $\overline L(\mu)$ is finite dimensional and isomorphic to the dual $L(-\mu)^\lor$, the algrbraic representation of $\mathrm{GSp}_4$ of highest weight $-\mu$.

For an admissible locally analytic representation $V$ of a $p$-adic locally analytic group $G$, by \cite{schn2002}, its dual $V^\lor$ is naturally a module over the distribution algebra $\mathcal{D}(G,E)$, which, equipped with the strong topology, is a coadmissible module over $\mathcal{D}(H,E)$ for a(ny) compact open subgroup $H$ of $G$. For admissible locally analytic representations $V_1$, $V_2$ of $G$, denote by $\mathrm{Ext}^i_G(V_1,V_2)$(or simply, $\mathrm{Ext}^i(V_1,V_2)$)$:= \mathrm{Ext}^i_{\mathcal{D}(G,E)}(V_2^\lor,V_1^\lor)$, where the latter is defined in the Abelian category of abstract $\mathcal D(G,E)$-modules. By \cite[Lem. 2.1.1]{Ext1}, $\mathrm{Ext}^1_G(V_1,V_2)$ is equal to the extension group of admissible locally analytic representations of $V_1$ by $V_2$. Note that any representation $\widetilde V$ in $\mathrm{Ext}^1(V,V)$ is equipped with a natural $E[\epsilon]/(\epsilon^2)$ structure where $\epsilon$ acts via $\widetilde V\twoheadrightarrow V\hookrightarrow \widetilde V$. If $V$ is locally algebraic, define $\mathrm{Ext}^1_g(V,V)$ to be the subgroup of locally algebraic extensions.

$\varepsilon: \mathbb{Q}_p^\times\rightarrow E^\times,\ x\mapsto x\vert x\vert_p$. $\eta=\vert p_1\vert_p^{-2}\vert p_2\vert_p^{-1}: T(\mathbb{Q}_p)\rightarrow E^\times$. $\delta_B$ is the modulus character of $B(\mathbb{Q}_p)$. Let $\phi=\phi_1(p_1)\phi_2(p_2)\phi_3(p_3):T(\mathbb{Q}_p)\rightarrow E^\times$ be a smooth character. We call $\phi$ is generic if $\phi_1,\phi_2,\phi_1\phi_2,\phi_1/\phi_2\neq 1, \vert\cdot\vert_p^{\pm1}$. Let $\pi_{\textrm{sm}}(\phi):=(\mathrm{Ind}^{\mathrm{GSp}_4(\mathbb{Q}_p)}_{\overline B(\mathbb{Q}_p)}\phi\eta)^\infty$, which is an absolutely irreducible smooth admissible representation of $\mathrm{GSp}_4(\mathbb{Q}_p)$ when $\phi$ is generic. Moreover, when $\phi$ is generic, $\pi_{\textrm{sm}}(\phi)\simeq\pi_{\textrm{sm}}(w(\phi))$ for all $w\in W$, which, after twisting by $\vert\mathrm{sim}\vert_p^{\frac{3}{2}}$, becomes the irreducible smooth representation corresponding to the $\mathrm{GSp}_4$-valued Weil-Deligne representation $\mathrm{diag}(\phi_3,\phi_1\phi_3,\phi_2\phi_3,\phi_1\phi_2\phi_3)\circ\mathrm{rec}^{-1}$ in the meaning of the classical local Langlands correspondance of $\mathrm{GSp}_4$.

We identify $L_P$ with $\mathrm{GL_2}\times \mathrm{GL_1}$ by the isomorphism \begin{equation}
\label{LP}\begin{bmatrix}M&\\&aM'\end{bmatrix}\mapsto (M,a).
\end{equation}
Assume $V_i$ is a locally analytic representation of compact type of $\mathrm{GL}_i(\mathbb{Q}_p)$, $i=1,2$. We view the (completed) tensor product $V_2\boxtimes(\hat\boxtimes)V_1$ as a representation of $L_P(\mathbb{Q}_p)$ through \eqref{LP}. Through \eqref{LP}, $T$ is identified with $\mathbb{G}_m^2\times\mathbb{G}_m$ by 
\begin{equation}\label{TLP}
    T\rightarrow \mathbb{G}_m^2\times\mathbb{G}_m,
\mathrm{diag}(a,b,cb^{-1}.ca^{-1})\mapsto ((a,b),c).
\end{equation}
Denote by $s$ the non-trivial element in $W_{\mathrm{GL_2}}$ acting on $\mathbb{G}_m^2$. Then the action of $s$ on $\mathbb{G}_m^2\times\mathbb{G}_m$ 
coincides with $s_1$ on $T$ through \eqref{TLP}. For $i=1,2$, let $\phi^P_i$ be the character of $\mathbb{G}_m^i$, such that $\phi=\phi^P_2\boxtimes\phi^P_1$ in the meaning of \eqref{TLP}. We have $(s_1(\phi))^P_2=s(\phi^P_2)$, $(s_1(\phi))^P_1=\phi^P_1$. Similarly, we identify $L_Q$ with $\mathrm{GL_1}\times\mathrm{GL_2}$ by the isomorphism \begin{equation}
\label{LQ}\begin{bmatrix}a&&\\&M&\\&&\mathrm{det}(M)a^{-1}\end{bmatrix}\mapsto (a,M).\end{equation}Assume $V_i$ is a locally analytic representation of compact type of $\mathrm{GL}_i(\mathbb{Q}_p)$, $i=1,2$. We view the (completed) tensor product $V_1\boxtimes(\hat\boxtimes)V_2$ as a representation of $L_Q(\mathbb{Q}_p)$ through \eqref{LQ}. Through \eqref{LQ}, $T$ is identified with $\mathbb{G}_m^2\times\mathbb{G}_m$ by 
\begin{equation}\label{TLQ}
    T\rightarrow \mathbb{G}_m\times\mathbb{G}_m^2,
\mathrm{diag}(a,b,cb^{-1}.ca^{-1})\mapsto (a,(b,cb^{-1})).
\end{equation}
The action of $s$ on $\mathbb{G}_m\times\mathbb{G}_m^2$ coincides with $s_2$ on $T$ through \eqref{TLQ}. For $i=1,2$, let $\phi^Q_i$ be the character of $\mathbb{G}_m^i$, such that $\phi=\phi^Q_1\boxtimes\phi^Q_2$ in the meaning of \eqref{TLQ}. We have $(s_2(\phi))^Q_2=s(\phi^Q_2)$, $(s_2(\phi))^Q_1=\phi^Q_1$.

    \subsubsection{Principal series}
We collect some facts on the locally analytic principal series of $\mathrm{GSp}_4(\mathbb{Q}_p)$. This subsection is almost a copy of \cite[\S 3.1.2]{Ding}.

Let $\boldsymbol h=(h_1,h_2,h_3,h_4)\in \mathbb{Z}^4$, such that $h_1>h_2>h_3>h_4$, $h_1+h_4=h_2+h_3$. Put $\lambda:=p_1^{h_1-h_3-2}p_2^{h_1-h_2-1}p_3^{h_4}=\mathscr L(z^{\boldsymbol h})p_1^{-2}p_2^{-1}$, which is a dominant weight of $T$. Assume $\phi:T(\mathbb{Q}_p)\rightarrow E^\times$ is a generic smooth character. Define $\pi_{\textrm{alg}}(\phi,\boldsymbol h)=\pi_{\textrm{sm}}(\phi)\otimes L(\lambda)$, which is a locally algebraic representation of $\mathrm{GSp}_4(\mathbb{Q}_p)$. $\pi_{\textrm{alg}}(\phi,\boldsymbol h)\simeq\pi_{\textrm{alg}}(w(\phi),\boldsymbol h)$, for all $w\in W$. For $w\in W$, define $\mathrm{PS}(w(\phi),\boldsymbol h):=(\mathrm{Ind}^{\mathrm{GSp}_4(\mathbb{Q}_p)}_{\overline B(\mathbb{Q}_p)}w(\phi)\eta\lambda)^{\mathrm{an}}$. We have (where $\mathcal{F}_{\overline B}^\mathrm{GSp}(-,-)$ denotes the Orlik-Strauch functor\cite{OS}, and $w\cdot\lambda$ denotes the dot action\cite[\S1.8]{BGG}):
\begin{prop}\label{4.1}
    (i) The irreducible constituents of $\mathrm{PS}(w(\phi),\boldsymbol h)$ are given by $$\{\mathscr C(w,u):=\mathcal{F}_{\overline B}^\mathrm{GSp}(\overline{L}(-u\cdot\lambda),w(\phi)\eta)\ \vert\ u\in W\},$$
    which are pairwisely distinct. Moreover, if $\mathrm{lg}(u)=1$, then $\mathscr C(w,u)$ has multiplicity one.
    
    (ii) $\mathrm{soc}\ \mathrm{PS}(w(\phi),\boldsymbol h)\simeq\pi_{\textrm{alg}}(\phi,\boldsymbol h)$.
    
    (iii) $\mathrm{soc}\ (\mathrm{PS}(w(\phi),\boldsymbol h)/\pi_{\textrm{alg}}(\phi,\boldsymbol h))\simeq \mathscr C(w,s_1)\oplus\mathscr C(w,s_2)$.
    
    (iv) For $w,w'\in W$, $u,u'\in \{s_1,s_2\}$, $\mathscr C(w,u)=\mathscr C(w',u')$ iff $u=u'=s_i$ and $w(w')^{-1}=s_j$, for $j\neq i\in\{1,2\}$.
\end{prop}
\proof
 (i) and (iv) follow from \cite{OS} (together with some standard facts on the constituents of the
 Verma module, see for example \cite{BGG}). (ii) (iii) follows from \cite[Cor. 2.5]{Br1} or \cite[Prop. 4.1.2]{OS1}.\qed

For $i\in\{1,2\}$, let $I\subset\{1,2,3,4\}$ be a subset of cardinality $i$, such that the sum of elements in $I$ is not $5$. We see that all the representations $\mathscr C(w,s_i)$ with $w^{-1}(\{1,i\})=I$ are isomorphic, which we denote by $\mathscr C(I,s_i)$. By Prop. \ref{4.1}(iv), $\mathscr C(I,s_i)$ are pairwisely distinct for different $I$. For $w\in W$ such that $w^{-1}(\{1,i\})=I$, we have (by \cite[Cor. 2.5]{Br1} or \cite[Prop. 4.1.2]{OS1})
\begin{equation}
    \label{socle}
    \mathscr C(I,s_i)\simeq \mathrm{soc}\ (\mathrm{Ind}^{\mathrm{GSp}_4(\mathbb{Q}_p)}_{\overline B(\mathbb{Q}_p)}w(\phi)\eta(s_i\cdot\lambda))^{\mathrm{an}}.
\end{equation}
\begin{lem}\label{4.2}
    Let $w\in W,i\in\{1,2\}$, such that $w^{-1}\{1,i\}=I$.

    (i) We have $\mathrm{Hom}_{T(\mathbb{Q}_p)}((s_i\cdot\lambda)w(\phi)\eta\delta_B,J_B(\mathscr C(I,s_i)))\simeq E$. Here $J_B$ is the Emerton's Jacquet functor \cite{Jacq1}.

    (ii) We have $I^{\mathrm{GSp}_4(\mathbb{Q}_p)}_{\overline B(\mathbb{Q}_p)}((s_i\cdot\lambda)w(\phi)\eta)\simeq\mathscr C(I,s_i)$. Here $I^{\mathrm{GSp}_4(\mathbb{Q}_p)}_{\overline B(\mathbb{Q}_p)}$ is Emerton’s induction functor \cite{Jacq2}.
\end{lem}
\proof
(i) follows from \cite{Jacq1}, \cite[Prop. 2.1.2]{Jacq2}, \cite[Thm. 4.3, Rem. 4.4(i)]{Br2} and an argument same to \cite[Lem. 3.2(1)]{Ding}.

(ii) follows from \cite{Jacq1}, \cite[Prop. 2.1.2]{Jacq2}, \cite[Thm. 4.3, Rem. 4.4(i)]{Br2}, \cite[Prop. 4.1.2]{OS1} and an argument same to \cite[Lem. 3.2(2)]{Ding}.
\qed

Let $\mathrm{PS}_1(w(\phi),\boldsymbol h)$ be the unique subrepresentation of $\mathrm{PS}(w(\phi),\boldsymbol h)$ of socle $\pi_{\textrm{alg}}(\phi,\boldsymbol h)$ and cosocle $\mathscr C(w,s_1)\oplus\mathscr C(w,s_2)$. Consider the amalgamated sum $\mathop\oplus\limits_{\substack{\pi_{\textrm{alg}}(\phi,\boldsymbol h)\\w\in W}}\mathrm{PS}_1(w(\phi),\boldsymbol h)$. It admits a unique quotient, denoted by $\pi_1(\phi,\boldsymbol{h})$ of socle $\pi_{\textrm{alg}}(\phi,\boldsymbol h)$ (by Lem. \ref{4.4}(i)). By Prop. \ref{4.1}(iii)(iv), $\pi_1(\phi,\boldsymbol h)$ is an extension of $\oplus_{\substack{i=1,2,I\subset\{1,2,3,4\},\\\#I=i,\Sigma_{j\in I}j\neq 5}}\mathscr C(I,s_i)$ (which are 8 in total) by $\pi_{\textrm{alg}}(\phi,\boldsymbol h)$. We fix the tautological injection $\pi_{\textrm{alg}}(\phi,\boldsymbol h)\hookrightarrow\pi_1(\phi,\boldsymbol h)$. We study
 the extension group of $\pi_1(\phi,\boldsymbol h)$ by $\pi_{\textrm{alg}}(\phi,\boldsymbol h)$. 
 
 Let $\mathrm{Hom}_{g'}(T(\mathbb{Q}_p),E):=\mathrm{Hom}_{\textrm{sm}}(T(\mathbb{Q}_p),E)+\mathrm{Hom}(\mathbb{Q}_p^\times,E)\circ\mathrm{sim}$. We have:
\begin{prop}\label{4.3}
    (i) For $w\in W$, the following map is a bijection:
    $$\zeta_w: \mathrm{Hom}_{g'}(T(\mathbb{Q}_p),E)\mathop\rightarrow\limits^{\sim}\mathrm{Ext}^1(\pi_{\textrm{alg}}(\phi,\boldsymbol h),\pi_{\textrm{alg}}(\phi,\boldsymbol h)),\ \psi\mapsto I^{\mathrm{GSp}_4(\mathbb{Q}_p)}_{\overline B(\mathbb{Q}_p)}(w(\phi)\eta\lambda(1+\psi\epsilon)),$$ and induces a bijection $\mathrm{Hom}_{\textrm{sm}}(T(\mathbb{Q}_p),E)\mathop\rightarrow\limits^{\sim}\mathrm{Ext}^1_{\textrm{lalg}}(\pi_{\textrm{alg}}(\phi,\boldsymbol h),\pi_{\textrm{alg}}(\phi,\boldsymbol h))$. In particular, we have $\mathrm{dim}_E\mathrm{Ext}^1_{\textrm{lalg}}(\pi_{\textrm{alg}}(\phi,\boldsymbol h),\pi_{\textrm{alg}}(\phi,\boldsymbol h))=3$ and $\mathrm{dim}_E\mathrm{Ext}^1(\pi_{\textrm{alg}}(\phi,\boldsymbol h),\pi_{\textrm{alg}}(\phi,\boldsymbol h))\allowbreak =4$.

    (ii) For $w_1,w_2\in W$, the following square commutes:
    $$\xymatrix{
    \mathrm{Hom}_{g'}(T(\mathbb{Q}_p),E)\ar[r]^-{\zeta_{w_1}}_-\sim\ar[d]_{w_2w_1^{-1}}^\sim&\mathrm{Ext}^1(\pi_{\textrm{alg}}(\phi,\boldsymbol h),\pi_{\textrm{alg}}(\phi,\boldsymbol h))\ar@{=}[d]\\
    \mathrm{Hom}_{g'}(T(\mathbb{Q}_p),E)\ar[r]^-{\zeta_{w_2}}_-\sim&\mathrm{Ext}^1(\pi_{\textrm{alg}}(\phi,\boldsymbol h),\pi_{\textrm{alg}}(\phi,\boldsymbol h))
    }.$$
\end{prop}
\proof
It follows from the definition of $I^{\mathrm{GSp}_4(\mathbb{Q}_p)}_{\overline B(\mathbb{Q}_p)}$ \cite{Jacq2}, \cite[Prop. 4.7]{SSS}, the analogue of \cite[Lem. 3.16]{DingGL3} for $\mathrm{GSp}_4$, and an argument similar to \cite[Prop. 3.3]{Ding}.
\qed
\begin{lem}\label{4.4}
    For any $\mathscr C(I,s_i)$, we have:

(i) $\mathrm{dim}_E\mathrm{Ext}^1(\pi_{\textrm{alg}}(\phi,\boldsymbol h),\mathscr C(I,s_i))=\mathrm{dim}_E\mathrm{Ext}^1(\mathscr C(I,s_i),\pi_{\textrm{alg}}(\phi,\boldsymbol h))=1$.

(ii) Let $\tilde\pi_{\textrm{alg}}(\phi,\boldsymbol h)\in\mathrm{Ext}^1(\pi_{\textrm{alg}}(\phi,\boldsymbol h),\pi_{\textrm{alg}}(\phi,\boldsymbol h))$ be non-split, then the pullback map (via $\tilde\pi_{\textrm{alg}}(\phi,\boldsymbol h)\twoheadrightarrow\pi_{\textrm{alg}}(\phi,\boldsymbol h)$)
$$\mathrm{Ext}^1(\pi_{\textrm{alg}}(\phi,\boldsymbol h),\mathscr C(I,s_i))\rightarrow\mathrm{Ext}^1(\tilde\pi_{\textrm{alg}}(\phi,\boldsymbol h),\mathscr C(I,s_i))$$ is a bijection.
\end{lem}
\proof
(i) The first equality follows from Schraen’s spectral sequence\cite[Cor. 4.9]{SSS}, \eqref{socle} and \cite[Lem. 2,26]{DingGLn}. The second equality is proven in \cite[Cor. 5.2.6]{breuil2023probleme}.

(ii) follows from \cite[Cor. 4.9]{SSS}, \eqref{socle} and \cite[Lem. 2,26]{DingGLn} and an argument similar to \cite[Lem. 3.5]{Ding}
\qed

For $w\in W$, consider the nature map
$$\mathrm{Hom}(T(\mathbb Q_p),E)\rightarrow\mathrm{Ext}^1(\mathrm{PS}(w(\phi),\boldsymbol h),\mathrm{PS}(w(\phi),\boldsymbol h)),$$$$\psi\mapsto(\mathrm{Ind}^{\mathrm{GSp}_4(\mathbb{Q}_p)}_{\overline B(\mathbb{Q}_p)}w(\phi)\eta\lambda(1+\epsilon\psi))^{\mathrm{an}}.$$
Composed with the pull-back map for an injection $\jmath:\pi_{\textrm{alg}}(\phi,\boldsymbol h)\hookrightarrow\mathrm{PS}_1(w(\phi),\boldsymbol h)$ and using \cite[Lem. 2.26]{DingGLn}, it induces a map
\begin{equation}
    \label{zetaw}
    \mathrm{Hom}(T(\mathbb{Q}_p),E)\rightarrow\mathrm{Ext}^1(\pi_{\textrm{alg}}(\phi,\boldsymbol h),\mathrm{PS}_1(w(\phi),\boldsymbol h)).
\end{equation}
Composed furthermore with the push-forward map via the injection $\mathrm{PS}_1(w(\phi),\boldsymbol h)\hookrightarrow\pi_1(\phi,\boldsymbol h)$ (associated to $\jmath$), we finally obtain a map
$$\zeta_w:\mathrm{Hom}(T(\mathbb{Q}_p),E)\rightarrow\mathrm{Ext}^1(\pi_{\textrm{alg}}(\phi,\boldsymbol h),\pi_1(\phi,\boldsymbol h)).$$
Note that the map $\zeta_w$ does not depend on the choice of $\jmath$.

\begin{prop}\label{4.5}
    (i) \eqref{zetaw} is bijective. 
    $\mathrm{dim}_E\mathrm{Ext}^1(\pi_{\textrm{alg}}(\phi,\boldsymbol h),\mathrm{PS}_1(w(\phi),\boldsymbol h))=6.$

    (ii) $\zeta_w\vert_{\mathrm{Hom}_{g'}(T(\mathbb{Q}_p),E)}$ is equal to the composition of the map of same name ,which defined in Prop. \ref{4.3}, with the injective push-forward map $\mathrm{Ext}^1(\pi_{\textrm{alg}}(\phi,\boldsymbol h),\pi_{\textrm{alg}}(\phi,\boldsymbol h))\hookrightarrow\mathrm{Ext}^1(\pi_{\textrm{alg}}(\phi,\boldsymbol h),\pi_1(\phi,\boldsymbol h))$.
\end{prop}
\proof
(i) follows from \cite[Cor. 4.9]{SSS} and \cite[Lem. 2,26]{DingGLn}. (ii) is clear (see the remark below).
\qed
\begin{rem}\label{4.6}
    The map $\zeta_w$ can also be obtained by using Emerton’s functor $I^{\mathrm{GSp}_4(\mathbb{Q}_p)}_{\overline B(\mathbb{Q}_p)}$. In fact, by definition (and using \cite[Lem. 2,26]{DingGLn}), it is not difficult to see for $\psi\in\mathrm{Hom}(T(\mathbb{Q}_p),E)$, $I^{\mathrm{GSp}_4(\mathbb{Q}_p)}_{\overline B(\mathbb{Q}_p)}(w(\phi)\eta\lambda(1+\psi\epsilon))\subset(\mathrm{Ind}^{\mathrm{GSp}_4(\mathbb{Q}_p)}_{\overline B(\mathbb{Q}_p)}w(\phi)\eta\lambda(1+\epsilon\psi))^{\mathrm{an}}$ is an extension of $\pi_{\textrm{alg}}(\phi,\boldsymbol h)$ by a certain subrepresentation $V$ of  $\mathrm{PS}_1(w(\phi),\boldsymbol h)$. Then $\zeta_w(\psi)$ is the image of $I^{\mathrm{GSp}_4(\mathbb{Q}_p)}_{\overline B(\mathbb{Q}_p)}(w(\phi)\eta\lambda(1+\psi\epsilon))$ via the push-forward map $V\hookrightarrow \mathrm{PS}_1(w(\phi),\boldsymbol h)\hookrightarrow\pi_1(\phi,\boldsymbol h)$.
\end{rem}
\begin{prop}\label{4.7}
    (i) We have exact sequences
     $$0\rightarrow \mathrm{Ext}^1(\pi_{\textrm{alg}}(\phi,\boldsymbol h),\pi_{\textrm{alg}}(\phi,\boldsymbol h))\rightarrow\mathrm{Ext}^1(\pi_{\textrm{alg}}(\phi,\boldsymbol h),\mathrm{PS}_1(w(\phi),\boldsymbol h))$$
     $$\rightarrow\mathrm{Ext}^1(\pi_{\textrm{alg}}(\phi,\boldsymbol h),\mathscr C(w,s_i))\oplus\mathrm{Ext}^1(\pi_{\textrm{alg}}(\phi,\boldsymbol h),\mathscr C(w,s_i))\rightarrow 0$$
     for $w\in W$, and
    $$0\rightarrow \mathrm{Ext}^1(\pi_{\textrm{alg}}(\phi,\boldsymbol h),\pi_{\textrm{alg}}(\phi,\boldsymbol h))\rightarrow\mathrm{Ext}^1(\pi_{\textrm{alg}}(\phi,\boldsymbol h),\pi_1(\phi,\boldsymbol h))$$
    $$\rightarrow\oplus_{\substack{i=1,2,I\subset\{1,2,3,4\},\\\#I=i,\Sigma_{j\in I}j\neq 5}}\mathrm{Ext}^1(\pi_{\textrm{alg}}(\phi,\boldsymbol h),\mathscr C(I,s_i))\rightarrow 0.$$
    $\mathrm{dim}_E\mathrm{Ext}^1(\pi_{\textrm{alg}}(\phi,\boldsymbol h),\pi_1(\phi,\boldsymbol h))=12.$
    
    (ii) The nature map $$\oplus_{w\in W}\mathrm{Ext}^1(\pi_{\textrm{alg}}(\phi,\boldsymbol h),\mathrm{PS}_1(w(\phi),\boldsymbol h))\rightarrow\mathrm{Ext}^1(\pi_{\textrm{alg}}(\phi,\boldsymbol h),\pi_1(\phi,\boldsymbol h))$$ is surjective.
\end{prop}
    \proof
 It follows from Lem. \ref{4.4}(i), Prop. \ref{4.5}(i) and Prop. \ref{4.3}(i) and an argument similar to \cite[Prop. 3.8]{Ding}
    \qed

Respectively, denote by $\mathrm{Ext}^1_g(\pi_{\textrm{alg}}(\phi,\boldsymbol h),\pi_1(\phi,\boldsymbol h))$, $\mathrm{Ext}^1_{g'}(\pi_{\textrm{alg}}(\phi,\boldsymbol h),\pi_1(\phi,\boldsymbol h))$ and $\mathrm{Ext}^1_w(\pi_{\textrm{alg}}(\phi,\boldsymbol h),\pi_1(\phi,\boldsymbol h))$ the image of $\mathrm{Ext}^1_{\textrm{lalg}}(\pi_{\textrm{alg}}(\phi,\boldsymbol h),\pi_{\textrm{alg}}(\phi,\boldsymbol h))$, $\mathrm{Ext}^1(\pi_{\textrm{alg}}(\phi,\boldsymbol h),\allowbreak \pi_{\textrm{alg}}(\phi,\boldsymbol h))$ and $\mathrm{Ext}^1(\pi_{\textrm{alg}}(\phi,\boldsymbol h),\mathrm{PS}_1(w(\phi),\boldsymbol h))$ in $\mathrm{Ext}^1(\pi_{\textrm{alg}}(\phi,\boldsymbol h),\pi_1(\phi,\boldsymbol h))$, for $w\in W$.
We have hence an isomorphism 
\begin{equation}\label{zeta}
\zeta_w:\mathrm{Hom}(T(\mathbb{Q}_p),E)\mathop\rightarrow\limits^\sim\mathrm{Ext}^1_w(\pi_{\textrm{alg}}(\phi,\boldsymbol h),\pi_1(\phi,\boldsymbol h)).
\end{equation}
By Prop. \ref{4.3}(ii), the following square commutes:
\begin{equation}\label{zetas}
    \xymatrix{
    \mathrm{Hom}_{g'}(T(\mathbb{Q}_p),E)\ar[r]^-{\zeta_{w_1}}_-\sim\ar[d]_{w_2w_1^{-1}}^\sim&\mathrm{Ext}^1_{g'}(\pi_{\textrm{alg}}(\phi,\boldsymbol h),\pi_1(\phi,\boldsymbol h))\ar@{=}[d]\\
    \mathrm{Hom}_{g'}(T(\mathbb{Q}_p),E)\ar[r]^-{\zeta_{w_2}}_-\sim&\mathrm{Ext}^1_{g'}(\pi_{\textrm{alg}}(\phi,\boldsymbol h),\pi_1(\phi,\boldsymbol h))
    }.
\end{equation}

    \subsubsection{Parabolic inductions}
This subsection is an analogue of \cite[\S3.1.3]{Ding}. It may be a little more complicated in form because of the structure of parabolic subgroups of $\mathrm{GSp_4}$.

Through \eqref{TLP}, $\lambda$ is identified with $z^{(h_1-h_3,h_1-h_2)}\boxtimes z^{h_4}$. Let $\boldsymbol{h}^P_2=(h_1-h_3,h_1-h_2)\in \mathbb{Z}^2$ and $\boldsymbol{h}^P_1=h_4\in \mathbb{Z}$.
Assume $I_P\subset\{1,2,3,4\}$, $\#I_P=2$, $\Sigma_{j\in I_P}j\neq 5$. Choose the shorter $w_{I_P}\in W$ such that $w_{I_P}^{-1}(\{1,2\})=I_P$. (The only other choice is $s_1w_{I_P}$, which is longer.) Consider the parabolic induction
\begin{equation}\label{IndP}
    (\mathrm{Ind}_{\overline{P}(\mathbb{Q}_p)}^{\mathrm{GSp}_4(\mathbb{Q}_p)}(\pi_1((w_{I_P}(\phi))_2^P,\boldsymbol{h}_2^P)\boxtimes (w_{I_P}(\phi))_1^Pz^{\boldsymbol h_1^P})(\varepsilon^{-1}\circ(\mathrm{det\boxtimes 1})))^{\mathrm{an}}
    \supset \pi_{\textrm{alg}}(\phi,\boldsymbol h).
\end{equation}
(Here, $\pi_1((w_{I_P}(\phi))_2^P,\boldsymbol{h}_2^P)$ is a representation of $\mathrm{GL}_2(\mathbb{Q}_p)$ defined in \cite[\S3]{Ding}.)
\begin{lem}\label{4.8}
    For $i\in\{1,2\}$ and $I\subset\{1,2,3,4\}$, $\#I=i$, $\Sigma_{j\in I}j\neq 5$, $\mathscr C(I,s_i)$ appears as an irreducible constituent of \eqref{IndP} iff one of the following conditions holds:

    (1) $i=1$, $I\subsetneq I_P$;

    (2) $i=2$, $I=I_P$.

Moreover, each of such constituents has multiplicity one, and lies in the socle of $$ (\mathrm{Ind}_{\overline{P}(\mathbb{Q}_p)}^{\mathrm{GSp}_4(\mathbb{Q}_p)}(\pi_1((w_{I_P}(\phi))_2^P,\boldsymbol{h}_2^P)\boxtimes (w_{I_P}(\phi))_1^Pz^{\boldsymbol h_1^P})(\varepsilon^{-1}\circ(\mathrm{det\boxtimes 1})))^{\mathrm{an}}
/\pi_{\textrm{alg}}(\phi,\boldsymbol h).$$
\end{lem}
\proof
It follows from \cite{OS}, \cite[Prop. 4.1.2]{OS1} and an argument similar to \cite[Lem. 3.10]{Ding}
\qed

Denote by $S_{I_P}$ the subset of the constituents $\mathscr C(I,s_i)$ those that satisfy one of the conditions in Lem. \ref{4.8}. Then \eqref{IndP} contains a unique subrepresentation $\pi_{I_P}(\phi,\boldsymbol h)$, such that $\mathrm{soc}\pi_{I_P}(\phi,\boldsymbol h)=\oplus_{\mathscr C\in S_{I_P}}\mathscr C$. It is easy to see the (tautological injection) $\pi_{\textrm{alg}}(\phi,\boldsymbol h)\hookrightarrow\pi_{I_P}(\phi,\boldsymbol h)$ uniquely extends to $\pi_{I_P}(\phi,\boldsymbol h)\hookrightarrow\pi_1(\phi,\boldsymbol h)$.

Through \eqref{TLQ}, $\lambda$ is identified with $z^{h_1-h_3}\boxtimes z^{(h_3,h_4)}$. Let $\boldsymbol h^Q_1=h_1-h_3\in\mathbb{Z}$, $\boldsymbol h^Q_2=(h_3,h_4)\in\mathbb{Z}^2$. Assume $I_Q\subset\{1,2,3,4\}$, $\#I_P=1$. Choose the shorter $w_{I_Q}\in W$ satiesfy $w_{I_Q}^{-1}(\{1\})=I_Q$. (The only other choice of $s_2w_{I_Q}$, which is longer.) Consider the parabolic induction
\begin{equation}\label{IndQ}
    (\mathrm{Ind}_{\overline{Q}(\mathbb{Q}_p)}^{\mathrm{GSp}_4(\mathbb{Q}_p)}((w_{I_Q}(\phi))_1^Qz^{\boldsymbol h_1^Q}\boxtimes\pi_1((w_{I_Q}(\phi))_2^Q,\boldsymbol{h}_2^Q))\varepsilon^{-1}\circ(z^2\boxtimes 1))^{\mathrm{an}}
    \supset \pi_{\textrm{alg}}(\phi,\boldsymbol h).
\end{equation}
(Here, $\pi_1((w_{I_Q}(\phi))_2^Q,\boldsymbol{h}_2^Q)$ is a representation of $\mathrm{GL}_2(\mathbb{Q}_p)$ defined in \cite[\S3]{Ding}.)
Similarly, we have
\begin{lem}\label{4.9}
    For $i\in\{1,2\}$ and $I\subset\{1,2,3,4\}$, $\#I=i$, $\Sigma_{j\in I}j\neq 5$, $\mathscr C(I,s_i)$ appears as an irreducible constituent of \eqref{IndQ} iff one of the following conditions holds:

    (1) $i=2$, $I\supsetneq I_Q$;

    (2) $i=2$, $I=I_Q$.

Moreover, each of such constituents has multiplicity one, and lies in the socle of $$(\mathrm{Ind}_{\overline{Q}(\mathbb{Q}_p)}^{\mathrm{GSp}_4(\mathbb{Q}_p)}((w_{I_Q}(\phi))_1^Qz^{\boldsymbol h_1^Q}\boxtimes\pi_1((w_{I_Q}(\phi))_2^Q,\boldsymbol{h}_2^Q))\varepsilon^{-1}\circ(z^2\boxtimes 1))^{\mathrm{an}}
/\pi_{\textrm{alg}}(\phi,\boldsymbol h).$$
\end{lem}\qed

Similarly, denote by $S_{I_Q}$ the subset of the constituents $\mathscr C(I,s_i)$ those that satisfy one of the conditions in Lem. \ref{4.9}. Then \eqref{IndQ} contains a unique subrepresentation $\pi_{I_Q}(\phi,\boldsymbol h)$, such that $\mathrm{soc}\pi_{I_Q}(\phi,\boldsymbol h)=\oplus_{\mathscr C\in S_{I_Q}}\mathscr C$. The (tautological injection) $\pi_{\textrm{alg}}(\phi,\boldsymbol h)\hookrightarrow\pi_{I_Q}(\phi,\boldsymbol h)$ uniquely extends to $\pi_{I_Q}(\phi,\boldsymbol h)\hookrightarrow\pi_1(\phi,\boldsymbol h)$.

For $X\in\{P,Q\}$, we have:
\begin{prop}\label{4.10}
$$\mathrm{dim}_E\mathrm{Ext}^1(\pi_{\textrm{alg}}(\phi,\boldsymbol h),\pi_{I_X}(\phi,\boldsymbol h))=7.$$
And the following push-forward map is injective:
$$\mathrm{Ext}^1(\pi_{\textrm{alg}}(\phi,\boldsymbol h),\pi_{I_X}(\phi,\boldsymbol h))\rightarrow\mathrm{Ext}^1(\pi_{\textrm{alg}}(\phi,\boldsymbol h),\pi_1(\phi,\boldsymbol h)).$$
\end{prop}
\proof
It follows from Prop. \ref{4.7}(i), Lem. \ref{4.4}(i), Prop. \ref{4.3}(i) and an argument similar to \cite[Prop. 3.11]{Ding}
\qed

For $X\in\{P,Q\}$, set $\mathrm{Ext}^1_{I_X}(\pi_{\textrm{alg}}(\phi,\boldsymbol h),\pi_1(\phi,\boldsymbol h))$ to be the image of the map above. 

Note by \cite[Lem. 2.26]{DingGLn}, the following natural maps are bijective:
{\fontsize{8pt}{9.6pt}\selectfont\begin{equation}
    \label{P1}
    \begin{aligned}
    &\mathrm{Ext}^1(\pi_{\textrm{alg}}(\phi,\boldsymbol h),\pi_{I_P}(\phi,\boldsymbol h))\mathop\rightarrow\limits^{\sim}\\
    \mathrm{Ext}^1(\pi_{\textrm{alg}}(\phi,\boldsymbol h),
    (\mathrm{Ind}&_{\overline{P}(\mathbb{Q}_p)}^{\mathrm{GSp}_4(\mathbb{Q}_p)}(\pi_1((w_{I_P}(\phi))_2^P,\boldsymbol{h}_2^P)\boxtimes (w_{I_P}(\phi))_1^Pz^{\boldsymbol h_1^P})(\varepsilon^{-1}\circ(\mathrm{det\boxtimes 1})))^{\mathrm{an}}),
    \end{aligned}
\end{equation}}
and
{\fontsize{8pt}{9.6pt}\selectfont\begin{equation}
    \label{Q1}
     \begin{aligned}
    &\mathrm{Ext}^1(\pi_{\textrm{alg}}(\phi,\boldsymbol h),\pi_{I_Q}(\phi,\boldsymbol h))\mathop\rightarrow\limits^{\sim}\\
    \mathrm{Ext}^1(\pi_{\textrm{alg}}(\phi,\boldsymbol h),
    (\mathrm{Ind}&_{\overline{Q}(\mathbb{Q}_p)}^{\mathrm{GSp}_4(\mathbb{Q}_p)}((w_{I_Q}(\phi))_1^Qz^{\boldsymbol h_1^Q}\boxtimes\pi_1((w_{I_Q}(\phi))_2^Q,\boldsymbol{h}_2^Q))\varepsilon^{-1}\circ(z^2\boxtimes 1))^{\mathrm{an}}).
    \end{aligned}
\end{equation}}
By Schraen’s spectral sequence(\cite[Cor. 4.9]{SSS}), there are bijections
{\fontsize{8pt}{9.6pt}\selectfont\begin{equation}
    \label{P2}
    \begin{aligned}
    \mathrm{Ext}^1_{L_P(\mathbb{Q}_p)}(\pi_{\textrm{alg}}((w_{I_P}(\phi))_2^P,\boldsymbol{h}_2^P)\boxtimes (w_{I_P}(\phi))_1^Pz^{\boldsymbol h_1^P},\pi_1((w_{I_P}(\phi))_2^P,\boldsymbol{h}_2^P)\boxtimes (w_{I_P}(\phi))_1^Pz^{\boldsymbol h_1^P})\\\mathop\rightarrow\limits^\sim
    \mathrm{Ext}^1(\pi_{\textrm{alg}}(\phi,\boldsymbol h),
    (\mathrm{Ind}_{\overline{P}(\mathbb{Q}_p)}^{\mathrm{GSp}_4(\mathbb{Q}_p)}(\pi_1((w_{I_P}(\phi))_2^P,\boldsymbol{h}_2^P)\boxtimes (w_{I_P}(\phi))_1^Pz^{\boldsymbol h_1^P})(\varepsilon^{-1}\circ(\mathrm{det\boxtimes 1})))^{\mathrm{an}})
    ,\end{aligned}
\end{equation}}
and
{\fontsize{8pt}{9.6pt}\selectfont\begin{equation}
    \label{Q2}
    \begin{aligned}
    \mathrm{Ext}^1_{L_Q(\mathbb{Q}_p)}((w_{I_Q}(\phi))_1^Qz^{\boldsymbol h_1^Q}\boxtimes\pi_{\textrm{alg}}((w_{I_Q}(\phi))_2^Q,\boldsymbol{h}_2^Q),(w_{I_Q}(\phi))_1^Qz^{\boldsymbol h_1^Q}\boxtimes\pi_1((w_{I_Q}(\phi))_2^Q,\boldsymbol{h}_2^Q))\\\mathop\rightarrow\limits^\sim
    \mathrm{Ext}^1(\pi_{\textrm{alg}}(\phi,\boldsymbol h),
    (\mathrm{Ind}_{\overline{Q}(\mathbb{Q}_p)}^{\mathrm{GSp}_4(\mathbb{Q}_p)}((w_{I_Q}(\phi))_1^Qz^{\boldsymbol h_1^Q}\boxtimes\pi_1((w_{I_Q}(\phi))_2^Q,\boldsymbol{h}_2^Q))\varepsilon^{-1}\circ(z^2\boxtimes 1))^{\mathrm{an}})
    .\end{aligned}
\end{equation}}
(Here, $\pi_1((w_{I_P}(\phi))_2^P,\boldsymbol{h}_2^P)$ and $\pi_{\textrm{alg}}((w_{I_P}(\phi))_2^P,\boldsymbol{h}_2^P)$ are representations of $\mathrm{GL}_2(\mathbb{Q}_p)$ defined in \cite[\S3]{Ding}.)
Using \eqref{P1},\eqref{P2},\eqref{Q1},\eqref{Q2}, the following natural maps
{\fontsize{8pt}{9.6pt}\selectfont\begin{equation}
    \label{P3}
    \begin{aligned}
    \mathrm{Ext}^1_{\mathrm{GL}_2(\mathbb{Q}_p)}&(\pi_{\textrm{alg}}((w_{I_P}(\phi))_2^P,\boldsymbol{h}_2^P),\pi_1((w_{I_P}(\phi))_2^P,\boldsymbol{h}_2^P))
    \\&\times\\\mathrm{Ext}^1_{\mathrm{GL}_1(\mathbb{Q}_p)}&((w_{I_P}(\phi))_1^Pz^{\boldsymbol h_1^P},(w_{I_P}(\phi))_1^Pz^{\boldsymbol h_1^P})
    \\
    \rightarrow\mathrm{Ext}^1_{L_P(\mathbb{Q}_p)}(\pi_{\textrm{alg}}((w_{I_P}(\phi))_2^P,\boldsymbol{h}_2^P)&\boxtimes (w_{I_P}(\phi))_1^Pz^{\boldsymbol h_1^P},\pi_1((w_{I_P}(\phi))_2^P,\boldsymbol{h}_2^P)\boxtimes (w_{I_P}(\phi))_1^Pz^{\boldsymbol h_1^P})
    \end{aligned}
\end{equation}}
and
{\fontsize{8pt}{9.6pt}\selectfont\begin{equation}
    \label{Q3}
    \begin{aligned}
    \mathrm{Ext}^1_{\mathrm{GL}_1(\mathbb{Q}_p)}&((w_{I_Q}(\phi))_1^Qz^{\boldsymbol h_1^Q},(w_{I_Q}(\phi))_1^Qz^{\boldsymbol h_1^Q})
    \\&\times\\\mathrm{Ext}^1_{\mathrm{GL}_2(\mathbb{Q}_p)}&(\pi_{\textrm{alg}}((w_{I_Q}(\phi))_2^Q,\boldsymbol{h}_2^Q),\pi_1((w_{I_Q}(\phi))_2^Q,\boldsymbol{h}_2^Q))
    \\
    \rightarrow\mathrm{Ext}^1_{L_Q(\mathbb{Q}_p)}((w_{I_Q}(\phi))_1^Qz^{\boldsymbol h_1^Q}&\boxtimes\pi_{\textrm{alg}}((w_{I_Q}(\phi))_2^Q,\boldsymbol{h}_2^Q),(w_{I_Q}(\phi))_1^Qz^{\boldsymbol h_1^Q}\boxtimes\pi_1((w_{I_Q}(\phi))_2^Q,\boldsymbol{h}_2^Q))
    \end{aligned}
\end{equation}}
sending a pair of representations to there tensor product (over $E[\epsilon]/(\epsilon^2)$), we finally obtain two maps
{\fontsize{8pt}{9.6pt}\selectfont\begin{equation}
    \label{zetaP}
    \zeta_{I_P}:\begin{aligned}
    &\mathrm{Ext}^1_{\mathrm{GL}_2(\mathbb{Q}_p)}(\pi_{\textrm{alg}}((w_{I_P}(\phi))_2^P,\boldsymbol{h}_2^P),\pi_1((w_{I_P}(\phi))_2^P,\boldsymbol{h}_2^P))
    \\\times&\mathrm{Ext}^1_{\mathrm{GL}_1(\mathbb{Q}_p)}((w_{I_P}(\phi))_1^Pz^{\boldsymbol h_1^P},(w_{I_P}(\phi))_1^Pz^{\boldsymbol h_1^P})
    \end{aligned}\rightarrow\mathrm{Ext}^1_{I_P}(\pi_{\textrm{alg}}(\phi,\boldsymbol h),\pi_1(\phi,\boldsymbol h))    
\end{equation}}
and
{\fontsize{8pt}{9.6pt}\selectfont\begin{equation}
    \label{zetaQ}
    \zeta_{I_Q}:\begin{aligned}
    &\mathrm{Ext}^1_{\mathrm{GL}_1(\mathbb{Q}_p)}((w_{I_Q}(\phi))_1^Qz^{\boldsymbol h_1^Q},(w_{I_Q}(\phi))_1^Qz^{\boldsymbol h_1^Q})
    \\\times&\mathrm{Ext}^1_{\mathrm{GL}_2(\mathbb{Q}_p)}(\pi_{\textrm{alg}}((w_{I_Q}(\phi))_2^Q,\boldsymbol{h}_2^Q),\pi_1((w_{I_Q}(\phi))_2^Q,\boldsymbol{h}_2^Q))
    \end{aligned}\rightarrow\mathrm{Ext}^1_{I_Q}(\pi_{\textrm{alg}}(\phi,\boldsymbol h),\pi_1(\phi,\boldsymbol h)).
\end{equation}}  

Let $\delta_X:=\begin{cases}1,X=P\\2,X=Q\end{cases}$. For $X\in\{P,Q\}$, $w\in\{w_{I_X},s_{\delta_X}w_{I_X}\}$, it's clear that $\mathrm{PS}_1(w(\phi),\boldsymbol h)$ is a subrepresentation of $\pi_{I_X}(\phi,\boldsymbol h)$ (by comparing constituents and using Lem. \ref{4.4}(i)), hence $\mathrm{Ext}^1_w(\pi_{\textrm{alg}}(\phi,\boldsymbol h),\pi_1(\phi,\boldsymbol h))\hookrightarrow\mathrm{Ext}^1_{I_X}(\pi_{\textrm{alg}}(\phi,\boldsymbol h),\pi_1(\phi,\boldsymbol h))$.
\begin{prop}\label{4.11}
    $\zeta_P$, $\zeta_Q$ are bijective. Moreover, the following two diagrams are commutative (where the left columns come from \cite[\S3.1.2]{Ding}):
{\fontsize{8pt}{9.6pt}\selectfont$$\xymatrix{
\genfrac{}{}{0pt}{0}{\mathrm{Ext}^1_{\textrm{id}}(\pi_{\textrm{alg}}((w_{I_P}(\phi))_2^P,\boldsymbol{h}_2^P),\pi_1((w_{I_P}(\phi))_2^P,\boldsymbol{h}_2^P))}{\times\mathrm{Ext}^1_{\mathrm{GL}_1(\mathbb{Q}_p)}((w_{I_P}(\phi))_1^Pz^{\boldsymbol h_1^P},(w_{I_P}(\phi))_1^Pz^{\boldsymbol h_1^P})}\ar[r]&\mathrm{Ext}^1_{w_{I_P}}(\pi_{\textrm{alg}}(\phi,\boldsymbol h),\pi_1(\phi,\boldsymbol h))\\
\mathrm{Hom}(\mathbb{G}_m^2(\mathbb{Q}_p),E)\times\mathrm{Hom}(\mathbb{G}_m(\mathbb{Q}_p),E)\ar[d]_-\sim\ar[u]^-\sim\ar[r]_-{\sim}^-{\eqref{TLP}}&\mathrm{Hom}(T(\mathbb{Q}_p),E)\ar[d]_-\sim\ar[u]^-\sim\\
\genfrac{}{}{0pt}{0}{\mathrm{Ext}^1_s(\pi_{\textrm{alg}}((w_{I_P}(\phi))_2^P,\boldsymbol{h}_2^P),\pi_1((w_{I_P}(\phi))_2^P,\boldsymbol{h}_2^P))}{\times\mathrm{Ext}^1_{\mathrm{GL}_1(\mathbb{Q}_p)}((w_{I_P}(\phi))_1^Pz^{\boldsymbol h_1^P},(w_{I_P}(\phi))_1^Pz^{\boldsymbol h_1^P})}\ar[r]&\mathrm{Ext}^1_{s_1w_{I_P}}(\pi_{\textrm{alg}}(\phi,\boldsymbol h),\pi_1(\phi,\boldsymbol h))
}$$}
and
{\fontsize{8pt}{9.6pt}\selectfont$$\xymatrix{
\genfrac{}{}{0pt}{0}{\mathrm{Ext}^1_{\mathrm{GL}_1(\mathbb{Q}_p)}((w_{I_Q}(\phi))_1^Qz^{\boldsymbol h_1^Q},(w_{I_Q}(\phi))_1^Qz^{\boldsymbol h_1^Q})}{\times\mathrm{Ext}^1_{\textrm{id}}(\pi_{\textrm{alg}}((w_{I_Q}(\phi))_2^Q,\boldsymbol{h}_2^Q),\pi_1((w_{I_Q}(\phi))_2^Q,\boldsymbol{h}_2^Q))}\ar[r]&\mathrm{Ext}^1_{w_{I_Q}}(\pi_{\textrm{alg}}(\phi,\boldsymbol h),\pi_1(\phi,\boldsymbol h))\\
\mathrm{Hom}(\mathbb{G}_m(\mathbb{Q}_p),E)\times\mathrm{Hom}(\mathbb{G}_m^2(\mathbb{Q}_p),E)\ar[d]_-\sim\ar[u]^-\sim\ar[r]_-{\sim}^-{\eqref{TLQ}}&\mathrm{Hom}(T(\mathbb{Q}_p),E)\ar[d]_-\sim\ar[u]^-\sim\\
\genfrac{}{}{0pt}{0}{\mathrm{Ext}^1_{\mathrm{GL}_1(\mathbb{Q}_p)}((w_{I_Q}(\phi))_1^Qz^{\boldsymbol h_1^Q},(w_{I_Q}(\phi))_1^Qz^{\boldsymbol h_1^Q})}{\times\mathrm{Ext}^1_s(\pi_{\textrm{alg}}((w_{I_Q}(\phi))_2^Q,\boldsymbol{h}_2^Q),\pi_1((w_{I_Q}(\phi))_2^Q,\boldsymbol{h}_2^Q))}\ar[r]&\mathrm{Ext}^1_{s_2w_{I_Q}}(\pi_{\textrm{alg}}(\phi,\boldsymbol h),\pi_1(\phi,\boldsymbol h))
}.$$}
\end{prop}
\proof
Similar to \cite[Prop. 3.12]{Ding}.
\begin{rem}\label{4.12}
    We see \eqref{P3}, \eqref{Q3} actually bijective.
\end{rem}
We denote by $\mathrm{Ext}^1_{I_P,g'}(\pi_{\textrm{alg}}(\phi,\boldsymbol h),\pi_1(\phi,\boldsymbol h))$ (resp.,$\mathrm{Ext}^1_{I_Q,g'}(\pi_{\textrm{alg}}(\phi,\boldsymbol h),\pi_1(\phi,\boldsymbol h))$) the image of {\fontsize{7.5pt}{9pt}\selectfont$$\mathrm{Ext}^1_{g'}(\pi_{\textrm{alg}}((w_{I_P}(\phi))_2^P,\boldsymbol{h}_2^P),\pi_1((w_{I_P}(\phi))_2^P,\boldsymbol{h}_2^P))\times\mathrm{Ext}^1_{\mathrm{GL}_1(\mathbb{Q}_p)}((w_{I_P}(\phi))_1^Pz^{\boldsymbol h_1^P},(w_{I_P}(\phi))_1^Pz^{\boldsymbol h_1^P})$$} (resp.,{\fontsize{7.5pt}{9pt}\selectfont$$\mathrm{Ext}^1_{\mathrm{GL}_1(\mathbb{Q}_p)}((w_{I_Q}(\phi))_1^Qz^{\boldsymbol h_1^Q},(w_{I_Q}(\phi))_1^Qz^{\boldsymbol h_1^Q})\times\mathrm{Ext}^1_{g'}(\pi_{\textrm{alg}}((w_{I_Q}(\phi))_2^Q,\boldsymbol{h}_2^Q),\pi_1((w_{I_Q}(\phi))_2^Q,\boldsymbol{h}_2^Q))$$}) via $\zeta_P$ (resp., $\zeta_Q$). We denote by $\mathrm{Hom}_{P,g'}(T(\mathbb{Q}_p),E)$ (resp. $\mathrm{Hom}_{Q,g'}(T(\mathbb{Q}_p),E)$) the image of $\mathrm{Hom}_{g'}(\mathbb{G}_m^2(\mathbb{Q}_p),E)\times\mathrm{Hom}(\mathbb{G}_m(\mathbb{Q}_p),E)$ (resp., $\mathrm{Hom}(\mathbb{G}_m(\mathbb{Q}_p),E)\times\mathrm{Hom}_{g'}\allowbreak(\mathbb{G}_m^2(\mathbb{Q}_p),E)$) via the middle horizontal arrow in the first (resp., second) diagram in Prop. \ref{4.11}. For the definition of those subspaces with subscript $g'$ on $\mathrm{GL}_2$ side, see \cite[\S3.1.2]{Ding}.
\begin{prop}\label{4.13}
For $X\in\{P,Q\}$, we have

(i) $\mathrm{dim}_E\mathrm{Ext}^1_{I_X,g'}(\pi_{\textrm{alg}}(\phi,\boldsymbol h),\pi_1(\phi,\boldsymbol h))=5$,

(ii) {\fontsize{7.8pt}{9.36pt}\selectfont$$\mathrm{Ext}^1_{I_X,g'}(\pi_{\textrm{alg}}(\phi,\boldsymbol h),\pi_1(\phi,\boldsymbol h))=\mathrm{Ext}^1_{w_{I_X}}(\pi_{\textrm{alg}}(\phi,\boldsymbol h),\pi_1(\phi,\boldsymbol h))\cap\mathrm{Ext}^1_{s_{\delta_X}w_{I_X}}(\pi_{\textrm{alg}}(\phi,\boldsymbol h),\pi_1(\phi,\boldsymbol h)),$$}

(iii) $$\xymatrix{
\mathrm{Hom}_{X,g'}(T(\mathbb{Q}_p),E)\ar[r]_-\sim^-{\zeta_{w_{I_X}}}\ar[d]
^-{s_{\delta_X}}_-\sim&\mathrm{Ext}^1_{I_X,g'}(\pi_{\textrm{alg}}(\phi,\boldsymbol h),\pi_1(\phi,\boldsymbol h))\ar@{=}[d]\\
\mathrm{Hom}_{X,g'}(T(\mathbb{Q}_p),E)\ar[r]_-\sim^-{\zeta_{s_{\delta_X}w_{I_X}}}&\mathrm{Ext}^1_{I_X,g'}(\pi_{\textrm{alg}}(\phi,\boldsymbol h),\pi_1(\phi,\boldsymbol h))\\
}$$commutes.
\end{prop}
\proof
(i) follows from \cite[Prop. 3.3(1)]{Ding}. (ii) follows from (i), Prop. \ref{4.10}, Prop. \ref{4.5}(i) and dimension comparison. (iii) follows from \cite[Prop. 3.3(2)]{Ding}, proposition 4.11 and an argument similar to \cite[(3.30)]{Ding}.
\qed

    \subsection{Hodge parameters in $\mathrm{GSp}_4(\mathbb{Q}_p)$-representations}
    \subsubsection{Construction and properties}
Define $\ell: \mathrm{Hom}(\mathbb{Q}_p^\times,\mathfrak{t}(E))\rightarrow \mathrm{Hom}(T(\mathbb{Q}_p),E)$ as the unique bijection such that $\mathrm{exp}_E\ell=\mathscr L\mathrm{exp}_{\mathfrak t(E)}$. It's easy to see $\ell w$=$\check w\ell$, for any $w\in W$. The image of $\mathrm{Hom}_*(\mathbb{Q}_p^\times,\mathfrak{t}(E))$ via $\ell$ is exactly $\mathrm{Hom}_*(T(\mathbb{Q}_p),E)$, for $*=sm,g,g'$. The image of $\mathrm{Hom}_{X,g'}(\mathbb{Q}_p^\times,\mathfrak{t}(E))$ via $\ell$ is exactly $\mathrm{Hom}_{Y,g'}(T(\mathbb{Q}_p),E)$, for $\{X,Y\}=\{P,Q\}$.

Assume $\boldsymbol{\alpha}=(\alpha_1,\alpha_2,\alpha_3,\alpha_4)\in (E^{\times})^4$ satiesfies $\alpha_1\alpha_4=\alpha_2\alpha_3$. Assume $D\in\mathrm{GSp}_4\text{-}\Phi\Gamma_{\mathrm{nc}}(\boldsymbol\alpha,\boldsymbol h)$. Assume $\phi=\mathscr L(\mathrm{unr}(\boldsymbol\alpha))$.

\begin{thm}\label{4.14}
    There is a unique surjective map $t_D:\mathrm{Ext}^1(\pi_{\textrm{alg}}(\phi,\boldsymbol h),\pi_1(\phi,\boldsymbol h))\twoheadrightarrow\overline{\mathrm{Ext}}^G(D,D)$ making the following diagram commute:
     $$\xymatrix{
    \bigoplus\limits_{w\in W}\overline{\mathrm{Ext}}^G_{\check w}(D,D)
    \ar[d]_{\boxplus_{w\in W}\zeta_w\ell\kappa_{\check w}}^-{\sim}\ar@{->>}[r]
    &\overline{\mathrm{Ext}}^G(D,D)
    \\
    \bigoplus\limits_{w\in W}\mathrm{Ext}^1_w(\pi_{\textrm{alg}}(\phi,\boldsymbol h),\pi_1(\phi,\boldsymbol h))
    \ar@{->>}[r]
    &\mathrm{Ext}^1(\pi_{\textrm{alg}}(\phi,\boldsymbol h),\pi_1(\phi,\boldsymbol h))
    \ar@{-->>}^{t_D}[u]
    }$$
\end{thm}
\proof It's enough two prove that if $(x_w)_w\in \mathrm{ker}:\bigoplus\limits_{w\in W}\mathrm{Ext}^1_w(\pi_{\textrm{alg}}(\phi,\boldsymbol h),\pi_1(\phi,\boldsymbol h))\rightarrow\mathrm{Ext}^1(\pi_{\textrm{alg}}(\phi,\boldsymbol h),\pi_1(\phi,\boldsymbol h))$, then, along clockwise direction, it maps to $0$ in $\overline{\mathrm{Ext}}^G(D,D)$.
By Prop. \ref{4.7}(i), comparing constituents and Prop. \ref{4.13}(ii), it's easy to see that $\mathrm{ker}:\bigoplus\limits_{w\in W}\mathrm{Ext}^1_w(\pi_{\textrm{alg}}(\phi,\boldsymbol h),\pi_1(\phi,\boldsymbol h))\rightarrow\mathrm{Ext}^1(\pi_{\textrm{alg}}(\phi,\boldsymbol h),\pi_1(\phi,\boldsymbol h))$ generates by elements in the form:
$x-y$, where there exist $X\in\{P,Q\}$ and $I_X\subset\{1,2,3,4\}$, such that $\#I_X=\delta_X$, $\Sigma_{i\in I_X}\neq 5$, $x_{w_{I_X}}=y_{s_{\delta_X}w_{I_X}}\in \mathrm{Ext}^1_{I_X}(\pi_{\textrm{alg}}(\phi,\boldsymbol h),\pi_1(\phi,\boldsymbol h))$ and all other components vanish.
Then $x-y$ maps to $0$ in $\overline{\mathrm{Ext}}^G(D,D)$ because of Prop. \ref{4.13}(iii), Prop. \ref{3.15}(ii) and properties of $\ell$.
\qed

We define $\mathcal{L}(D):=\mathrm{ker}t_D$. By comparing dimensions (Cor. \ref{3.9}, Cor. \ref{3.13}, Prop. \ref{4.7}(i)), $\mathrm{dim}_E\mathcal{L}(D)=2$. We define $\pi_{\textrm{min}}(D)$ to be the extention of $\mathcal{L}(D)\otimes_E\pi_{\textrm{alg}}(\phi,\boldsymbol h)$ by $\pi_1(\phi,\boldsymbol h)$ associated to $\mathcal{L}(D)$. The following lemma is clear.
\begin{lem}\label{4.15}
    $\mathcal{L}(D)\cap\mathrm{Ext}^1_w(\pi_{\textrm{alg}}(\phi,\boldsymbol h),\pi_1(\phi,\boldsymbol h))=0$, for any $w\in W$.
\end{lem}\qed

We have an exact sequence
\begin{equation}\label{4.a}\begin{aligned}
    0\rightarrow\mathrm{Hom}(\pi_{\textrm{alg}}(\phi,\boldsymbol h),\mathcal{L}(D)\otimes_E\pi_{\textrm{alg}}(\phi,\boldsymbol h))\rightarrow\mathrm{Ext}^1&(\pi_{\textrm{alg}}(\phi,\boldsymbol h),\pi_1(\phi,\boldsymbol h))\\&\rightarrow\mathrm{Ext}^1(\pi_{\textrm{alg}}(\phi,\boldsymbol h),\pi_{\textrm{min}}(D)).
\end{aligned}\end{equation}
By Lem. \ref{4.4}(ii), the last map of \eqref{4.a} is surjective. For $*=g,g',w,I_P,I_Q,(I_P,g'),\allowbreak (I_Q,g')$, denote by $\mathrm{Ext}^1_*(\pi_{\textrm{alg}}(\phi,\boldsymbol h),\pi_{\textrm{min}}(D))$ the image of $\mathrm{Ext}^1_*(\pi_{\textrm{alg}}(\phi,\boldsymbol h),\pi_1(\phi,\boldsymbol h))$ via the last map of \eqref{4.a}. 
\begin{cor}\label{4.16}
    The map $t_D$ induces $\mathrm{Ext}^1_g(\pi_{\textrm{alg}}(\phi,\boldsymbol h),\pi_{\textrm{min}}(D))\simeq\overline{\mathrm{Ext}}^G_g(D,D)$ and $\mathrm{Ext}^1_w(\pi_{\textrm{alg}}(\phi,\boldsymbol h),\pi_{\textrm{min}}(D))\simeq\overline{\mathrm{Ext}}^G_w(D,D)$ for $w\in W$.
\end{cor}
\proof By Lem. \ref{4.15}, $\mathrm{Ext}^1_*(\pi_{\textrm{alg}}(\phi,\boldsymbol h),\pi_1(D))\simeq\mathrm{Ext}^1_g(\pi_{\textrm{alg}}(\phi,\boldsymbol h),\pi_{\textrm{min}}(D))$. The corollary then follows from the definition of $t_D$, \eqref{kappa}, \eqref{zeta}, Prop. \ref{4.3}(i).\qed

The following proposition is then a direct consequence of Thm. \ref{4.14}.
\begin{cor}\label{4.17}
    $\pi_{\textrm{min}}(D)$ is the unique extension of $\pi_{\textrm{alg}}(\phi,\boldsymbol h)^{\oplus 2}$ by $\pi_1(\phi,\boldsymbol h)$ satisfying the following properties:

    $\bullet \mathrm{soc}\pi_{\textrm{min}}(D)=\pi_{\textrm{alg}}(\phi,\boldsymbol h)$ and $\mathrm{soc}(\pi_{\textrm{min}}(D)/\pi_{\textrm{alg}}(\phi,\boldsymbol h))=\mathrm{soc}(\pi_1(\phi,\boldsymbol h)/\pi_{\textrm{alg}}(\phi,\boldsymbol h))$.

    $\bullet$ There is a bijection $t_D:\overline{\mathrm{Ext}}^G(D,D)\rightarrow \mathrm{Ext}^1(\pi_{\textrm{alg}}(\phi,\boldsymbol h),\pi_{\textrm{min}}(D))$  which is compatible with trianguline deformations: for any $w\in W$, we have
    $$\xymatrix{\overline{\mathrm{Ext}}^G_{\check{w}}(D,D)\ar[r]^-{t_D}_-\sim\ar[d]_-{\kappa_{\check{w}}}^-\sim
    &\mathrm{Ext}^1_w(\pi_{\textrm{alg}}(\phi,\boldsymbol{h}),\pi_{\textrm{min}}(D))\\
    \mathrm{Hom}(\mathbb{Q}_p^\times,\mathfrak{t}(E))\ar[r]^-\ell_-\sim
    &\mathrm{Hom}(T(\mathbb{Q}_p),E) \ar[u]^-{\zeta_w}_-\sim}
    $$
    commutes.
\end{cor}
\qed
\begin{thm}\label{4.18}
    Up to isomorphism, $\pi_{\textrm{min}}(D)$ uniquely determines $D$.
\end{thm}
\proof Note that $\mathrm{End}_{\mathrm{GSp}_4(\mathbb{Q}_p)}(\pi_1(\phi,\boldsymbol h))=\mathrm{End}_{\mathrm{GSp}_4(\mathbb{Q}_p)}(\pi_{\textrm{alg}}(\phi,\boldsymbol h))=E$, so $\pi_{\textrm{min}}(D)$ determines $\mathcal{L}(D)$. Then the theorem follows from the definition of $t_D$ (Thm. \ref{4.14}) and Thm. \ref{3.17}.\qed

\subsubsection{Universal extensions}
In this subsection, we give a reformulation of Thm. \ref{4.14} using deformation rings of $(\varphi,\Gamma)$-modules with $\mathrm{GSp}_4$-structure. Denote by $R_{D,*}$ the deformation rings for $\mathrm{GSp}_4\textrm{-}X_{D,*}$, $*\in\{\emptyset,w,g\}$. For $\delta\in \mathrm{Hom}(T(\mathbb{Q}_p),E^\times)$, denote by $R_\delta$ the deformation ring of $\delta$ over local Artinian $E$-algebras. If $\gamma\in \mathrm{Hom}(\mathbb{Q}_p^\times,T(E))$, denote by $R_\gamma$ the deformation ring of $\gamma$. Moreover, if $\delta$(resp., $\gamma$) is locally algebraic, denote by $R_{\delta,g}$(resp., $R_{\gamma,g}$) the deformation ring of locally algebraic deformations of $\delta$(resp., $\gamma$).
For a complete local Noetherian $E$-algebra $R$, we use $\mathfrak{m}_R$ to denote its maximal ideal and we will use $\mathfrak{m}$ for simplicity when it does not cause confusion.

We have nature surjections $R_D\twoheadrightarrow R_{D,w}\twoheadrightarrow R_{D,g}$, $R_\delta\twoheadrightarrow R_{\delta,g}$, $R_\gamma\twoheadrightarrow R_{\gamma,g}$. We have natural isomorphisms of $E$-vector spaces for the tangent spaces $(\mathfrak{m}_{R_*}/\mathfrak{m}_{R_*}^2)^\vee\simeq\mathrm{Ext}^G_*(D,D)$ for $*\in\{\emptyset,g,w\}$, $(\mathfrak{m}_{R_\delta}/\mathfrak{m}_{R_\delta}^2)^\vee\simeq\mathrm{Hom}(T(\mathbb{Q}_p),E)$, $(\mathfrak{m}_{R_{\delta,g}}/\mathfrak{m}_{R_{\delta,g}}^2)^\vee\simeq\mathrm{Hom}_{\textrm{sm}}(T(\mathbb{Q}_p),E)$, $(\mathfrak{m}_{R_\gamma}/\mathfrak{m}_{R_\gamma}^2)^\vee\simeq\mathrm{Hom}(\mathbb{Q}_p^\times,\mathfrak{t}(E))$, $(\mathfrak{m}_{R_{\gamma,g}}/\mathfrak{m}_{R_{\gamma,g}}^2)^\vee\simeq\mathrm{Hom}_{\textrm{sm}}(\mathbb{Q}_p^\times,\allowbreak \mathfrak{t}(E))$. For $w\in W$, by Prop. \ref{3.11},Prop. \ref{3.12}(i), $\kappa_w$ induces a commutative Cartesian diagram (of local Artinian $E$-algebras):
\begin{equation}\label{RD}
\xymatrix{
R_{w(\mathrm{unr}(\boldsymbol\alpha))z^{\boldsymbol h}}/\mathfrak{m}^2\ar@{->>}[r]\ar@{^{(}->}[d]
&R_{w(\mathrm{unr}(\boldsymbol\alpha))z^{\boldsymbol h},g}/\mathfrak{m}^2\ar@{^{(}->}[d]\\
R_{D,w}/\mathfrak{m}^2\ar@{->>}[r]
&R_{D,g}/\mathfrak{m}^2
}.
\end{equation}
By Prop. \ref{3.12}(iii), the image of the right vertical map in \eqref{RD} is independent on $w\in W$ and we denote the image by $A_{D,g}\subset R_{D,g}/\mathfrak{m}^2$. Define $A_{D,*}:=R_{D,*}/\mathfrak{m}^2\times_{R_{D,g}/\mathfrak{m}^2}A_{D,g}$ for $*\in\{\emptyset,w\}$. We have $A_{D,w}\simeq R_{w(\mathrm{unr}(\boldsymbol\alpha))z^{\boldsymbol h}}/\mathfrak{m}^2$ by \eqref{4.a}. We have $\mathfrak{m}_{A_{D,*}}=\overline{\mathrm{Ext}}^G_*(D,D)^\vee$ for $*\in{\emptyset,g,w}$. Let $\mathcal{I}_w$ be the kernel of $A_D\twoheadrightarrow A_{D,w}$. By Thm. \ref{3.16}, $A_D\rightarrow\prod_w A_{D,w}$ is injective.

Denote by $\pi_1(\phi,\boldsymbol h)^{\textrm{univ}}_*$ be the tautological extension of $\mathrm{Ext}^1_*(\pi_{\textrm{alg}}(\phi,\boldsymbol h),\pi_1(\phi,\boldsymbol h))\otimes\allowbreak\pi_{\textrm{alg}}(\phi,\boldsymbol h)$ by $\pi_1(\phi,\boldsymbol h)$, $*\in\{\emptyset.w\}$. For $w\in W$, denote by $\delta_w:=w(\phi)\eta\lambda$, and $\widetilde{\delta}_w^{\textrm{univ}}$ the tautological extension of $\mathrm{Ext}^1_{T(\mathbb{Q}_p)}(\delta_w,\delta_w)\otimes\delta_w$($\simeq\mathrm{Hom}(T(\mathbb{Q}_p),E)\otimes\delta_w$) by $\delta_w$.
\begin{lem}\label{4.19}
    The induced representation $I_{\overline B(\mathbb{Q}_p)}^{\mathrm{GSp}_4(\mathbb{Q}_p)}\widetilde{\delta}_w^{\textrm{univ}}$ is the universal extension of $\mathrm{Ext}^1(\pi_{\textrm{alg}}(\phi,\boldsymbol h),\mathrm{PS}_1(w(\phi),\boldsymbol h))\otimes\pi_{\textrm{alg}}(\phi,\boldsymbol h)$ by $\mathrm{PS}_1(\phi,\boldsymbol h)$
\end{lem}
\proof
Use Rem. \ref{4.6}, Prop. \ref{4.5}(i), the surjectivity of the last map of \eqref{4.a}, and a similar argument of the proof of \cite[Lem. 3.36]{Ding}.\qed

There is an $A_{D,\check w}$-action on $\widetilde{\delta}_w^{\textrm{univ}}$ where an element $x\in \mathfrak{m}_{A_{D,\check w}}\simeq\overline{\mathrm{Ext}}^G_{\check w}(D,D)^\vee$ acts via $$\widetilde{\delta}_w^{\textrm{univ}}\twoheadrightarrow \mathrm{Hom}(T(\mathbb{Q}_p),E)\otimes\delta_w\xrightarrow{(\kappa_{\check w}^{-1}\ell^{-1})^\vee x}\delta_w\hookrightarrow\widetilde{\delta}_w^{\textrm{univ}}.$$
It induces an $A_{D,\check w}$-action in $I_{\overline B(\mathbb{Q}_p)}^{\mathrm{GSp}_4(\mathbb{Q}_p)}\widetilde{\delta}_w^{\textrm{univ}}$.  Similarly, there is an $A_{D,\check w}$-action on $\pi_1(\phi,\boldsymbol h)^{\textrm{univ}}_w$ where an element $x\in \mathfrak{m}_{A_{D,\check w}}\simeq\overline{\mathrm{Ext}}^G_{\check w}(D,D)^\vee$ acts via
$$\begin{aligned}\pi_1(\phi,\boldsymbol h)^{\textrm{univ}}_w\twoheadrightarrow\mathrm{Ext}^1_w(\pi_{\textrm{alg}}(\phi,\boldsymbol h),\pi_1(\phi,\boldsymbol h))&\otimes\pi_{\textrm{alg}}(\phi,\boldsymbol h)\\&\xrightarrow{(\kappa_{\check w}^{-1}\ell^{-1}\zeta_w^{-1})^\vee x}\pi_{\textrm{alg}}(\phi,\boldsymbol h)\hookrightarrow\pi_1(\phi,\boldsymbol h)^{\textrm{univ}}_w.\end{aligned}$$
By Lem. \ref{4.18}, we have an injection $I_{\overline B(\mathbb{Q}_p)}^{\mathrm{GSp}_4(\mathbb{Q}_p)}\widetilde{\delta}_w^{\textrm{univ}}\hookrightarrow\pi_1(\phi,\boldsymbol h)^{\textrm{univ}}_w$ of $\mathrm{GSp}_4(\mathbb{Q}_p)$-representations. It's easy to see this injection is also $A_{D,\check w}$-equivariant.

The following theorem is a reformulation of Thm. \ref{4.14}.
\begin{thm}\label{4.20}
    There is a unique $A_D$-action on $\pi_1(\phi,\boldsymbol h)^{\textrm{univ}}$ such that for all $w\in W$, we have an $A_{D,\check w}\times\mathrm{GSp}_4(\mathbb{Q}_p)$-equivariant injection $\pi_1(\phi,\boldsymbol h)^{\textrm{univ}}_w\hookrightarrow\pi_1(\phi,\boldsymbol h)^{\textrm{univ}}[\mathcal{I}_{\check w}]$. Moreover, this action is given by the following way. For $x\in \mathfrak{m}_{A_D}\simeq\overline{\mathrm{Ext}}^G(D,D)^\vee$ acts via
    $$\pi_1(\phi,\boldsymbol h)^{\textrm{univ}}\twoheadrightarrow\mathrm{Ext}^1(\pi_{\textrm{alg}}(\phi,\boldsymbol h),\pi_1(\phi,\boldsymbol h))\otimes\pi_{\textrm{alg}}(\phi,\boldsymbol h)\xrightarrow{t_D^\vee x}\otimes\pi_{\textrm{alg}}(\phi,\boldsymbol h)\hookrightarrow\pi_1(\phi,\boldsymbol h)^{\textrm{univ}}.$$
\end{thm}
\proof Notice that by the definition of $t_D$ (Thm. \ref{4.14}), $t_D$ coincides with $\kappa_{\check w}^{-1}\ell^{-1}\zeta_w^{-1}$ on $\mathrm{Ext}^1_w(\pi_{\textrm{alg}}(\phi,\boldsymbol h),\pi_1(\phi,\boldsymbol h))$. Hence the $A_D$-action described in the theorem satisfies the requirement of the theorem. The uniqueness follows from the fact that $\pi_1(\phi,\boldsymbol h)^{\textrm{univ}}$ is generated by $\pi_1(\phi,\boldsymbol h)^{\textrm{univ}}_w$, $w\in W$.\qed

By the construction of $\pi_{\textrm{min}}(D)$, we have:
\begin{cor}\label{4.21}
$\pi_{\textrm{min}}(D)\simeq\pi_1(\phi,\boldsymbol h)^{\textrm{univ}}[\mathfrak{m}_{A_D}]$.
\end{cor}\qed
\subsection{A slight variant of the construction}
In this subsection, we prove some slight variant of results in the last subsection, which will be useful in the next chapter.

Let $\mathrm{Ext}_U^G(D,D)$ be a certain subspace of $\mathrm{Ext}^G(D,D)$.  For $*\in\{w,g'\}$, $\mathrm{Ext}_{U,*}^G(D,D)\allowbreak:=\mathrm{Ext}_U^G(D,D)\cap\mathrm{Ext}_*^G(D,D)$. We assume the following hypotheses.
\begin{hyp}\label{4.22}
    (i) $\mathrm{Ext}_U^G(D,D)\cap\mathrm{Ext}_g^G(D,D)=0$.

    (ii) $\mathrm{dim}_E\mathrm{Ext}_{U,w}^G(D,D)=3$.
\end{hyp}
\begin{cor}\label{4.23}
    (i) The nature map $\oplus_{w\in W} \mathrm{Ext}_{U,w}^G(D,D)\rightarrow\mathrm{Ext}_U^G(D,D)$ is surjective.

    (ii) $\mathrm{dim}_E\mathrm{Ext}_U^G(D,D)=7$
\end{cor}
\proof Consider the following commutative diagram:
$$\xymatrix{
\oplus_w\mathrm{Ext}_{U,w}^G(D,D)\ar[d]^\sim\ar[r]
&\mathrm{Ext}_U^G(D,D)\ar@{^{(}->}[d]\\
\oplus_w\mathrm{Ext}_w^G(D,D)/\mathrm{Ext}_g^G(D,D)\ar@{->>}[r]
&\mathrm{Ext}^G(D,D)/\mathrm{Ext}_g^G(D,D)
}.$$
The vertical maps are injective by Hyp. \ref{4.22}(i) hence the left one is bijective by comparing dimensions (Hyp. \ref{4.22}(ii), Cor. \ref{3.9}), the surjectivity of the bottom map follows from Thm. \ref{3.16}. We deduce the top and right maps are also surjective. Hence (i) is proved and the right map is bijective. (ii) follows then from Cor. \ref{3.9} and Cor. \ref{3.13}.\qed

Denote by $\mathrm{Ext}^1_{U,w}(\pi_{\textrm{alg}}(\phi,\boldsymbol h),\pi_1(\phi,\boldsymbol h))$ the image of $\mathrm{Ext}_{U,\check{w}}^G(D,D)$ via $\zeta_w\ell\kappa_{\check w}$. Recalling that $\mathrm{ker}\kappa_{\check w}=\mathrm{Ext}_0^G(D,D)\subset \mathrm{Ext}_g^G(D,D)$, hence by Hyp. \ref{4.22}(i), $\zeta_w\ell\kappa_{\check w}$ induces a bijection  $\mathrm{Ext}_{U,\check{w}}^G(D,D)\simeq\mathrm{Ext}^1_{U,w}(\pi_{\textrm{alg}}(\phi,\boldsymbol h),\pi_1(\phi,\boldsymbol h))$. Denote by $\mathrm{Ext}^1_{U,g'}(\pi_{\textrm{alg}}(\phi,\boldsymbol h),\allowbreak \pi_1(\phi,\boldsymbol h))$ the image of $\mathrm{Ext}_{U,\check{w}}^G(D,D)$ via $\zeta_w\ell\kappa_{\check w}$. By Prop. \ref{3.14} and \eqref{zetas}, the definition of $\mathrm{Ext}_{U,\check{w}}^G(D,D)$ via $\zeta_w\ell\kappa_{\check w}$ isn't depend on the choice of $w$ and $\zeta_w\ell\kappa_{\check w}$ induces a bijection $\mathrm{Ext}_{U,g'}^G(D,D)\simeq\mathrm{Ext}^1_{U,g'}(\pi_{\textrm{alg}}(\phi,\boldsymbol h),\pi_1(\phi,\boldsymbol h))$ which isn't depend on the choice of $w\in W$, neither. Denote by $\mathrm{Ext}^1_U(\pi_{\textrm{alg}}(\phi,\boldsymbol h),\pi_1(\phi,\boldsymbol h))$ the subspace of $\mathrm{Ext}^1(\pi_{\textrm{alg}}(\phi,\boldsymbol h),\pi_1(\phi,\boldsymbol h))$ generated by $\mathrm{Ext}^1_{U,w}(\pi_{\textrm{alg}}(\phi,\boldsymbol h),\pi_1(\phi,\boldsymbol h))$, $w\in W$. We then have:

\begin{cor}\label{4.24}
    (i) $$0\rightarrow\mathrm{Ext}^1_{U,g'}(\pi_{\textrm{alg}}(\phi,\boldsymbol h),\pi_1(\phi,\boldsymbol h))\rightarrow\mathrm{Ext}^1_{U,w}(\pi_{\textrm{alg}}(\phi,\boldsymbol h),\pi_1(\phi,\boldsymbol h))$$
    $$\rightarrow\mathrm{Ext}^1(\pi_{\textrm{alg}}(\phi,\boldsymbol h),\mathscr C(w,s_i))\oplus\mathrm{Ext}^1(\pi_{\textrm{alg}}(\phi,\boldsymbol h),\mathscr C(w,s_i))\rightarrow 0$$
    for $w\in W$ are exact. $\mathrm{dim}_E\mathrm{Ext}_{U,g'}^G(D,D)=\mathrm{dim}_E\mathrm{ker}\mathrm{Ext}^1_{U,g'}(\pi_{\textrm{alg}}(\phi,\boldsymbol h),\pi_1(\phi,\boldsymbol h))=1$.

    (ii)$$0\rightarrow\mathrm{Ext}^1_{U,g'}(\pi_{\textrm{alg}}(\phi,\boldsymbol h),\pi_1(\phi,\boldsymbol h))\rightarrow\mathrm{Ext}^1_U(\pi_{\textrm{alg}}(\phi,\boldsymbol h),\pi_1(\phi,\boldsymbol h))$$
    $$\rightarrow\oplus_{\substack{i=1,2,I\subset\{1,2,3,4\},\\\#I=i,\Sigma_{j\in I}j\neq 5}}\mathrm{Ext}^1(\pi_{\textrm{alg}}(\phi,\boldsymbol h),\mathscr C(I,s_i))\rightarrow 0$$ is exact.
    $\mathrm{dim}_E\mathrm{ker}\mathrm{Ext}^1_U(\pi_{\textrm{alg}}(\phi,\boldsymbol h),\pi_1(\phi,\boldsymbol h))=9$.
\end{cor}
\proof The two sequences are left exact by Prop. \ref{4.7}(i). By Hyp. \ref{4.22}(i) and dimension comparing (Prop. \ref{4.3}(i)), we have $\mathrm{dim}_E\mathrm{ker}\mathrm{Ext}^1_{U,g'}(\pi_{\textrm{alg}}(\phi,\boldsymbol h),\pi_1(\phi,\boldsymbol h))\leq1$. Then (i) follows by dimension comparing (Lem. \ref{4.4}(i), Hyp. \ref{4.22}(ii)). Then the second map in the second sequence is surjective because the one in the first sequence does so. Hence the second sequence is exact and (ii) follows by dimension comparing ((i) and Lem. \ref{4.4}(i)).\qed

\begin{cor}\label{4.25}
    (i) There exists a unique surjection $$t_{D,U}:\mathrm{Ext}^1_U(\pi_{\textrm{alg}}(\phi,\boldsymbol h),\pi_1(\phi,\boldsymbol h))\rightarrow\mathrm{Ext}_U^G(D,D)$$ such that the following diagram commutes.
    $$\xymatrix{
    \bigoplus\limits_{w\in W}\mathrm{Ext}^G_{U,\check w}(D,D)
    \ar[d]_{\boxplus_{w\in W}\zeta_w\ell\kappa_{\check w}}^-{\sim}\ar@{->>}[r]
    &\mathrm{Ext}^G_U(D,D)
    \\
    \bigoplus\limits_{w\in W}\mathrm{Ext}^1_{U,w}(\pi_{\textrm{alg}}(\phi,\boldsymbol h),\pi_1(\phi,\boldsymbol h))
    \ar@{->>}[r]
    &\mathrm{Ext}^1_U(\pi_{\textrm{alg}}(\phi,\boldsymbol h),\pi_1(\phi,\boldsymbol h))
    \ar@{-->>}^{t_{D,U}}[u]
    }$$
    (ii) $\mathrm{ker}t_{D,U}=\mathrm{ker}t_D$.
\end{cor}
\proof 
First, if such $t_{D,U}$ exists, it's obvious it's unique and surjective. By the definition of $t_D$, the image of $\mathrm{Ext}^1_U(\pi_{\textrm{alg}}(\phi,\boldsymbol h),\pi_1(\phi,\boldsymbol h))$ via $t_D$ is in $\mathrm{Ext}^G_U(D,D)\hookrightarrow \overline{\mathrm{Ext}}^G(D,D)$. Hence by Thm. \ref{4.14}, one can take $t_{D,U}$ as the restriction of $t_D$ on $\mathrm{Ext}^1_U(\pi_{\textrm{alg}}(\phi,\boldsymbol h),\pi_1(\phi,\boldsymbol h))$. Hence (i) is proved and $\mathrm{ker}t_{D,U}\subset\mathrm{ker}t_D$. By Cor. \ref{4.24}(ii) and Cor. \ref{4.23}(ii), $\mathrm{dim}_E\mathrm{ker}t_{D,U}\geq 2=\mathrm{dim}_E\mathrm{ker}t_D$. (ii) is proved. \qed

Denote by $\pi_1(\phi,\boldsymbol h)^{\textrm{univ}}_{U,*}$ the tautological extension of $\mathrm{Ext}^1_{U,*}(\pi_{\textrm{alg}}(\phi,\boldsymbol h),\pi_1(\phi,\boldsymbol h))\allowbreak\otimes\pi_{\textrm{alg}}(\phi,\boldsymbol h)$ by $\pi_1(\phi,\boldsymbol h)$, $*\in\{\emptyset.w\}$. Let $\mathrm{Ext}^1_U(\delta_w,\delta_w)\subset\mathrm{Ext}^1(\delta_w,\delta_w)$($\simeq \mathrm{Hom}(T(\mathbb{Q}_p),E)$) be the image of $\mathrm{Ext}^G_{U,\check w}(D,D)$ via $\ell\kappa_{\check w}$. Let $\widetilde{\delta}_{U,w}^{\textrm{univ}}$ be the tautological extension of $\mathrm{Ext}^1_U(\delta_w,\delta_w)\otimes\delta_w$($\simeq\mathrm{Hom}(T(\mathbb{Q}_p),E)\otimes\delta_w$) by $\delta_w$. Let $A_{D,U,*}$ be the quotient of $A_{D,*}$ associated to $\mathrm{Ext}^G_{U,*}(D,D)\hookrightarrow \overline{\mathrm{Ext}}^G_*(D,D)$, $*\in\{\emptyset,w\}$. Let $A_{U,w}$ be the quotient of $R_{w(\mathrm{unr}(\boldsymbol \alpha)z^{\boldsymbol h}}/\mathfrak{m}^2$ associated to the image of $\mathrm{Ext}^G_{U,w}(D,D)$. The map $\kappa_w$ then induces an isomorphism of $E$-algebras $A_{U,w}\simeq A_{D,U,w}$ for $w\in W$. We equip $A_{U,\check w}$ with the $T(\mathbb{Q}_p)$-action given by $$T(\mathbb{Q}_p)\xrightarrow{\mathrm{univ}}R_{\delta_w}\xrightarrow{\sim} R_{\check{w}(\mathrm{unr}(\boldsymbol \alpha)z^{\boldsymbol h}}\rightarrow A_{U,\check w},$$where the first map is given by the universal extension and the middle map is induced by the composition of twisting by $\varepsilon^{-1}\circ p_1^2p_2$ and \eqref{L}. The $T(\mathbb{Q}_p)$-representation $A_{U,\check w}^\vee$ is then isomorphism to $\widetilde{\delta}_{U,w}^{\textrm{univ}}$.
Reciprocally, the $T(\mathbb{Q}_p)$-representation $\widetilde{\delta}_{U,w}^{\textrm{univ}}$ is equipped with a natural $A_{U,\check w}$-action (hence an $A_{D,U,\check w}$-action) as in the
discussion above Thm. \ref{4.20}. Note the natural map $E[T(\mathbb{Q}_p)]\rightarrow R_{\delta_w}/\mathfrak{m}^2$ is surjective. Thus the action of $T(\mathbb{Q}_p)$ and $A_{U,\check w}$ actually determine each other.

By a proof similar to Lem. \ref{4.19}, $I_{\overline B(\mathbb{Q}_p)}^{\mathrm{GSp}_4(\mathbb{Q}_p)}\widetilde{\delta}_{U,w}^{\textrm{univ}}$ is the universal extension of $\mathrm{Ext}^1_U(\pi_{\textrm{alg}}(\phi,\boldsymbol h),\mathrm{PS}_1(w(\phi),\boldsymbol h))\otimes\pi_{\textrm{alg}}(\phi,\boldsymbol h)$ by $\mathrm{PS}_1(\phi,\boldsymbol h)$. Moreover, similarly as in the discussion above Thm. \ref{4.20}, we have a $\mathrm{GSp}_4(\mathbb{Q}_p)\times A_{D,U,\check w}$-equivariant injection $I_{\overline B(\mathbb{Q}_p)}^{\mathrm{GSp}_4(\mathbb{Q}_p)}\widetilde{\delta}_{U,w}^{\textrm{univ}}\hookrightarrow\pi_1(\phi,\boldsymbol h)^{\textrm{univ}}_{U,w}$. Here the $A_{U,\check w}$($\simeq A_{D,U,\check w}$)-action is given in the similar way as in the discussion above Thm. \ref{4.20}. Set $\mathcal I_{U,w}:=\mathrm{ker}:A_{D,U}\twoheadrightarrow A_{D,U,w}$. Using Cor. \ref{4.25}, similarly to Thm. \ref{4.20} and Cor. \ref{4.21}, we have:
\begin{cor}\label{4.26}
    (i) There is a unique $A_{D,U}$-action on $\pi_1(\phi,\boldsymbol h)^{\textrm{univ}}_U$ such that we have an $A_{U,\check w}\times\mathrm{GSp}_4(\mathbb{Q}_p)$-equivariant injection $\pi_1(\phi,\boldsymbol h)^{\textrm{univ}}_{U,w}\hookrightarrow\pi_1(\phi,\boldsymbol h)^{\textrm{univ}}_U[\mathcal{I}_{U,\check w}]$ for all $w\in W$. Moreover, this action is given by the following way. For $x\in \mathfrak{m}_{A_{D,U}}\simeq\mathrm{Ext}^G_U(D,D)^\vee$ acts via
    $$\pi_1(\phi,\boldsymbol h)^{\textrm{univ}}_U\twoheadrightarrow\mathrm{Ext}^1_U(\pi_{\textrm{alg}}(\phi,\boldsymbol h),\pi_1(\phi,\boldsymbol h))\otimes\pi_{\textrm{alg}}(\phi,\boldsymbol h)\xrightarrow{t_{D,U}^\vee x}\pi_{\textrm{alg}}(\phi,\boldsymbol h)\hookrightarrow\pi_1(\phi,\boldsymbol h)^{\textrm{univ}}_U.$$

    (ii) $\pi_{\textrm{min}}(D)\simeq\pi_1(\phi,\boldsymbol h)^{\textrm{univ}}_U[\mathfrak{m}_{A_{D,U}}]$.
\end{cor}\qed

\section{Global theory: local-global compatibility}
\subsection{Classical theory: automorphic Galois representations}
Let $F\subset \mathbb{C}$ be an imaginary quadratic extension of $\mathbb{Q}$. We fix an isomorphism $\imath:\mathbb{C}\simeq\overline{\mathbb{Q}}_p$, and hence a prime $\wp$ of $F$ lying over $p$. For any $\mathbb{Q}$-algebra $R$, there is a conjugation on $R\otimes_\mathbb{Q}F$ induced by the conjugation on $F$. Consider the unitary general symplectic group $G/\mathbb{Q}$ attached to the quadratic extension $F$. More explicitly, $G$ is defined by
$$G(R)=\{M\in \mathrm{GL}_4(R\otimes_\mathbb{Q}F)\vert M^*M=I,\exists \mu(M)\in R\otimes_\mathbb{Q}F,s.t.\ M^TJM=\mu(M)J\}.$$ Notice that it implies $\mu(M)\overline{\mu(M)}=1$. We have $G\times_\mathbb{Q}F\simeq \mathrm{GSp}_4$ and $G(\mathbb{R})$ is compact. If $l=xx'$ is a prime splits in $F$, then $x$(resp.$x'$)$:F\hookrightarrow{\mathbb{Q}_l}$ induces an isomorphism $\iota_x$(resp.$\iota_{x'}$)$:G(\mathbb{Q}_l)\simeq\mathrm{GSp}_4(F_x)$(resp.$\mathrm{GSp}_4(F_{x'})$). The composition of $\iota_x\iota_{x'}^{-1}$ with $\mathrm{GSp}_4(F_x)\simeq\mathrm{GSp}_4(\mathbb{Q}_l)\simeq\mathrm{GSp}_4(F_{x'})$ equals taking the inverse transpose. The inclusion $F\subset \mathbb{C}$ induces an isomorphism $G(\mathbb{C})\simeq\mathrm{GSp}_4(\mathbb{C})$ and an isomorphism $\iota_\infty:G(\mathbb{R})\simeq\mathrm{GSp}_4(\mathbb{C})\cap \mathrm{U}_4(\mathbb{R})\triangleq\mathrm{UGSp}_4(\mathbb{R})$. We denote by $\mathbb{A}$ the ring of $\mathbb{Q}$-ad\`{e}les and $\mathbb{A}_f$ the ring of finite ad\`{e}les.

The space of automorphic forms of $G$ is the representation of $G(\mathbb{A})$ by right translations on the space $\mathcal{A}(G)$ of complex functions on  $G(\mathbb{Q})\backslash G(A)$ which are smooth and $G(\mathbb{R})$-finite. Denote $\mathrm{Irr}(G(\mathbb{R}))$ the set of isomorphic classes of irreducible continuous complex representations of $G(\mathbb{R})$. For $W\in \mathrm{Irr}(G(\mathbb{R}))$, we define $\mathcal{A}(G,W)$ to be the $G (\mathbb{A}_f)$ representation by right translations on the space of smooth vector valued functions $f:G(\mathbb{A}_f)\rightarrow W$ such that $f(\gamma g)=\gamma_\infty f(g)$ for all $g\in G(\mathbb{A}_f)$ and $\gamma\in G(\mathbb{Q})$. Lem. \ref{5.1} and Cor. \ref{5.2} follows by a same proof to \cite[Lem. 6.2.5, Cor. 6.2.7]{Selmergp}.

\begin{lem}\label{5.1}
    $\mathcal{A}(G)$ is admissible. $\mathcal{A}(G,W)$ is semi-simple for $W\in\mathrm{Irr}(G(\mathbb{R}))$. Thus $\mathcal{A}(G)=\bigoplus\limits_{W\in\mathrm{Irr}(G(\mathbb{R}))}W\otimes\mathcal{A}(G,W)$ is semi-simple.
\end{lem}\qed

An irreducible representation of $G(\mathbb{A})$ is said to be automorphic if it contributes to $\mathcal{A}(G)$.
\begin{cor}\label{5.2}
    If $\pi=\pi_\infty\otimes\pi_f$ is an automorphic representation of $G(\mathbb{A})$, then $\pi_f$ is defined over a number field.
\end{cor}\qed

Let $\boldsymbol{h}=(h_1,h_2,h_3,h_4)\in \mathbb Z^4$ such that $h_1>h_2>h_3>h_4$ and $h_1+h_4=h_2+h_3$. Then $\lambda:=p_1^{h_1-h_3-2}p_2^{h_1-h_2-1}p_3^{h_4}$ is a dominant weight of $T$. We define $W(\lambda)$ as the irreducible complex representation of $G(\mathbb{R})$ obtained by taking pullback of the irreducible representation of $\mathrm{GSp}_4(\mathbb{C})$  attached to the dominant weight $\lambda$ via $\iota_\infty$. Actually all irreducible representations of $G(\mathbb{R})$ can be obtained in this way. For a prime $l=xx'$ splits in $F$, denote by $K_l:=\iota_x^{-1}\mathrm{GSp}_4(\mathbb{Z}_l)$ and the Hecke operators $T_{l,i}\in\mathbb{C}[G(\mathbb Q_l)//K_l)]$ for $i=\{0,1,2\}$ are given by
$$T_{l,0}=[K_l\iota_x^{-1}\mathrm{diag}(l,l,l,l)K_l],$$
$$T_{l,1}=[K_l\iota_x^{-1}\mathrm{diag}(l,l,1,1)K_l],$$
$$T_{l,2}=[K_l\iota_x^{-1}\mathrm{diag}(l^2,l,l,1)K_l].$$
Recall that we have $\mathbb{C}[G(\mathbb Q_l)//K_l)]=\mathbb{C}[T_{l,0}^{\pm 1},T_{l,1},T_{l,2}]$.

From now on, we assume that $G$ satisfies the following conjecture.
\begin{conj}\label{5.3}
    If $\pi=\pi_\infty\otimes\pi_f$, then up to isomorphism, there exists a unique continuous semi-simple Galois representation with $\mathrm{GSp}_4$-structure $$\rho_\pi:\mathrm{Gal}_F\rightarrow\mathrm{GSp}_4(\overline{\mathbb{Q}}_p)$$ such that the following properties are satisfied:

    (Gal1) If prime $l=xx'$ splits in $F$ and $\iota_x^*\pi_l$ is unramified, then $\rho_\pi$ is unramified above $l$. Moreover, if $c_{l,i}$ is the eigenvalue of $T_{l,i}$ on $\pi_l^{K_l}$ for $i\in\{0,1,2\}$, then
    $$\imath^{-1}\mathrm{det}(T\mathrm{id}-\rho_\pi(\mathrm{Frob}_x))=T^4-c_{l,1}T^3+((l^3+l)c_{l,0}+lc_{l,2})T^2-l^3c_{l,0}c_{l,1}T+l^6c_{l,0}^2$$
    and $$\imath^{-1}\mathrm{sim}(\rho_\pi(\mathrm{Frob}_x))=l^3c_{l,0}.$$
    
    (Gal2) If prime $l=xx'$ splits in $F$, then $\imath^{-1}\mathrm{WD}(\rho_\pi\vert_{\mathrm{Gal}_{F_x}})$ coincides with $\iota_x^*(\pi_l\otimes\vert\mathrm{sim}\vert^\frac{3}{2})$ in the meaning of classical local Langlands correspondence.
    
    (Gal3) If $p=\wp\wp'$ splits in $F$, then the $p$-adic representation $\rho_\pi\vert_{\mathrm{Gal}_{F_\wp}}$ is de Rham. $\imath^{-1}\mathrm{WD}(\rho_\pi\vert_{\mathrm{Gal}_{F_\wp}})$ coincides with $\iota_\wp^*(\pi_p\otimes\vert\mathrm{sim}\vert^\frac{3}{2})$ in the meaning of classical local Langlands correspondence. Moreover, if $\pi_\infty\simeq W(\lambda)$, then $\rho_\pi\vert_{\mathrm{Gal}_{F_\wp}}$ has Hodge-Tate weights $h_1>h_2>h_3>h_4$.
\end{conj}
\begin{rem}\label{5.4}
    By the Cebotarev density theorem, and since the primes of $F$ splitting above $\mathbb{Q}$ have density $1$, (Gal1) alone determines the $4$-dimensional representation $\rho_\pi$ (after forgetting the $\mathrm{GSp}_4$-structure). Denote $\rho^c(g)=\rho(cgc)$ for a $\mathrm{Gal}_F$-representation $\rho$, $g\in\mathrm{Gal}_F$, where $c$ replaces the complex conjugation. It implies that after forgetting the $\mathrm{GSp}_4$-structure, we have
    \begin{equation}\label{U}
        \rho_\pi^c\simeq\rho_\pi^\vee(3),
    \end{equation}where $(3)$ means twisting by the cyclotomic character $3$ times.
    Moreover if $\rho_\pi$ is irreducible, it is uniquely determined by (Gal1) since the probable $\mathrm{GSp}_4$-structure is unique up to isomorphisms by Schur's lemma, hence \eqref{U} is also an isomorphism of representations with $\mathrm{GSp}_4$-structure.
\end{rem}

\subsection{$p$-adic theory: $p$-adic automorphic representations and eigenvarieties}
We maintain the notations in \S 5.1. Assume $p=\wp\wp'$ splits in $F$. Let $U^p=\prod_{l\neq p}U_l$ be a sufficiently small (that is, for some place $l$ its projection to $G(\mathbb Q_l)$ contains only one element of finite order, namely $1$) compact open subgroup of $G(\mathbb{A}_f^{\infty,p})$. We also assume that $U_l$ is hyperspecial if $l$ is inert in $F$. Let $S$ be the union of $\{p\}$ and the set of the primes $l$ such that $U_l$ is not hyperspecial.

Let $E$ be a finite extension of $\mathbb{Q}_p$ containing all square roots of integers and $\varpi_E$ be a uniformizer of $E$. For $k\in \mathbb Z_{\geq 1}$ and a compact open subgroup $U_p$ of $G(\mathbb{Z}_p)$, consider the $\mathcal{O}_E/\varpi_E^k$-module
$$S(U^pU_p,\mathcal{O}_E/\varpi_E^k)=\{f:G(\mathbb{Q}_p)\backslash G(\mathbb{A}^\infty)/U^pU_p\rightarrow\mathcal{O}_E/\varpi_E^k\}.$$
Put $$\widehat{S}(U^p,\mathcal{O}_E):=\lim\limits_{\substack{\longleftarrow\\k}}S(U^p,\mathcal{O}_E/\varpi_E^k):=\lim\limits_{\substack{\longleftarrow\\k}}\lim\limits_{\substack{\longrightarrow\\U_p}}S(U^pU_p,\mathcal{O}_E/\varpi_E^k)$$ and $\widehat{S}(U^p,E):=\widehat{S}(U^p,\mathcal{O}_E)\otimes_{\mathcal{O}_E}E$. $\widehat{S}(U^p,E)$ is the space of $p$-adic automorphisms of $G$ with tame level $U^p$ valued in $E$. $\widehat{S}(U^p,E)$ is an admissible unitary Banach representation of $G(\mathbb{Q}_p)\mathop{\simeq}\limits^{\iota_\wp}\mathrm{GSp}_4(\mathbb{Q}_p)$. Recall that $\widehat{S}(U^p,E)$ is equipped with a natural action of $\mathbb{T}(U^p)=E[T_{l,0}(U^p_l)^{\pm 1},T_{l,1}(U^p_l),T_{l,2}(U^p_l)\vert l\notin S\ splits\ in\ F]$. Here $U^p_l$ is the projection of $U^p$ on $G(\mathbb{Q}_l)$ and $T_{l,i}(U^p_l)$ is the element in $G(\mathbb{Q}_l)//U^p_l$ which is conjugated to $T_{l,i}\in G(\mathbb{Q}_l)//K_l$ for $i\in \{0,1,2\}$ and $l\notin S$ splitting in $F$. Let $F^S$ be the maximal algebraic extension of $F$ unramified outside the places dividing those in $S$, and $\mathrm{Gal}_{F,S}:=\mathrm{Gal}(F^S/F)$. Let $\rho:\mathrm{Gal}_{F,S}\rightarrow \mathrm{GSp}_4(E)$ be a continuous representation with $\mathrm{GSp}_4$-structure such that $\rho^c\simeq\rho^\vee(3)$ (see Rem. \ref{5.4} for $\rho^c$). Compared with (Gal1) in Conj. \ref{5.3}, to $\rho$, one can associate a maximal ideal $\mathfrak m_\rho$ of $\mathbb{T}(U^p)$ generated by $T_{l,0}(U^p_l)-l^{-3}\mathrm{sim}(\rho(\mathrm{Frob}_x)),T_{l,1}-\mathrm{tr}(\rho(\mathrm{Frob}_x)),T_{l,2}-l^{-1}a_2(\rho(\mathrm{Frob}_x))+(l^{-1}+l^{-3})\mathrm{sim}(\rho(\mathrm{Frob}_x))$, where $a_2$ denotes the coefficient of the quadratic term in the characteristic polynomial. Let $\omega_\rho$ denote the morphism $\mathbb{T}(U^p)\twoheadrightarrow\mathbb{T}(U^p)/\mathfrak{m}_\rho\simeq E$. Let $\boldsymbol{\alpha}=(\alpha_1,\alpha_2,\alpha_3,\alpha_4)\in (E^{\times})^4$ satisfy $\alpha_1\alpha_4=\alpha_2\alpha_3$ and $\alpha_i\alpha_j^{-1}\neq 1,p^{\pm 1}$ for $i\neq j$. Put $\phi:=\mathscr{L}(\mathrm{unr}(\boldsymbol{\alpha}))$. Let $\boldsymbol{h}=(h_1,h_2,h_3,h_4)\in \mathbb{Z}^4$ satisfy $h_1> h_2> h_3> h_4$ and $h_1+h_4=h_2+h_3$. Put $\delta_w:=w(\phi)\eta\mathscr L(z^{\boldsymbol h})p_1^{-2}p_2^{-1}$ fow $w\in W$. Assume $\rho$ is absolutely irreducible, $\widehat{S}(U^p,E)^\mathrm{lalg}[\mathfrak m_\rho]\neq 0$ and $D:=D_\mathrm{rig}(\rho\vert_{\mathrm{Gal}_{F_\wp}})\in\mathrm{GSp}_4$-$\Phi\Gamma_{\mathrm{nc}}(\boldsymbol{\alpha},\boldsymbol{h})$. The main goal of this chapter is to detect the structure of $\widehat{S}(U^p,E)^{\mathrm{an}}[\mathfrak m_\rho]$.

\begin{lem}\label{5.5}
    The action of $\mathbb{T}(U^p)\times \mathrm{GSp}_4(\mathbb{Q}_p)$ on $\widehat{S}(U^p,E)^\mathrm{lalg}$ is semi-simple.
\end{lem}
\proof
Assume $W$ is an irreducible algebraic representation of $\mathrm{GSp}_4$. Denote by $W_p$ the representation of $\mathrm{GSp}_4(\mathbb Q_p)$ over $E$ obtained by $W$. Denote by $W_\infty$ the pullback via $\iota_\infty$ of the complex representation of $\mathrm{GSp}_4(\mathbb C)$ obtained by $W$. Denote by $\hat W$ the contragredient representation of $W$. Then, by \cite[Thm. 2.2.16, Cor. 2.2.25, Prop. 3.2.2]{interpolation}, there is an $\mathbb{T}(U^p)\times \mathrm{GSp}_4(\mathbb{Q}_p)$-equivariant injection 
$$\widehat{S}(U^p,E)^{W-\mathrm{lalg}}\otimes_{E}\overline{\mathbb{Q}}_p\hookrightarrow \imath_*\mathcal{A}(G,\hat W_\infty)^{U^p}\otimes_E W_p.$$
Then $\widehat{S}(U^p,E)^{W-\mathrm{lalg}}$ is semi-simple by Lem. \ref{5.1}, for all irreducible representation $W$ of $\mathrm{GSp}_4$, which proves this lemma. 
\qed
\begin{lem}\label{5.6}
    There exists an integer $r\in \mathbb N$, such that $\widehat{S}(U^p,E)^\mathrm{lalg}[\mathfrak m_\rho]=\pi_\mathrm{alg}(\phi,\boldsymbol{h})^{\oplus r}$.
\end{lem}
\proof
    By Lem. \ref{5.5}, it's enough to prove that if a locally algebraic irreducible representation $\pi=\pi_{\mathrm{alg}}\otimes\pi_{\mathrm{sm}}$ contributes to $\widehat{S}(U^p,E)^\mathrm{lalg}[\mathfrak m_\rho]$, then both of the algebraic part $\pi_{\mathrm{alg}}$ and the smooth part $\pi_{\mathrm{sm}}$ of it are isomorphic to those of $\pi_\mathrm{alg}(\phi,\boldsymbol{h})$. By Rem. \ref{5.4}, and Conj. \ref{5.3}(Gal1), $\rho$ is just the Galois representation with $\mathrm{GSp}_4$-structure attached by Conj. \ref{5.3} to the automorphic representation attached to the classical $p$-adic automorphic representation generated by $\pi_\mathrm{alg}(\phi,\boldsymbol{h})\subset\widehat{S}(U^p,E)^\mathrm{lalg}[\mathfrak m_\rho]$ in the meaning of \cite[Prop. 3.2.4]{interpolation}. Then this lemma was proved by Conj. \ref{5.3}(Gal3).
\qed

To detect more information about $\widehat{S}(U^p,E)^{\mathrm{an}}[\mathfrak m_\rho]$, we need to introduce a technical geometric tool. Using Emerton’s method (ref. \cite[(2.3)]{interpolation}), one can construct an eigenvariety $\mathcal{E}(U^p)$ from $J_B(\widehat{S}(U^p,E)^\mathrm{an})$. The strong dual $J_B(\widehat{S}(U^p,E)^\mathrm{an})^\vee$ gives rise to a coherent sheaf $\mathcal{M}(U^p)$ over $\mathcal{E}(U^p)$. An $E$-point of $\mathcal{E}(U^p)$ can be parameterized by $(\delta,\omega)$ where $\delta$ is a continuous character of $T(\mathbb{Q}_p)$, and $\omega$ is a morphism of $E$-algebras $\mathbb T(U^p)\twoheadrightarrow E$ which corresponds to a maximal ideal $\mathfrak{m}_\omega$ of $\mathbb T(U^p)$. Moreover, $(\delta,\omega)\in \mathcal{E}(U^p)$ if and only if $\mathrm{Hom}_{T(\mathbb Q_p)}(\delta,J_B(\widehat{S}(U^p,E)^\mathrm{an})[\mathfrak{m}_\omega])\neq0$. Let $\widehat T$(resp. $\mathcal W$) be the rigid space which is parameterized by continuous characters of $T(\mathbb{Q}_p)$(resp. $T(\mathbb{Z}_p)$). Denote by $\kappa$(resp. $\kappa_0$) the natural morphism of rigid spaces from $\mathcal{E}(U^p)$ to $\widehat T$(resp. $\mathcal W$).
\begin{prop}\label{5.7}
    (i) $\mathcal{E}(U^p)$ is equidimensional of dimension $3$.
    
    (ii) There exists an admissible affinoid covering $(U_i)_{i\in I}$ of $\mathcal{E}(U^p)$, such that for every $i$, there exists an affinoid open subset $W_i$ of $\mathcal{W}$ such that $\kappa_0$ induces, when restricted to each irreducible component of $U_i$, a finite surjective morphism onto $W_i$. The image of any irreducible component of $\mathcal{E}(U^p)$ under $\kappa_0$ is Zariski open in $\mathcal W$.

    (iii) The coherent sheaf $\mathcal{M}(U^p)$ is Cohen-Macaulay over $\mathcal{E}(U^p)$.
\end{prop}
\proof
By the same argument of \cite[Lem. 6.1]{DingGL3}, for a compact open subgroup $H$ of $\mathrm{GSp}_4(\mathbb Z_p)$, we have $\widehat{S}(U^p,E)\vert_H\simeq \mathcal{C}(H,E)^{\oplus s}$ for some $s\geq1$ (where $\mathcal C(H,E)$ denotes the space of continuous functions on H). Then the proposition follows by the same arguments of \cite[Lem. 3.10, Prop. 3.11, Cor. 3.12, Cor 3.13]{breuil2017} and \cite[Lem. 3.8]{breuil2017smoothness}, applying \cite[\S5.2]{breuil2017} to $\widehat{S}(U^p,E)^{\mathrm{an}}$.
\qed

Assume $\lambda$ be a dominant weight of $T$. We call $(\delta,\omega)\in \mathcal{E}(U^p)(E)$ a classical point of weight $\lambda$ if $\delta=\delta_0\eta\delta_B\lambda$, such that $\delta_0$ is generic and unramified and 
\begin{equation}
    \mathrm{Hom}_{\mathrm{GSp}_4(\mathbb Q_p)}(L(\lambda)\otimes (\mathrm{Ind}_{\overline B(\mathbb{Q}_p)}^{\mathrm{GSp_4}(\mathbb{Q}_p)}\delta_0\eta)^\infty,\widehat{S}(U^p,E)^{\mathrm{an}}[\mathfrak m_\omega])\neq0.
\end{equation}
We call $(\delta,\omega)\in \mathcal{E}(U^p)(E)$ a very classical point of weight $\lambda$ if $\delta=\delta_0\eta\delta_B\lambda$, such that $\delta_0$ is generic and unramified and $L(\lambda)\otimes(\mathrm{Ind}_{\overline B(\mathbb{Q}_p)}^{\mathrm{GSp_4}(\mathbb{Q}_p)}\delta_0\eta)^\infty$ is the only irreducible subquotient of $(\mathrm{Ind}_{\overline B(\mathbb{Q}_p)}^{\mathrm{GSp_4}(\mathbb{Q}_p)}\delta_0\eta\lambda)^{\mathrm{an}}$ on which a $\mathrm{GSp}_4(\mathbb{Q}_p)$-invariant norm exists. (In the definition of (very) classical point, it is more natural to require $\delta_0$ to be smooth. However, for convenience, we require $\delta_0$ to be unramified, which is a little stronger.)
\begin{lem}\label{5.8}
    A very classical point is always a classical point.
\end{lem}
\proof
Assume $(\delta,\omega)$ is a very classical point. By \cite[Thm. 4.3]{Br2}, we have 
\begin{equation*}
\begin{aligned}
   &\mathrm{Hom}_{\mathrm{GSp}_4(\mathbb Q_p)}(\mathcal{F}_{\overline{B}}^{\mathrm{GSp}}((U(\mathfrak{g})\otimes_{U(\overline{\mathfrak{b}})}(-\lambda))^\vee,\delta_0\eta),\widehat{S}(U^p,E)^{\mathrm{an}}[\mathfrak m_\omega])\\&=\mathrm{Hom}_{T(\mathbb{Q}_p)}(\delta,J_B(\widehat{S}(U^p,E)^{\mathrm{an}})[\mathfrak m_\omega])\neq 0.
\end{aligned}    
\end{equation*}
By \cite{BGG} and \cite[Thm.]{OS}, $\mathcal{F}_{\overline{B}}^{\mathrm{GSp}}((U(\mathfrak{g})\otimes_{U(\overline{\mathfrak{b}})}(-\lambda))^\vee,\delta_0\eta)$ shares the same irreducible subquotients with $(\mathrm{Ind}_{\overline B(\mathbb{Q}_p)}^{\mathrm{GSp_4}(\mathbb{Q}_p)}\delta_0\eta\lambda)^{\mathrm{an}}$. Also, $L(\lambda)\otimes(\mathrm{Ind}_{\overline B(\mathbb{Q}_p)}^{\mathrm{GSp_4}(\mathbb{Q}_p)}\delta_0\eta)^\infty$ is the unique irreducible quotient of $\mathcal{F}_{\overline{B}}^{\mathrm{GSp}}((U(\mathfrak{g})\otimes_{U(\overline{\mathfrak{b}})}(-\lambda))^\vee,\delta_0\eta)$ and appears as its subquotient only in this way. Recall that $\widehat{S}(U^p,E)^{\mathrm{an}}$ is unitary, so every nonzero morphism from $\mathcal{F}_{\overline{B}}^{\mathrm{GSp}}((U(\mathfrak{g})\otimes_{U(\overline{\mathfrak{b}})}(-\lambda))^\vee,\delta_0\eta)$ to $\widehat{S}(U^p,E)^{\mathrm{an}}[\mathfrak m_\omega])$ factors through $(\mathrm{Ind}_{\overline B(\mathbb{Q}_p)}^{\mathrm{GSp_4}(\mathbb{Q}_p)}\delta_0\eta)^\infty\otimes L(\lambda)$. Hence $(\delta,\omega)$ is a classical point.
\qed

Consider a point $y=(\delta_y,\omega_y)\in \mathcal{E}(U^p)(E)$. Assume $\delta_y=\delta_{y,0}\eta\delta_B\lambda_y$, where $\delta_y$ is smooth and $\lambda_y$ is algebraic. Define $\boldsymbol{\alpha}_y=(\alpha_{y,1},\alpha_{y,2},\alpha_{y,3},\alpha_{y,4})\in E^4$ by 
$$\alpha_{y,1}=\delta_{y,0}(\mathrm{diag}(p,p,1,1)),\alpha_{y,2}=\delta_{y,0}(\mathrm{diag}(p,1,p,1)),$$
$$\alpha_{y,3}=\delta_{y,0}(\mathrm{diag}(1,p,1,p)),\alpha_{y,4}=\delta_{y,0}(\mathrm{diag}(1,1,p,p)),$$
i.e., $\delta_{y,0}=\mathscr{L}(\mathrm{unr}(\boldsymbol{\alpha}_y))$. Define $\boldsymbol{h}_y=(h_{y,1},h_{y,2},h_{y,3},h_{y,4})\in\mathbb Z^4$ by 
$$p^{h_{y,1}}=p^3\lambda_y(\mathrm{diag}(p,p,1,1))),p^{h_{y,2}}=p^2\lambda_y(\mathrm{diag}(p,1,p,1)),$$
$$p^{h_{y,3}}=p\lambda_y(\mathrm{diag}(1,p,1,p)),p^{h_{y,4}}=\lambda_y(\mathrm{diag}(1,1,p,p)),$$
i.e., $\lambda_y=\mathscr{L}(z^{\boldsymbol{h}_y})p_1^{-2}p_2^{-1}$. If $y$ is a classical point, we can attach a Galois representation $\rho_y:\mathrm{Gal}_{F,S}\rightarrow\mathrm{GSp_4}(E)$ to $y$ by \cite[Prop. 3.2.4]{interpolation} and Conj. \ref{5.3}. Up to isomorphism, $\rho_y$ is unique after forgetting the $\mathrm{GSp}_4$-structure. Moreover, if $\rho_y$ is absolutely irreducible, then it is unique up to isomorphism. By Conj. \ref{5.3}(Gal1), $\mathfrak{m}_{\rho_y}=\mathfrak{m}_{\omega_y}$. By Conj. \ref{5.3}(Gal3), $\rho_y\vert_{\mathrm{Gal}_{F_\wp}}$ is crystalline with Hodge-Tate weights $h_{y,1}>h_{y,2}>h_{y,3}>h_{y,4}$ and distinct $\varphi$-eigenvalues $\alpha_{y,1},\alpha_{y,2},\alpha_{y,3},\alpha_{y,4}$. We say a classical point $y$ is non-critical if the refinement $\boldsymbol{\alpha}_y$ of $\rho_y$ is non-critical. By the definition of $\mathcal{E}(U^p)$ and Lem. \ref{5.6}, $y_w:=(\delta_w\delta_B,\omega_\rho)\in\mathcal{E}(U^p)(E)$ for $w\in W$ is a classical point and $\rho_{y_w}=\rho$.
\begin{prop}\label{5.9}
    The classical points form a Zariski-dense and accumulation subset of $\mathcal{E}(U^p)$. More precisely, for every $E$-point $y_0$ such that $\kappa_0(y_0)$ is an algebraic character, and $X$ an irreducible component of $\mathcal{E}(U^p)$ containing $y_0$, there exists an affinoid neighborhood basis $\mathcal U$ of $y_0$ in $X$ such that those non-critical very classical points is Zariski-dense in $U$ for each $U\in\mathcal{U}$.
\end{prop}
\proof
Assume $U$ is an irreducible affinoid open neighborhood of $y_0$ such that $\kappa_0(U)$ is a affinoid open neighborhood of $\kappa_0(y_0)$ and $\kappa_0$ induces a finite surjective map from $U$ onto $\kappa_0(U)$. By Prop. \ref{5.7}(ii), we only need to prove the subset of non-critical very classical points is Zariski-dense in $U$.

For $i\in\{1,2,3,4\}$ and an $E$-point $y=\{\delta_y,\omega_y\}$ with locally algebraic $\delta_y$, $\alpha_{y,i}p^{h_{y,i}}$ equals to a product of a scalar and an interpolation of $\delta_y$ at a fixed element in $T(\mathbb{Q}_p)$. Hence for $i\in\{1,2,3,4\}$, $y\mapsto\alpha_{y,i}p^{h_{y,i}}$ extends to an analytic function $f_i$ on $\mathcal{E}(U^p)$. Fix $C>0$ such that $-C\leq\mathrm{val}_p(f_i(y))\leq C$ for $y\in U$, $i\in\{1,2,3,4\}$. Define
\begin{equation*}
    \begin{aligned}
        Z:=\{\lambda\in \kappa_0(U^p)(E)\vert\ \exists \boldsymbol{n}=(n_1,n_2&,n_3,n_4)\in\mathbb{Z}^4, s.t. n_1+n_4=n_2+n_3,\\ n_1>n_2>n_3>n_4, \min\limits_{i\in\{1,2,3\}}\{n_i-n_{i+1}\}>&20170901C+20260630, \lambda=\mathscr L(z^{\boldsymbol{n}})p_1^{-2}p_2^{-1}\}.
    \end{aligned}
\end{equation*}
$Z$ is obviously Zariski-dense in $\kappa_0(U)$, so $\kappa_0^{-1}(Z)\cap U$ is Zariski-dense in $U$. Now let's prove the points in $\kappa_0^{-1}(Z)\cap U$ are non-critical very classical points.

Assume $y=(\delta_y,\omega_y)\in\kappa_0^{-1}(Z)\cap U$, $\delta_y:=\delta_{0,y}\eta\delta_B\lambda_y$. Then we have $\delta_{0,y}$ is unramified and 
\begin{equation}\label{condition}
    \begin{cases}
        -C\leq h_{y,i}+\mathrm{val}_p(\alpha_{y,1})\leq C,\ \  i\in\{1,2,3,4\};\\
        h_{y,i}-h_{y,i+1}>20170901C + 20260630,\ \ i\in\{1,2,3\}.
    \end{cases}
\end{equation}
\eqref{condition} guarantees $\mathrm{val}(\alpha_i\alpha_j^{-1})>1$, for $1\leq i<j\leq4$, hence $\delta_{0,y}$ is generic.

By \cite[Thm. ]{OS} and \cite{BGG}, any irreducible subquotient of $(\mathrm{Ind}_{\overline B(\mathbb{Q}_p)}^{\mathrm{GSp_4}(\mathbb{Q}_p)}\delta_{0,y}\eta\lambda_y)^{\mathrm{an}}$ is in the form of $\mathcal{F}_{\overline B}^\mathrm{GSp}(\overline{L}(-w\cdot\lambda_y),\delta_{0,y}\eta)$, for some $w\in W$. If it admits a $\mathrm{GSp}_4(\mathbb{Q}_p)$-invariant norm, by \cite[Cor. 3.5]{Br1}, we have 
$(w\cdot\lambda_y)\delta_{o,y}\eta(t)\in \mathcal{O}_E$ for all $t\in T^+(\mathbb{Q}_p)$. Here $T^+(\mathbb{Q}_p):=\{\mathrm{diag}(t_1,t_2,t_3,t_4)\in T(\mathbb{Q}_p)\vert \mathrm{val}_p(t_1)\geq\mathrm{val}_p(t_2)\geq\mathrm{val}_p(t_3)\geq\mathrm{val}_p(t_4)\}$. By replacing $t$ by a couple of generators of $T^+(\mathbb{Q}_p)$ and direct calculation, one see this condition is equivalant to
\begin{equation}\label{condition2}
    \begin{cases}
        \sum\limits_{j=1}^ih_{y,\check w^{-1}(i)}+\mathrm{val}_p(\alpha_{y,i})\geq 0 ,\ i\in\{1,2,3\};\\
        \sum\limits_{j=1}^4h_{y,\check w^{-1}(i)}+\mathrm{val}_p(\alpha_{y,i})=0.
    \end{cases}
\end{equation}
Combining \eqref{condition} and \eqref{condition2}, $w$ must be trival. Hence $y$ is a very classical point. If the order of Hodge-Tate weights determined by the refinement $\boldsymbol{\alpha}_y$ is $(h_{y,\check w^{-1}(1)},$ $h_{y,\check w^{-1}(2)},$ $h_{y,\check w^{-1}(3)},$ $h_{y,\check w^{-1}(4)})$ (in the meaning of \cite[\S2.4.1]{Selmergp}). Use the weakly admissibility of $D_{\mathrm{cris}}(\rho_y)$ (\cite[\S5]{font}), we have
$$\begin{cases}
t_N(\bigoplus\limits_{j=1}^iD_{\mathrm{cris}}(\rho_y)^{\varphi=\alpha_{y,i}})\geq t_H(\bigoplus\limits_{j=1}^iD_{\mathrm{cris}}(\rho_y)^{\varphi=\alpha_{y,i}}),\ i\in\{1,2,3\};\\
t_N(D_{\mathrm{cris}}(\rho_y))= t_H(D_{\mathrm{cris}}(\rho_y)),
\end{cases}$$ which is still equivalent to \eqref{condition2}. This means $w=1$, i.e., the refinement $\boldsymbol{\alpha}_y$ of $\rho_y$ is non-critical.
\qed

\begin{prop}\label{5.10}
    $\mathcal{E}(U^p)$ is reduced.
\end{prop}
\proof
Assume $(\delta,\omega)\in\mathcal{E}(U^p)(E)$ is very classical. Then the fiber of $\kappa_*\mathcal{M}$ at the point parameterized by $\delta$ is
\begin{equation*}
    \begin{aligned}
        &\mathrm{Hom}_{T(\mathbb{Q}_p)}(\delta,J_B(\widehat{S}(U^p,E)^{\mathrm{an}}))\\
        =&\mathrm{Hom}_{\mathrm{GSp}_4(\mathbb Q_p)}(\mathcal{F}_{\overline{B}}^{\mathrm{GSp}}((U(\mathfrak{g})\otimes_{U(\overline{\mathfrak{b}})}(-\lambda))^\vee,\delta_0\eta),\widehat{S}(U^p,E)^{\mathrm{an}})\\
        =&\mathrm{Hom}_{\mathrm{GSp}_4(\mathbb Q_p)}(\mathcal{F}_{\overline{B}}^{\mathrm{GSp}}((U(\mathfrak{g})\otimes_{U(\overline{\mathfrak{b}})}(-\lambda))^\vee,\delta_0\eta),\widehat{S}(U^p,E)^{\mathrm{lalg}})\\
        =&\mathrm{Hom}_{T(\mathbb{Q}_p)}(\delta,J_B(\widehat{S}(U^p,E)^{\mathrm{alg}})).
    \end{aligned}
\end{equation*}
The first and third equation is because of \cite[Thm. 4.3]{Br2} and the second equation is because of the definition of very classical points. Hence by Lem. \ref{5.5}, this fiber is semi-simple as a $\mathbb{T}(U^p)$-module. Then this proposition follows by a similar argument of \cite[Prop. 3.9]{che2005}(see also \cite[Cor. 3.20]{breuil2017}), using this fact and Prop. \ref{5.9}. 
\qed

\begin{prop}\label{5.11}
    (i) There exists a (unique) continuous pseudo-character $\mathcal{T}:\mathrm{Gal}_{F,S}\allowbreak\rightarrow\mathcal{O}(\mathcal{E}(U^p))$ such that $\mathcal T_{y}=\mathrm{tr}(\rho_y)$ for all classical points $y\in \mathcal{E}(U^p)(E)$, where $\mathcal{T}_y$ replaces for the evaluation of $\mathcal{T}$ at $y$.
    
    (ii) Let $y_0\in \mathcal{E}(U^p)(E)$ satisfy that $\kappa_0(y_0)$ is an algebraic character and $\mathcal{T}_{y_0}$ is absolutely irreducible. Then there exists an affinoid open neighborhood $X$ of $y_0$ in which classical points form a Zariski-dense subset and a continuous representation with $\mathrm{GSp}_4$-structure $\rho_X:\mathrm{Gal}_{F,S}\rightarrow\mathrm{GSp}_4(\mathcal{O}(X))$ such that the evaluation of $\rho_X$ at $y$ is isomorphic to $\rho_y$ for each classical points $y\in X$ and after forgetting the $\mathrm{GSp}_4$-structure, $\rho_X^c\simeq\rho_X^\vee(3)$.
\end{prop}
\proof
Notice that those Hecke operators $T_{l,i}(U^p_l)$ preserves $\widehat{S}(U^p,\mathcal{O}_E)$, so the image of them in $\mathcal{O}(\mathcal{E}(U^p))$ is contained in $\mathcal{O}(\mathcal E(U^p))^0$, which is a compact subring by \cite[Lem. 7.2.11]{Selmergp}. Then the existence of $\mathcal{T}$ follows by \cite[Prop. 7.1.1]{che2004}. Also, there exists a continuous character $\mathcal{S}:\mathrm{Gal}_{F,S}\rightarrow \mathcal{O}(\mathcal{E}(U^p))^\times$ whose evaluation $\mathcal S_y=\mathrm{sim}(\rho_y)$ at all classical points $y\in \mathcal{E}(U^p)(E)$. For $g\in \mathrm{Gal}_{F,S}$, we have $\mathcal{T}(g^{-1})S(g)=\mathcal{T}(g)$ because this equation holds at all classical points. 

By a completely same argument of \cite[Lem. 5.5]{galoisfam}, there exists an affinoid open neighborhood $X$ of $y_0$ in which classical points form a Zariski-dense subset and a continuous representation $\rho_X: \mathrm{Gal}_{F,S}\rightarrow\mathrm{GL}(V_X)$, where $V_X$ is a free $\mathcal{O}(X)$-module with rank $4$, such that $\mathrm{tr}(\rho_X)=\mathcal{T}_X$ and the evaluation of $\rho_X$ is absolutely irreducible at each closed point. Then the two representation $V_X^\vee\otimes \mathcal{S}_X$ and $V_X$ share a same trace map. By \cite[Thm. 1]{pseudo} and the fact that $\mathcal{O}_{X,{y_0}}$ is henselian, after shrinking $X$ if necessary, there exists an isomorphism of $\mathrm{Gal}_{F,S}$-representations $s_X: V_X\rightarrow V_X^\vee\otimes \mathcal{S}_X$. For a classical point $y\in X(E)$, the evaluation of $s_X$ at $y$ must coinside with the $\mathrm{GSp}_4$-structure of $\rho_y$ up to isomorphism because $\rho_y$ is irreducible. Then the bilinear form $r_X: V_X\times V_X\rightarrow\mathcal{O}(X)$ defined by $r_X(v,w)=(s_X(w))(v)$ is alternate because its evaluation at each classical point is so. Shrinking $X$ if needed, one may assume under a suitable ordered basis, $r_X$ is represented by the matrix $J=\begin{bmatrix}&&&1\\&&1&\\&-1&&\\-1&&&\end{bmatrix}$. Then under this ordered basis, $\rho_X$ takes value in $\mathrm{GSp}_4(\mathcal{O}(X))$. Similarly, for $g\in\mathrm{Gal}_{F,S}$, $\mathcal{T}_X(cgc)=\mathcal{T}_X(g^{-1})\chi_{\mathrm{cyc}}^3(g)$ because this equation holds at all classical points. Hence the two representation $\rho_X^c$ and $\rho_X^\vee(3)$ shares a same trace map and then they are isomorphic after shrinking $X$ if necessary. 
\qed

\begin{prop}\label{5.12}
   (i) $(\delta, \omega_\rho)\in\mathcal{E}(U^p)(E)$ iff $\delta=\delta_w\delta_B$ for some $w\in W$.

   (ii) If $C$ is an irreducible subquotient of $\mathrm{PS}(w(\phi),\boldsymbol{h})$, $w\in W$. Then $C$ can be embedded into $\widehat{S}(U^p,E)^\mathrm{an}[\mathfrak{m}_\rho]$ iff $C\simeq\pi_\mathrm{alg}(\phi,\boldsymbol{h})$.
\end{prop}
\proof
(i): Define following analytic functions on $\mathcal{E}(U^p)$:
$$g_1(y):=\mathrm{wt}((\kappa_0(\delta_y))\circ\mathrm{diag}(z,z,1,1))+3:=\frac{\partial}{\partial z}_{\vert z=1}((\kappa_0(\delta_y))\circ\mathrm{diag}(z,z,1,1))+3,$$
$$g_2(y):=\mathrm{wt}((\kappa_0(\delta_y))\circ\mathrm{diag}(z,1,z,1))+2:=\frac{\partial}{\partial z}_{\vert z=1}((\kappa_0(\delta_y))\circ\mathrm{diag}(z,1,z,1))+2,$$
$$g_3(y):=\mathrm{wt}((\kappa_0(\delta_y))\circ\mathrm{diag}(1,z,1,z))+1:=\frac{\partial}{\partial z}_{\vert z=1}((\kappa_0(\delta_y))\circ\mathrm{diag}(1,z,1,z))+1,$$
$$g_4(y):=\mathrm{wt}((\kappa_0(\delta_y))\circ\mathrm{diag}(1,1,z,z)):=\frac{\partial}{\partial z}_{\vert z=1}((\kappa_0(\delta_y))\circ\mathrm{diag}(1,1,z,z)),$$
for $y=(\delta_y,\omega_y)\in \mathcal{E}(U^p)(E)$. $g_1g_4=g_2g_3$. $\kappa(y)$ is locally algebraic iff $g_1(y),g_2(y),g_3(y),\allowbreak g_4(y)$ are all integers. If this is the case, then $g_i(y)=h_{y,i}$ for $i\in\{1,2,3,4\}$.

Assume $y_0=(\delta_{y_0},\omega_\rho)\in\mathcal{E}(U^p)(E)$. Notice that $\mathcal T_y$ is determined only by $\omega_y$ for $y=(\delta_y,\omega_y)$ by  definition, hence $\mathcal T_{y_0}=\mathrm{tr}\rho$. By \cite[7.5.12]{Selmergp}, there is a unique $u\in W$ such that $g_i(y_0)=h_{\check{u}^{-1}(i)}$ for $i\in\{1,2,3,4\}$. This means $\delta_{y_0}$ is locally algebraic and $h_{y_0,i}=h_{\check{u}^{-1}(i)}$ for $i\in\{1,2,3,4\}$. Let $\delta_{y_0}=\delta_{0,y_0}\eta\delta_B\lambda_{y_0}$, where $\delta_{0,y_0}$ is smooth and $\lambda_{y_0}=\mathscr{L}(z^{\boldsymbol{h}_{y_0}})p_1^{-2}p_2^{-1}$ is algebraic. For simplicity, put $\lambda_\rho:=\mathscr{L}(z^{\boldsymbol h})p_1^{-2}p_2^{-1}$, then $\lambda_{y_0}=u\cdot\lambda_\rho$. 

Let $\delta_X:T(\mathbb{Q}_p)\rightarrow \mathcal{O}(X)^\times$ be the natural character associated to $\kappa$. By \cite[Prop 2.4.1]{Selmergp}, for a classical point $y$, $\rho_y$ has a triangulation parameterized by the cocharacter $\mathrm{unr}(\boldsymbol{\alpha}_y)z^{\boldsymbol{h}_y}:\mathbb{Q}_p^\times\rightarrow T(E)$, which is the evaluation of $\mathscr{L}^{-1}(\delta_X\eta^{-1}\delta_B^{-1}p_1^2p_2):\mathbb{Q}_p^\times\rightarrow T(\mathcal{O}(\mathcal{E}(U^p))$. Take $X$ and $\rho_X$ as Prop. \ref{5.11}(ii). Using \cite[Thm. 6.3.13]{Xiao} on $\rho_X$, we deduce that $\rho$ has a triangulation parameterized by $\mathscr{L}^{-1}(\delta_{0,y_0})z^{\boldsymbol{h}_{y_0}+\boldsymbol{h}'}$, for some $\boldsymbol{h}'=(h'_1,h'_2,h'_3,h'_4)\in \mathbb Z^4$. Comparing the smooth part of the parameter of the triangulaion with \cite[Prop. 2.4.1]{Selmergp}, we have $\delta_{0,y_0}=w(\phi)=\mathscr{L}(\check w\mathrm{unr}(\boldsymbol{\alpha}))$ for some $w\in W$. Similarly, comparing the algebraic part, we have $\boldsymbol{h}_{y_0}+\boldsymbol{h}'=\boldsymbol{h}$. Moreover, for $i\in\{1,2,3,4\}$, applying \cite[Cor. 6.3.10(2')]{Xiao} on $\wedge^i\rho_X$, by a similar argument to \cite[Ex. 6.3.14]{Xiao}, we have $\sum\limits_{j=1}^ih'_j\geq 0$. Then it must be $\boldsymbol{h}_{y_0}=\boldsymbol{h}$, $\boldsymbol{h}'=0$, $u=1$ and $\lambda_{y_0}=\lambda_\rho$. Then $\delta_{y_0}=\delta_w\delta_B$. (i) was proved. 

(ii): By \cite[Thm. ]{OS}, $C\simeq\mathscr C(w,u):=\mathcal{F}_{\overline B}^\mathrm{GSp}(\overline{L}(-u\cdot\lambda),w(\phi)\eta)$ for some $u,w\in W$. Then by \cite[Thm. 4.3]{Br2}, $(w(\phi)\eta\delta_B(u\cdot\lambda_\rho),\omega_\rho)\in \mathcal{E}(U^p)$. Hence by (1), $u=1$ and then (ii) was proved.
\qed

Let ($\mathrm{GSp}_4$-)$X_\rho$ denotes the functor of deformations (with $\mathrm{GSp}_4$-structure) over local Artinian $E$-algebras of $\rho$. Replacing $(\varphi,\Gamma)$-module cohomology with group cohomology in Prop. \ref{3.3} and Lem. \ref{3.4}, we have that $\mathrm{GSp}_4$-$X_\rho$ is a subfunctor of $X_\rho$. Hence if a deformation with $\mathrm{GSp}_4$-structure $\widetilde\rho$ satisfies $\widetilde\rho\simeq(\widetilde\rho^c)^\vee(3)$ after forgetting the $\mathrm{GSp}$-4 structure, then $\widetilde\rho\simeq(\widetilde\rho^c)^\vee(3)$ also holds without forgetting the $\mathrm{GSp}$-4 structure. We denote by $\mathrm{GSp}_4$-$X_{\rho,U}$ the subfunctor of $\mathrm{GSp}_4$-$X_\rho$ of deformations with $\mathrm{GSp}_4$-structure satisfying the condition above. Similarly to \S 3.2.2, we denote the tangent space of $\mathrm{GSp}_4$-$X_{\rho}$ (resp. $\mathrm{GSp}_4$-$X_{\rho,U}$) by $\mathrm{Ext}^G(\rho,\rho)$ (resp. $\mathrm{Ext}^G_U(\rho,\rho)$). Denote by $\mathrm{Ext}^G_U(D,D)$ the image of $\mathrm{Ext}^G_U(\rho,\rho)$ via $\mathrm{Ext}^G(\rho,\rho)\rightarrow\mathrm{Ext}^G(\rho\vert_{\mathrm{Gal}_{F_\wp}},\rho\vert_{\mathrm{Gal}_{F_\wp}})\simeq \mathrm{Ext}^G(D,D)$.

We make the following vanishing hypothesis:
\begin{hyp}\label{5.13}
    Suppose $\mathrm{Ext}^G_U(\rho,\rho)\rightarrow\mathrm{Ext}^G(D,D)/\mathrm{Ext}^G_g(D,D)$ is injective. In particular $\mathrm{Ext}^G_U(D,D)\cap\mathrm{Ext}^G_g(D,D)=0$.
\end{hyp}
\begin{rem}\label{5.14}
    Philosophically speaking, if $\rho$ is also a Galois representation attached to an automorphic representation of the definite unitary group over $\mathbb{Q}$ attached to $F/\mathbb{Q}$, then by \cite[Rem. 4.16]{Ding}, Hyp. \ref{5.13} is known to hold in many cases.
\end{rem}

Let $R_{\rho,U}$ be the universal deformation ring of $\mathrm{GSp}_4$-$X_\rho$. Let $\mathfrak{a}_D\supset\mathfrak{m}_{R_D}^2$ be the ideal associated to $\mathrm{Ext}^G_U(D,D)$. By Hyp. \ref{5.13}, $\mathrm{Ext}^G_U(\rho,\rho)\simeq\mathrm{Ext}^G_U(D,D)$. The natural morphism $R_D\rightarrow R_{\rho,U}$ induces hence an isomorphism (of local Artinian $E$-algebras) $A_{D,U}:=R_D/\mathfrak{a}_D\xrightarrow{\sim}R_{\rho,U}/\mathfrak{m}_{R_{\rho,U}}^2$. Let $\widetilde{\rho}_U$ be the universal deformation of $\rho$ over $R_{\rho,U}/\mathfrak{m}_{R_{\rho,U}}^2$. We have a natural morphism $\mathbb{T}(U^p)\rightarrow R_{\rho,U}/\mathfrak{m}_{R_{\rho,U}}^2$, sending $T_{l,0}(U^p_l)$ to $l^{-3}\mathrm{sim}(\widetilde\rho_U(\mathrm{Frob_x}))$, $T_{l,1}(U^p_l)$ to $\mathrm{tr}(\widetilde\rho_U(\mathrm{Frob_x}))$, $T_{l,2}(U^p_l)$ to $l^{-1}a_2(\widetilde\rho_U(\mathrm{Frob_x}))-(l^{-1}+l^{-3})\mathrm{sim}(\widetilde\rho_U(\mathrm{Frob_x}))$, where $a_2$ denotes the coefficient of the quadratic term in the characteristic polynomial. Let $\mathfrak{a}_T$ be its kernel. The induced morphism $\mathbb{T}(U^p)/\mathfrak{a}_T\rightarrow R_{\rho,U}/\mathfrak{m}_{R_{\rho,U}}^2$ is an isomorphism. Indeed, it suffices to show the morphism sends $\mathfrak{m}_\rho$ to onto $\mathfrak{m}_{R_{\rho,U}}/\mathfrak{m}_{R_{\rho,U}}^2$. Consider the universal representation $\widetilde \rho$ over $R_{\rho,U}/\mathfrak{m}_\rho$. By the definition of the morphism, we deduce the pseudo-character $\mathrm{tr}(\widetilde\rho)$ takes values in $E$. This implies $\mathrm{tr}(\widetilde\rho)$ is a trivial deformation of $\mathrm{tr}(\rho)$. As $\rho$ is absolutely irreducible, deforming $\rho$ is equivalent to deforming $\mathrm{tr}(\rho)$ by \cite[Thm. 1]{pseudo}. We deduce $\widetilde\rho$ is a trivial deformation of $\rho$. Hence  $\mathfrak{m}_\rho(R_{\rho,U}/\mathfrak{m}_{R_{\rho,U}}^2)=\mathfrak{m}_{R_{\rho,U}}/\mathfrak{m}_{R_{\rho,U}}^2$.

\begin{prop}\label{5.15}
    Assume Hyp. \ref{5.13}, then $\mathcal{E}(U^p)$ is smooth at $y_w=(\delta_w\delta_B,\omega_\rho)$ for all $w\in W$. Moreover, $\mathrm{Ext}^G_U(D,D)$ satisfies Hyp. \ref{4.22}.
\end{prop}
\proof
Take an affinoid neighborhood $X$ of $y_w$ and $\rho_X:\mathrm{Gal}_{F,S}\rightarrow\mathrm{GSp}_4(\mathcal{O}(X))$ as Prop. \ref{5.11}(ii). Shrinking $X$ if needed, we may assume $y_{w'}\notin X$ for $w'\neq w$. Let $\mathscr{L}^{-1}(\delta_X\eta^{-1}\delta_B^{-1}p_1^2p_2)=\mathrm{diag}(\gamma_1,\gamma_2,\gamma_3,\gamma_4):\mathbb{Q}_p^\times\rightarrow T(\mathcal{O}(X))$, where $\delta_X:T(\mathbb{Q}_p)\rightarrow \mathcal{O}(X)^\times$ is the natural character associated to $\kappa$. Let $\mathcal{R}_{X}$ be the relative Robba ring over $\mathcal{O}(X)$ defined in \cite[Def. 2.2.2]{Xiao}. By \cite[Thm. 5.3]{globaltri} (and an easy induction argument), after shrinking $X$ if needed, we have,

(*): The ($\varphi$,$\Gamma$)-module $D_{\mathrm{rig}}(\rho_X\vert_{\mathrm{Gal}_{F_\wp}})$ is isomorphic to a successive extension of the rank one ($\varphi$,$\Gamma$)-modules (\cite[Def. 6.2.1]{Xiao}) $\mathcal{R}_X(\gamma_1)$, $\mathcal{R}_X(\gamma_2)$, $\mathcal{R}_X(\gamma_3)$, $\mathcal{R}_X(\gamma_4)$. 

Let $T_{y_w}$ be the tangent space of $\mathcal{E}(U^p)$ at the point $y_w$. By Prop. \ref{5.11}(ii), we then have a map $f:T_{y_w}\rightarrow\mathrm{Ext}^G_U(\rho,\rho)$, sending $t:\mathrm{Spec}E[\epsilon]/(\epsilon^2)\rightarrow \mathcal{E}(U^p)$ to $t^*\rho_X$. We are going to prove $f$ is injective. 

Denote by $f_{1}:T_{y_w}\rightarrow\mathrm{Ext}^1(\delta_w\delta_B,\delta_w\delta_B)$ the tangent map of $\kappa:\mathcal{E}(U^p)\rightarrow\widehat{T}$ at $y_w$. Denote by $f_2:T_{y_w}\rightarrow (\mathfrak{m}_\rho/\mathfrak{m}_\rho^2)^\vee$ the tangent map of $\mathbb{T}(U^p)\rightarrow \mathcal{O}(X)$ at $y_w$. Let $\nu\in\mathrm{ker}f$. Then $\nu^*\rho_X$ is the trivial deformation of $\rho$. This implies the composition $\mathbb{T}(U^p)\rightarrow\mathcal{O}(X)\xrightarrow{\nu}E[\epsilon]/(\epsilon^2)$ factors through $\mathfrak{m}_\rho$, i.e. $f_2(\nu)=\rho$. By (*), $\mathscr{L}^{-1}(f_1(\nu)\eta^{-1}\delta_B^{-1}p_1^2p_2)$ is the parameter of the triangulation of the trivial deformation $f(\nu)$. Hence $\mathscr{L}^{-1}(f_1(\nu)\eta^{-1}\delta_B^{-1}p_1^2p_2)$ is the trivial deformation of $\check{w}(\mathrm{unr}(\boldsymbol{\alpha}))z^{\boldsymbol{h}}=\mathscr{L}^{-1}(\delta_w\delta_B\eta^{-1}\delta_B^{-1}p_1^2p_2)$, i.e. $f_1(v)$ is the trivial deformation of $\delta_w\delta_B$. But by the construction of $\mathcal{E}(U^p)$, $\mathbb{T}(U^p)\otimes E[T(\mathbb{Q}_p)]$ is dense in $\mathcal{O}(X)$, hence the map $T_{y_w}\xrightarrow{(f_1,f_2)}\mathrm{Ext}^1(\delta_w\delta_B,\delta_w\delta_B)\times (\mathfrak{m}_\rho/\mathfrak{m}_\rho^2)^\vee$ is injective. We deduce $\nu$ is zero. Hence $f$ is injective.

By (*), the image of the composition $T_{y_w}\rightarrow\mathrm{Ext}^G_U(\rho,\rho)\rightarrow \mathrm{Ext}^G_U(D,D)$ is in $\mathrm{Ext}^G_{U,\check{w}}(D,D)$. Together with Hyp. \ref{5.13}, we have $\mathrm{dim}_ET_{y_w}\leq\mathrm{dim}_E\mathrm{Ext}^G_{U,\check{w}}(D,D)\leq\mathrm{dim}_E\mathrm{Ext}^G_{\check{w}}(D,D)/\mathrm{Ext}^G_g(D,D)=3$. As $\mathrm{dim}_E\mathcal{E}(U^p)=3$, we see $y_w$ is a smooth point and $\mathrm{dim}_E\mathrm{Ext}^G_{U,\check{w}}(D,D)=3$. This finishes the proof.
\qed

By Lem. \ref{5.6}, Prop. \ref{5.7}(iii) and Prop. \ref{5.15}, $\mathcal{M}(U^p)$ is locally free of rank $r$ at each $y_w$ for $w\in W$. Let $X$ be a sufficiently small smooth neighbourhood of $y_w$, and $\mathfrak{m}_{y_w}\subset\mathcal{O}(X)$ be the maximal ideal associated to $y_w$. We use the notation of \S4.3: $A_{D,U}$, $A_{D,U,w}$, $A_{U,w}$, etc. By the proof of Prop. \ref{5.15}, the composition $A_{D,U}=R_D/\mathfrak{a}_D\simeq R_{\rho,U}/\mathfrak{m}_{R_{\rho,U}}^2\mathop{\twoheadrightarrow}\limits^{f}\mathcal{O}(X)/\mathfrak{m}_{y_w}^2$ factors through a unique isomorphism $A_{D,U,\check{w}}\xrightarrow{\sim}\mathcal{O}(X)/\mathfrak{m}_{y_w}^2$. The map $A_{U,\check{w}}\xrightarrow{\sim} A_{D,U,\check{w}}\xrightarrow{\sim}\mathcal{O}(X)/\mathfrak{m}_{y_w}^2$ coincides with the composition $T_{y_w}\xrightarrow[\sim]{f_1}\mathrm{Ext}^1_U(\delta_w\delta_B,\delta_w\delta_B)\xrightarrow{\mathscr{L}^{-1}}\mathfrak{m}_{A_{U,\check{w}}}^\vee$. The map $\mathbb{T}(U^p)/\mathfrak{a}_T\xrightarrow{\sim} R_{\rho,U}/\mathfrak{m}_{R_{\rho,U}}^2\mathop{\twoheadrightarrow}\limits^{f}\mathcal{O}(X)/\mathfrak{m}_{y_w}^2$ coincides with the one induced by $\mathbb{T}(U^p)\rightarrow\mathcal{O}(X)$. We deduce a $T(\mathbb{Q}_p)\times A_{D,U,\check{w}}$-equivariant injection
\begin{equation}\label{final1}
    (\widetilde{\delta}_{U,w}^{\textrm{univ}}\delta_B)^{\oplus r}\simeq (\mathcal{M}(U^p)/\mathfrak{m}_{y_w}^2)^\vee\hookrightarrow J_B(\widehat{S}(U^p,E)^\mathrm{an})[\mathfrak{a}_T][\mathcal{I}_{\check{w}}]\{T(\mathbb{Q}_p)=\delta_w\delta_B\},
\end{equation}
where the $A_{D,U,\check{w}}$-action on the left hand side is given as in the discussion below Cor. \ref{4.25} and it acts on the right hand side via $A_{D,U,\check{w}}\simeq A_{D,U}/\mathcal{I}_{\check{w}}\twoheadleftarrow A_{D,U}\simeq R_{\rho,U}/\mathfrak{m}_{R_{\rho,U}}^2\simeq \mathbb{T}(U^p)/\mathfrak{a}_T$. Note as in the discussion below Cor. \ref{4.25}, the action of $A_{D,U,\check{w}}$ and $T(\mathbb{Q}_p)$ determine each other.

\begin{lem}\label{5.16}
    The map \eqref{final1} is balanced in the sense of \cite[Def. 0.8]{Jacq2}, hence (by \cite[Thm. 0.13]{Jacq2}) induces a $\mathrm{GSp}_4(\mathbb{Q}_p)\times A_{D,U,\check{w}}$-equivariant injection
    \begin{equation}
        \iota_w: (I_{\overline B(\mathbb{Q}_p)}^{\mathrm{GSp}_4(\mathbb{Q}_p)}\widetilde{\delta}_{U,w}^{\textrm{univ}})^{\oplus r}\hookrightarrow \widehat{S}(U^p,E)^\mathrm{an}[\mathfrak{a}_T][\mathcal{I}_{\check{w}}].
    \end{equation}
\end{lem}
\proof
Similar to \cite[Lem. 4.11]{DingGLn}, using Prop. \ref{5.12}.
\qed

Let $\widetilde\pi$ be the subrepresentation of $\widehat{S}(U^p,E)^\mathrm{an}[\mathfrak{a}_T]$ generated by $\mathrm{Im}\iota_w$ for all $w\in W$. Note $\widetilde\pi$ inherits from $\widehat{S}(U^p,E)^\mathrm{an}[\mathfrak{a}_T]$ an action of $A_{D,U}(\simeq\mathbb{T}(U^p)/\mathfrak{a}_T)$.

\begin{thm}\label{5.17}
    Suppose Hyp. \ref{5.13}. We have a $\mathrm{GSp}_4(\mathbb{Q}_p)\times A_{D,U}$-equivariant isomorphism $\widetilde\pi\simeq(\pi_1(\phi,\boldsymbol{h})^{\mathrm{univ}}_U)^{\oplus r}$. Consequently, we have $\pi_{\mathrm{min}}(D)^{\oplus r}\hookrightarrow \widehat{S}(U^p,E)^\mathrm{an}[\mathfrak{m}_\rho]$.
\end{thm}
\proof
By Lem. \ref{4.4}(i), Lem. \ref{5.6} and Prop. \ref{5.12}(ii), it is not difficult to see $\widetilde\pi\simeq(\pi_1(\phi,\boldsymbol{h})^{\mathrm{univ}}_U)^{\oplus r}$ as $\mathrm{GSp}_4(\mathbb{Q}_p)$ representations. It also preserves the $A_{D,U}$-action because its restriction at each $\mathrm{PS}_1(w(\phi),\boldsymbol{h})$ does so by the discussion below Cor. \ref{4.25}. The second part follows by Cor. \ref{4.26}.
\qed

Using Lem. \ref{4.4}(i), Lem. \ref{5.6} and Prop. \ref{5.12}(ii) again, it is not easy to see that for any subrepresentation $V$ of $\widehat{S}(U^p,E)^\mathrm{an}[\mathfrak{a}_T]$, if $V$ is isomorphic to an extension of $\pi_{\mathrm{alg}}(\phi,\boldsymbol{h})$ by $\pi_1(\phi,\boldsymbol{h})$, then $V\subset\widetilde\pi$. Hence we have:
\begin{cor}\label{5.18}
The representation $\pi_{\mathrm{min}}(D)^{\oplus r}$ is the maximal subrepresentation of $\widehat{S}(U^p,E)^\mathrm{an}[\mathfrak{m}_\rho]$, which is generated by extensions of $\pi_{\mathrm{alg}}(\phi,\boldsymbol{h})$ by $\pi_1(\phi,\boldsymbol{h})$.
\end{cor}
\qed
\begin{cor}\label{5.19}
    For $D'\in\mathrm{GSp}_4$-$\Phi\Gamma(\boldsymbol{\alpha},\boldsymbol{h})$, $\pi_{\mathrm{min}}(D')\hookrightarrow \widehat{S}(U^p,E)^\mathrm{an}[\mathfrak{m}_\rho]$ iff $D'\simeq D$. Hence $\widehat{S}(U^p,E)^\mathrm{an}[\mathfrak{m}_\rho]$ determines $\rho\vert_{\mathrm{Gal}_{F_\wp}}$.
\end{cor}
\qed
\newpage


\bibliographystyle{IEEEtran}
\bibliography{main.bib}


\end{document}